\documentclass[12pt]{article}

\usepackage[flushleft]{threeparttable}
\usepackage{multirow}
\usepackage{shortcuts}
\usepackage{stmaryrd}
\usepackage{lipsum}
\usepackage{graphicx}
\graphicspath{{./figures/}} 
\usepackage{hyperref}
\usepackage{subfigure}

\hypersetup{
    colorlinks,
    citecolor=black,
    filecolor=black,
    linkcolor=black,
    urlcolor=black
}
\newcommand{\pvphi}{\pmb{\vphi}}
\newcommand{\brbz}{\breve{\bz}}

\newcommand{\cE}{\mathcal{E}}

\newcommand{\bxi}{\pmb{\xi}}

\title{On structure-preserving discontinuous Galerkin methods for Hamiltonian partial differential equations: Energy conservation and multi-symplecticity}
\author{
	Zheng Sun\footnote{Department of Mathematics, The Ohio State University,
		Columbus, OH 43210, USA. E-mail: sun.2516@osu.edu.} \and Yulong Xing\footnote{Department of Mathematics, The Ohio State University,
		Columbus, OH 43210, USA. E-mail: xing.205@osu.edu. The work of this author is partially supported by the NSF grant DMS-1753581.}}
\date{}

\begin{document}
\maketitle 

\textbf{Abstract:} In this paper, we present and study discontinuous Galerkin (DG) methods for one-dimensional multi-symplectic Hamiltonian partial differential equations. We particularly focus on semi-discrete schemes with spatial discretization only, and show that the proposed DG methods can simultaneously preserve the multi-symplectic structure and energy conservation with a general class of numerical fluxes, which includes the well-known central and alternating fluxes. Applications to the wave equation, the Benjamin--Bona--Mahony equation, the Camassa--Holm equation, the Korteweg--de Vries equation and the nonlinear Schr\"odinger equation are discussed. Some numerical results are provided to demonstrate the accuracy and long time behavior of the proposed methods. Numerically, we observe that certain choices of numerical fluxes in the discussed class may help achieve better accuracy compared with the commonly used ones including the central fluxes.


\section{Introduction}
\setcounter{equation}{0}
In this paper, we present and study discontinuous Galerkin (DG) finite element methods for multi-symplectic Hamiltonian partial differential equations (HPDEs) in one dimension
\begin{equation}\label{eq-csv}
M \bz_t + K \bz_x = \nabla_\bz S(\bz).
\end{equation}
Here $\bz =\bz(x,t): \Omega\times [0,+\infty)\to \mathbb{R}^m$ is a vector-valued function, $M$ and $K$ are $m\times m$ real anti-symmetric matrices, and $S: \mathbb{R}^m \to \mathbb{R}$ is a smooth function.  For simplicity, periodic or compactly supported boundary conditions are considered. 

Symplectic integrator for Hamiltonian ordinary differential equations (HODEs)
\begin{equation}\label{eq-hodes}
M\bz_t = \nabla S(\bz), \qquad \bz = \bz(t),
\end{equation} 
is a well studied subject \cite{hairer2006geometric}, which has been widely used in applications such as rigid  body and  molecular dynamics.   The formulation \eqref{eq-csv} was proposed by Bridges and Reich 
for generalizing similar concepts to partial differential equations \cite{bridges2001multi}.
It applies to equations in various of fields such as classical mechanics, quantum physics and hydrodynamics, with examples including but not limited to the Sine-Gordon equation, the Hamiltonian wave equation, the Korteweg--de Vries (KdV) equation, the Camassa--Holm (CH) equation, the Benjamin--Bona--Mahony (BBM), the nonlinear Schr\"odinger (NLS) equation, Maxwell's equations, and the Dirac equation, etc. As analogues of Hamiltonian-preserving and symplecticity-preserving properties of the HODEs \eqref{eq-hodes}, it is well-known that the multi-symplectic HPDEs \eqref{eq-csv} admit the following conservation laws \cite{bridges2001multi}. 
\begin{enumerate}
	\item Multi-symplectic conservation law:
	\begin{equation}\label{eq-csvms}
	\omega_t+  \kappa_x = 0, \qquad  \omega = \omega(d\bz,d\brbz) = Md\bz\cdot d\brbz ,\quad \kappa = \kappa(d\bz,d\brbz) = Kd\bz\cdot d\brbz,
	\end{equation}
	where
	$d\bz$ and $d\brbz$ are one-forms satisfying the variational equation
	\begin{equation*}
	M (d\bz)_t + K(d\bz)_x = \nabla_{\bz\bz}S(\bz)d\bz. 
	\end{equation*}
	\item Local energy conservation law:
	\begin{equation}\label{eq-csve}
	E(\bz)_t + F(\bz)_x = 0,
	\qquad 
	E(\bz) = S(\bz)-\frac{1}{2}K\bz_x\cdot \bz, \qquad F(\bz) = \frac{1}{2} K\bz_t \cdot \bz.
	\end{equation}
	\item Local momentum conservation law:
	\begin{equation}\label{eq-csvm}
	I(\bz)_t +  G(\bz)_x = 0,
	\qquad 
	G(\bz) = S(\bz)-\frac{1}{2}M\bz_t\cdot \bz, \qquad I(\bz) = \frac{1}{2} M\bz_x \cdot \bz.
	\end{equation}
\end{enumerate}

Structure-preserving numerical methods, which preserve certain structures and invariants of the model in the discrete level, have gained many attention in the simulation of mathematical models. It is well-known that such schemes tend to achieve better long time behavior in terms of stability and accuracy for Hamiltonian dynamics \cite{leimkuhler2004simulating}. For the multi-symplectic HPDEs \eqref{eq-csv}, numerical methods preserving \eqref{eq-csvms} are referred as multi-symplectic integrators. There have been many studies on such methods, with examples including the Preissman box schemes \cite{bridges2001multi,zhao2000multisymplectic,ascher2004multisymplectic,cohen2008multi}, Euler box schemes \cite{moore2003backward,cohen2008multi}, diamond schemes \cite{mclachlan2015multisymplectic}, spectral methods \cite{bridges2001multiZK,chen2001multi}, multi-symplectic (partitioned) Runge--Kutta (RK) methods \cite{reich2000multi,hong2006multiPRK,ryland2008multisymplecticity} and recently, DG methods \cite{tang2017discontinuous,cai2018local}. 
These methods have been successfully applied to various equations including the Hamiltonian wave equation \cite{bridges2001multi,mclachlan2015multisymplectic}, the BBM equation \cite{sun2004multi,li2013new}, the CH equation \cite{cohen2008multi}, the KdV equation \cite{zhao2000multisymplectic,ascher2004multisymplectic},  the Schr\"odinger equation \cite{chen2001multi,tang2017discontinuous} and the Dirac equation \cite{hong2006multi}. 
Recently, there have been increasing interests in designing local energy conserving numerical schemes for the continuous dynamical systems.
In the case that the associated energy and momentum functionals of \eqref{eq-csv} are quadratic, some of these multi-symplectic methods also preserve the local conservation of energy \eqref{eq-csve} and momentum \eqref{eq-csvm} at the discrete level, while in general this is not the case. Usually, particular discretizations have to be constructed to preserve the energy or the momentum of \eqref{eq-csv}, such as averaged vector fields (AVF) methods \cite{mclachlan1999geometric,quispel2008new,celledoni2012preserving} and invariant energy quantization (IEQ) methods \cite{cai2019two}.

In this paper, we investigate the structure-preserving property of the DG spatial discretization for the multi-symplectic HPDEs \eqref{eq-csv}. The DG method is a class of finite element methods using discontinuous piecewise polynomial spaces. It was first introduced by Reed and Hill \cite{reed1973triangular} for solving the transport equation and then received its major development in a series of work by Cockburn et al. in  \cite{rkdg1,rkdg2,rkdg3,rkdg4,rkdg5} for hyperbolic conservation laws. 
For equations containing high order spatial derivatives, the local discontinuous Galerkin (LDG) methods were proposed by Cockburn and Shu in \cite{CS1998SINUM}.
The DG method also finds its strength on preserving structures of the continuum equations (sometimes with suitable limiters), such as the positivity or other physical bounds \cite{zhang2010maximum}, the hydrostatic balance \cite{xing2006high,XZS2010}, the entropy inequality \cite{chen2017entropy,sun2018discontinuous} and the asymptotic limits \cite{guermond2010asymptotic}. Recently, there have been many studies in designing DG and LDG methods which can numerically preserve the energy or Hamiltonian structure of the model. Energy conserving LDG methods have been designed for the generalized KdV equation \cite{bona2013conservative, liu2016hamiltonian}, the acoustic wave equation \cite{CSX2014}, the CH equation \cite{liu2016invariant}, the Degasperis-Procesi equation \cite{HLY2014}, the nonlinear Schr\"odinger equation \cite{LKX2015}, the improved Boussinesq equation \cite{LSXC2020} and so on.

This paper can be considered as a step further along this track, and our goal is to design structure-preserving DG methods for the multi-symplectic HPDEs \eqref{eq-csv} which can preserve both the multi-symplectic structure and associated local energy conservation.  
This work is also a reinterpretation and generalization of the earlier framework by Tang et al. \cite{tang2017discontinuous} on DG methods with alternating fluxes for HPDEs. When $K$ has a blocked structure, it was shown in \cite{tang2017discontinuous} that DG methods for \eqref{eq-csv} with alternating fluxes can be rephrased as partitioned RK methods. Then with suitable quadrature rules and time integrators, the fully discrete methods can be interpreted as space-time partitioned RK methods, whose multi-symplecticity follows from \cite{ryland2008multisymplecticity}. In this paper, our attention is particularly on method-of-lines DG schemes. We derive the multi-symplecticity directly from the weak formulation. The analysis covers a very general class of numerical fluxes including, but not limited to, the well-known central and alternating fluxes studied in \cite{cai2018local}. 
We observe that, in the case when neither central nor alternating fluxes give optimal convergence rate, the proposed wider range of numerical fluxes selections may recover the optimal convergence numerically.
Furthermore, these DG methods are shown in \cite{cai2018local} to be energy conserving for two particular examples of HPDEs, namely the Schr\"odinger equation and KdV equation. Motivated by this, we are able to prove that semi-discrete DG methods simultaneously conserve the associated energy for general multi-symplectic HPDEs \eqref{eq-csv}. This result also indicates a rather general approach of designing DG schemes that conserves a certain Hamiltonian invariant of HPDEs. 
If one can construct a multi-syplectic system with the aimed invariant as the associated energy, it will be automatically conserved by the proposed DG approximation applied to such multi-symplectic formulation. For example, we can show that the proposed multi-symplectic DG scheme, with central fluxes, for the KdV equation retrieves the Hamiltonian preserving DG scheme presented in \cite{liu2016hamiltonian}. 

For applications, we especially consider the Hamiltonian wave equation, the BBM equation and the CH equation in this paper. The choices of numerical fluxes and the corresponding implementation procedures are discussed in details. Since the focus of the paper is on spatial discretization, 
for numerical tests, we usually apply high order RK time integrators with small time steps to reduce temporal error, so that fully discrete schemes faithfully approximates method-of-lines schemes. In principle, one can apply a suitable symplecitic RK method to preserve multi-symplecticity, or an AVF method to preserve the energy. However, preserving both properties simultaneously can be difficult for problems with non-quadratic energy functionals.
In the numerical tests, we also observe certain choices of numerical fluxes may help improve the accuracy while still conserving its corresponding energy. Other tests, such as multi-wave interactions, are also provided to illustrate the performance of the DG schemes.

The rest of the paper is organized as follows, in Section \ref{sec-DG} we state DG methods for multi-symplectic system \eqref{eq-csv}. The preservation of multi-symplectic structure and the local energy conservation are proved in Section \ref{sec-prop}. After that, application to various HPDEs and the implementation procedure are discussed in Section \ref{sec-examp}, and numerical tests are provided in Section \ref{sec-num}. Finally, we close the paper with conclusions in Section \ref{sec-concl}.
  
\section{The DG scheme}\label{sec-DG}
\setcounter{equation}{0}
Consider a quasi-uniform partition of the spatial domain $\Omega =\cup_{j=1}^N I_j$, where $I_j=[x_{j-\frac{1}{2}},x_{j+\frac{1}{2}}]$ for $j=1,2,...,N$. The center of each cell is $x_j=\frac{1}{2}(x_{j+\frac{1}{2}}+x_{j-\frac{1}{2}})$, and the mesh size is denoted by $h_j=x_{j+\frac{1}{2}}-x_{j-\frac{1}{2}}$ with $h=\max h_j$ for $j=1,2,...,N$ being the maximal mesh size. We use $P^k(I_j)$ to represent the linear space spanned by polynomials of degrees no more than $k$ on $I_j$. Let 
\begin{equation*}
	V_h = \{v_h\in L^2(\Omega): v_h\big|_{I_j}\in P^k(I_j), j=1,2,...,N\}
\end{equation*}
be the discontinuous piecewise polynomial space and 
$\bV_h = \prod_{l = 1}^m V_h$ be the product space. Since functions $v_h \in V_h$ (or $\bv_h \in \bV_h$) can be double-valued at cell interfaces $x_{j+\hf}$ for all $j$, we use $v_h^+$ and $v_h^-$ ($\bv_h^+$ and $\bv_h^-$) to represent the function limit from the right and left respectively. We denote the average and jump of the functions at the cell interfaces by 
$
\{v_h\} = \hf\left(v_h^++v_h^-\right)
$ ($
\{\bv_h\} = \hf\left(\bv_h^++\bv_h^-\right)
$) and $[v_h] = v_h^+ - v_h^-$ ($[\bv_h] = \bv_h^+ - \bv_h^-$). 

The semi-discrete DG approximation to \eqref{eq-csv} is given as follows:
Seek the numerical solutions $\bz_h \in \bV_h$, such that 
\begin{equation}\label{eq-DG}
\intj M (\bz_h)_t \cdot {\pvphi} dx - \intj
K \bz_h \cdot {\pvphi}_x dx + \left(\widehat{K\bz_h}\cdot{\pvphi}^-\right)_{j+\hf} - \left(\widehat{K\bz_h}\cdot\pvphi^+\right)_{j-\hf} = \intj\nabla_\bz S(\bz_h)\cdot\pvphi dx
\end{equation}
holds for all test functions $\pvphi\in \bV_h$. 
The hatted terms, $\widehat{K\bz_h}$, are the numerical fluxes defined on the element interfaces, and are the key component in designing the DG methods. In this paper, we choose the family of the numerical flux $\widehat{K\bz_h}$ to be
\begin{equation}\label{eq-flux}
	\widehat{K \bz_h } = K\{\bz_h\} + A[\bz_h] + B[\bz_h]_t,
\end{equation}
for any $m\times m$ real symmetric matrix $A$ and real anti-symmetric matrix $B$. 

Some comments on the choice of the numerical flux $\widehat{K\bz_h}$ are given below. We will show that this family of numerical flux includes the central flux and the alternating flux, well-known in the LDG methods applied to equations with high order derivatives. 
Since $K$ is anti-symmetric, there exists an orthogonal matrix $Q$ such that 
	\begin{equation} \label{K_decomp}
	K = 
	\begin{cases}
\displaystyle  Q^T \left(
	\begin{matrix}
	0&0&-\Lambda^T\\
	0&0&0\\
	\Lambda&0 & 0
	\end{matrix}\right)Q, \qquad \text{ if $m$ is odd,}\\
	\\
\displaystyle 	Q^T \left(
	\begin{matrix}
	0&-\Lambda^T\\
	\Lambda&0 
	\end{matrix}\right)Q,\qquad\quad\,\, \text{ if $m$ is even,}
\end{cases}
	\end{equation}
where $\Lambda $ is an $\lfloor{\frac{m}{2}}\rfloor \times \lfloor{\frac{m}{2}}\rfloor$ real matrix. Assume 
\begin{equation}\label{eq-uwv}
Q\bz_h = \left(
\begin{matrix}
\bu_h\\
w_h\\
\bv_h
\end{matrix}\right)\text{ for odd } m \quad \text{ or } \quad Q\bz_h = \left(
\begin{matrix}
\bu_h\\
\bv_h
\end{matrix}\right) \text{ for even } m,
\end{equation}
with $\bu_h$, $\bv_h \in (V_h)^{\lfloor \frac{m}{2}\rfloor}$. If we choose the matrices $A$ and $B$ as
	\begin{equation}\label{eq-A}
	A = 
	\begin{cases}
\displaystyle  \alpha Q^T \left(
	\begin{matrix}
	0&0&\Lambda^T\\
	0&0&0\\
	\Lambda&0 & 0
	\end{matrix}\right)Q, \qquad \text{ if $m$ is odd,}\\
	\\
\displaystyle 	\alpha Q^T \left(
	\begin{matrix}
	0&\Lambda^T\\
	\Lambda&0 
	\end{matrix}\right)Q,\qquad\quad\,\, \text{ if $m$ is even,}
\end{cases}
		\qquad B = 0,
	\end{equation}
with $\alpha \in [-1/2,\,1/2]$, the numerical flux in \eqref{eq-flux} reduces to
	\begin{equation*}
		\widehat{K\bz_h} = Q^T
		\left(
		\begin{matrix}
		-\Lambda^T \left(\{\bv_h\}-\alpha[\bv_h]\right) \\
		0\\
		\Lambda \left(\{\bu_h\}+\alpha[\bu_h]\right)
		\end{matrix}\right) \qquad \text{or} \qquad 
		\widehat{K\bz_h} =	Q^T
		\left(
		\begin{matrix}
		-\Lambda^T \left(\{\bv_h\}-\alpha[\bv_h]\right) \\
		\Lambda \left(\{\bu_h\}+\alpha[\bu_h]\right)
		\end{matrix}\right),
	\end{equation*}
which retrieves alternating fluxes with $\alpha = \pm \hf$, and central fluxes with $\alpha = 0$. 
Second, from the practical point of view, $B = 0$ should work for most situations. While in numerical tests, we notice that  a nonzero $B$ may help improve the accuracy when $A = 0$. We therefore include $B$ in our analysis for completion. Further comments on choices of numerical fluxes are postponed to Remark \ref{rem-choices} in the conclusion.

The DG scheme \eqref{eq-DG} for HPDEs \eqref{eq-csv} also relates to an LDG method for the associated scalar equation. The DG discretization is applied to a first order system reformulated from the original scalar equation and all auxiliary variables will be eliminated. In our cases, the first order system is particularly given as the multi-symplectic system and the elimination procedure will be detailed in Section \ref{sec-examp}.

\section{Properties of the DG scheme}\label{sec-prop}
\setcounter{equation}{0}

In this section, some properties of the proposed semi-discrete DG scheme \eqref{eq-DG} are investigated. More specifically, we are interested in demonstrating that the multi-sympleciticy and the energy conservation are both preserved by our methods. We focus on the spatial discretization only, and will comment on the appropriate temporal discretization at the end of this section.  

\subsection{Multi-symplecticity}
Applying the exterior derivative to the DG scheme \eqref{eq-DG} yields the following variational equation for the one-forms of $\bz_h$
\begin{equation}\label{eq-wgd}
\intj M (d\bz_h)_t \cdot {\pvphi} dx - \intj
K d\bz_h \cdot {\pvphi}_xdx + \left(\widehat{Kd\bz_h}\cdot{\pvphi}^-\right)_{j+\hf} - \left(\widehat{Kd\bz_h}\cdot\pvphi^+\right)_{j-\hf} = \intj\nabla_{\bz\bz} S(\bz_h)d\bz_h\cdot\pvphi dx.
\end{equation}

To facilitate our discussion, the following equalities are provided. 
\begin{LEM}\label{lem-flux}
	Let $\widehat{K\bz_h}$ be defined in \eqref{eq-flux}. For any $\bz_h, \brbz_h\in \bV_h$, we have
	\begin{align}
		K\bz_h^-\cdot \brbz_h^- - \widehat{K\bz_h}\cdot \brbz_h^- + \widehat{K\brbz_h}\cdot \bz_h^-= \cF(\bz_h,\brbz_h) - \hf(B[\brbz_h]\cdot[\bz_h] )_t, \label{eq-fl-mc}\\
	K\bz_h^+\cdot \brbz_h^+ - \widehat{K\bz_h}\cdot \brbz_h^+ + \widehat{K\brbz_h}\cdot \bz_h^+  = \cF(\bz_h,\brbz_h)+\hf(B[\brbz_h]\cdot[\bz_h] )_t,\label{eq-fl-pc}
	\end{align}
	where
	\begin{equation} \label{def:F}
	\begin{aligned}
	\cF(\bz_h,\brbz_h) =&\; \{K\bz_h\cdot \brbz_h\} - \widehat{K\bz_h}\cdot \{\brbz_h\} + \widehat{K\brbz_h}\cdot \{\bz_h\}.
	\end{aligned}
	\end{equation} 
\end{LEM}
\begin{proof}
	We provide the proof of \eqref{eq-fl-mc} only, and skip that of \eqref{eq-fl-pc} which follows along the similar line. 
	Introduce the notation of 
	\begin{equation*}
	\cD(\bz_h,\brbz_h) := K\bz_h^-\cdot \brbz_h^- - \widehat{K\bz_h}\cdot \brbz_h^- + \widehat{K\brbz_h}\cdot \bz_h^- - \cF(\bz_h,\brbz_h),
	\end{equation*} 
	and it yields that
	\begin{equation*}
	\cD(\bz_h,\brbz_h) = \hf\left(-[K\bz_h\cdot \brbz_h] +  \widehat{K\bz_h}\cdot [\brbz_h] - \widehat{K\brbz_h}\cdot [\bz_h]\right),
	\end{equation*} 
	by combining corresponding terms in the subtraction. Using the definition of $\widehat{ K\bz_h}$ in \eqref{eq-flux}, one can obtain
	\begin{equation}\label{eq-lem-1}
	\begin{aligned}
		\cD(\bz_h,\brbz_h)=&\;\hf\left(-[K\bz_h\cdot \brbz_h] +  K \{\bz_h\}\cdot [\brbz_h] - K\{\brbz_h\}\cdot [\bz_h]\right)\\
	 &+\hf\left( A[\bz_h]\cdot[\brbz_h]-A[\brbz_h]\cdot[\bz_h]\right)+ \hf\left(B[\bz_h]_t\cdot [\brbz_h] - B[\brbz_h]_t\cdot [\bz_h] \right).
	\end{aligned}
	\end{equation}
	It can be verified that 
	\begin{align}
		-[K\bz_h\cdot \brbz_h] + K \{\bz_h\}\cdot [\brbz_h] - K\{\brbz_h\}\cdot [\bz_h] =&\; 0,  \notag\\
		A[\bz_h]\cdot[\brbz_h]-A[\brbz_h]\cdot[\bz_h] = &\; 0, \label{eq-lem-2} \\
		B[\bz_h]_t\cdot [\brbz_h] - B[\brbz_h]_t\cdot [\bz_h] = 
		-B[\brbz_h]\cdot [\bz_h]_t -  B[\brbz_h]_t \cdot [\bz_h] = &\; -\left(B[\brbz_h]\cdot [\bz_h]\right)_t, \notag
	\end{align}
	since the matrices $K$, $B$ are anti-symmetric and $A$ is symmetric. 
	After substituting \eqref{eq-lem-2} into \eqref{eq-lem-1}, one can then obtain \eqref{eq-fl-mc}.
\end{proof}

Suppose $d\bz_h, d\brbz_h \in \bV_h$ both satisfy \eqref{eq-wgd}, then we have the following theorem.
\begin{THM}[Conservation of multi-symplecticity]\label{thm-ms}
For any one-forms $d\bz_h, d\brbz_h \in \bV_h$ satisfying the variational equation \eqref{eq-wgd}, with fluxes defined in \eqref{eq-flux}, we have the semi-discrete version of the multi-symplectic conservation laws
\begin{equation*}
\begin{aligned}
\frac{d}{dt} \omega_{h,j}
- {\cF(d\bz_h,d\brbz_h)_{j+\hf} + \cF(d\bz_h,d\brbz_h)_{j-\hf}}= 0,
\end{aligned}
\end{equation*}
where $\cF$ is defined in \eqref{def:F} and 
\begin{equation}\label{omegah}
\omega_{h,j} = \int_{I_j}Md\bz_h\cdot d\brbz_hdx+ \hf \left(B[d\brbz_h]\cdot[d\bz_h]\right)_{j+\hf}+ \hf \left(B[d\brbz_h]\cdot[d\bz_h]\right)_{j-\hf}.
\end{equation}
\end{THM} 
\begin{proof}
	Using the anti-symmetric property of the matrix $M$ and the variational equation \eqref{eq-wgd}, it can be shown that 
	\begin{equation*}
	\begin{aligned}
	 & \frac{d}{dt}\intj Md\bz_h\cdot d\brbz_h  dx\\
	=& \intj M(d\bz_h)_t\cdot d\brbz_h  dx
	+ \intj Md\bz_h\cdot (d\brbz_h)_t  dx
	= \intj M(d\bz_h)_t\cdot d\brbz_h  dx
	- \intj M(d\brbz_h)_t\cdot d\bz_h  dx\\
	=&\intj
	K d\bz_h \cdot {(d\brbz_h)}_x dx - \left(\widehat{Kd\bz_h}\cdot{d\brbz_h}^-\right)_{j+\hf} + \left(\widehat{Kd\bz_h}\cdot{d\brbz_h}^+\right)_{j-\hf} + \intj\nabla_{\bz\bz} S(\bz_h)d\bz_h\cdot d\brbz_h dx\\
	-&\intj
	K d\brbz_h \cdot {(d\bz_h)}_x dx+ \left(\widehat{Kd\brbz_h}\cdot{d\bz_h}^-\right)_{j+\hf} - \left(\widehat{Kd\brbz_h}\cdot{d\bz_h}^+\right)_{j-\hf} - \intj\nabla_{\bz\bz} S(\bz_h)d\brbz_h\cdot d\bz_h dx .  
	\end{aligned}
	\end{equation*}
	With the anti-symmetry of $K$ and the symmetry of $\nabla_{\bz\bz} S(\bz_h)$, we have
	\begin{equation*}\label{eq-wg1}
	\begin{aligned}
	&
	\frac{d}{dt}\intj Md\bz_h\cdot d\brbz_h  dx\\
	=&\intj (Kd\bz_h\cdot d\brbz_h)_x dx + \left(-\widehat{Kd\bz_h}\cdot{d\brbz_h}^-+\widehat{Kd\brbz_h}\cdot{d\bz_h}^-\right)_{j+\hf} 
	 + \left(\widehat{Kd\bz_h}\cdot{d\brbz_h}^+ - \widehat{Kd\brbz_h}\cdot{d\bz_h}^+\right)_{j-\hf} \\
	=& \left(Kd\bz_h^-\cdot d\brbz_h^- - \widehat{Kd\bz_h}\cdot d\brbz_h^- + \widehat{Kd\brbz_h}\cdot d\bz_h^-\right)_{j+\hf}
	-\left(Kd\bz_h^+\cdot d\brbz_h^+ - \widehat{Kd\bz_h}\cdot d\brbz_h^+ + \widehat{Kd\brbz_h}\cdot d\bz_h^+\right)_{j-\hf}.
	\end{aligned}
	\end{equation*}
	With fluxes defined in \eqref{eq-flux}, one can then apply Lemma \ref{lem-flux} to complete the proof. 
\end{proof}
\begin{REM}
	When the anti-symmetric matrix $B$ is chosen as zero, $\omega_{h,j}$ defined in \eqref{omegah} reduces to $\int_{I_j}Md\bz_h\cdot d\brbz_hdx$. 
	Also, when $\bz_h$, $\brbz_h$ are continuous across cell interfaces, we have 
	\begin{equation*}
		\cF(\bz_h,\brbz_h) = - K\bz_h\cdot\brbz_h = K \brbz_h\cdot \bz_h.
	\end{equation*}
\end{REM}

\begin{REM}
	As can be seen from the proof, other choices of fluxes besides \eqref{eq-flux} also preserve the mulit-symplecticity, as long as
	\begin{equation*}
	 \left(Kd\bz_h^-\cdot d\brbz_h^- - \widehat{Kd\bz_h}\cdot d\brbz_h^- + \widehat{Kd\brbz_h}\cdot d\bz_h^-\right)
	-\left(Kd\bz_h^+\cdot d\brbz_h^+ - \widehat{Kd\bz_h}\cdot d\brbz_h^+ + \widehat{Kd\brbz_h}\cdot d\bz_h^+\right)
	\end{equation*}
	is equal to the time derivative of a term consistent with zero. 
\end{REM}


\subsection{Local energy conservation}

In this subsection, the local energy conservation property of the proposed DG scheme \eqref{eq-DG} is explored. 
Recall that the continuous energy $E(\bz)= S(\bz) - \hf K\bz_x\cdot \bz$ is defined in \eqref{eq-csve}, and we have the following
theorem on the local conservation of the discrete energy. 

\begin{THM}[Conservation of energy]\label{thm-ec}
	The numerical solution $\bz_h$ of the semi-discrete DG scheme \eqref{eq-DG} satisfies the local energy conservation law in the form of 
	\begin{equation} \label{localenergy}
	\frac{d}{dt}\cE_{h,j}
	 + \hf \cF(\bz_h,(\bz_h)_t)_{j+\hf} - \hf \cF(\bz_h,(\bz_h)_t)_{j-\hf} = 0,
	\end{equation}
	where 
	\begin{equation*}
	\cE_{h,j} = \int_{I_j}E(\bz_h)dx -\left(\hf\widehat{K\bz_h}\cdot\bz_h^-+\frac{1}{4}\left(B[\bz_h]_t\cdot[\bz_h]\right)\right)_{j+\hf}+\left(\hf\widehat{K\bz_h}\cdot \bz_h^+-\frac{1}{4}B[\bz_h]_t\cdot[\bz_h]\right)_{j-\hf}.
	\end{equation*}
	As a consequence, it conserves the total energy 
	\begin{equation}\label{eq-energy}
	\cE_h = \int_\Omega E(\bz_h) dx + \hf \sum_{j}\big(\left(K\{\bz_h\}+A[\bz_h]\right)\cdot [\bz_h]\big)_{j+\hf}.
	\end{equation}
\end{THM}
\begin{proof} Following the definition of the energy $E(\bz)$, we have
	\begin{equation}\label{eq-ener-1}
	\begin{aligned}
	\frac{d}{dt}\int_{I_j} E(\bz_h) dx = 
	&\;\frac{d}{dt}\int_{I_j} \Big( S(\bz_h) - \hf K(\bz_h)_x\cdot \bz_h \Big) dx\\
	=&\;\int_{I_j} \nabla_\bz S(\bz_h)\cdot (\bz_h)_t dx-\hf\int_{I_j} \Big( K(\bz_h)_{xt}\cdot \bz_h+ K(\bz_h)_x\cdot(\bz_h)_{t} \Big) dx.
	\end{aligned}
	\end{equation}
	The numerical solution $\bz_h$ satisfies the variational equation \eqref{eq-DG}. Taking the test function $\pvphi=(\bz_h)_t$ leads to
	\begin{equation}
	\begin{aligned}\label{eq-ener-S}
	&\int_{I_j} \nabla_\bz	 S(\bz_h)\cdot (\bz_h)_tdx\\ =& \int_{I_j}M(\bz_h)_t\cdot(\bz_h)_t dx -\int_{I_j} K\bz_h \cdot (\bz_h)_{tx} dx + \left(\widehat{K\bz_h}\cdot (\bz_h^-)_t\right)_{j+\hf} - \left(\widehat{K\bz_h}\cdot (\bz_h^+)_t\right)_{j-\hf}\\
		=&\int_{I_j} K (\bz_h)_{tx}\cdot \bz_h dx + \left(\widehat{K\bz_h}\cdot (\bz_h^-)_t\right)_{j+\hf} - \left(\widehat{K\bz_h}\cdot (\bz_h^+)_t\right)_{j-\hf},
	\end{aligned}
	\end{equation}
	where the last equality follows from the anti-symmetry of $M$ and $K$. 
	By substituting \eqref{eq-ener-S} into \eqref{eq-ener-1}, we have
	\begin{equation*}
		\begin{aligned}
		&\frac{d}{dt}\int_{I_j} E(\bz_h) dx\\
		 =&\; \hf\int_{I_j} \Big( K(\bz_h)_{tx}\cdot \bz_h- K(\bz_h)_x\cdot(\bz_h)_{t}  \Big) dx 
		 + \left(\widehat{K\bz_h}\cdot (\bz_h^-)_t\right)_{j+\hf} - \left(\widehat{K\bz_h}\cdot (\bz_h^+)_t\right)_{j-\hf}\\
		=&\; -\frac{1}{2}\int_{I_j} \big( K\bz_h \cdot (\bz_h)_t \big)_x dx+ \left(\widehat{K\bz_h}\cdot (\bz_h^-)_t\right)_{j+\hf} - \left(\widehat{K\bz_h}\cdot (\bz_h^+)_t\right)_{j-\hf}\\
		=&\; \left(\widehat{K\bz_h}\cdot (\bz_h^-)_t-\hf K\bz_h^-\cdot (\bz_h^-)_t\right)_{j+\hf} - \left(\widehat{K\bz_h}\cdot (\bz_h^+)_t-\hf K\bz_h^+\cdot (\bz_h^+)_t\right)_{j-\hf}.
		\end{aligned}
	\end{equation*}
	On the other hand, using the product rule, it can be shown that
	\begin{align*}
		-\hf\left(\widehat{K\bz_h}\cdot \bz_h^-\right)_t =&\; -\hf\widehat{K(\bz_h)_t}\cdot \bz_h^- -\hf\widehat{K\bz_h}\cdot (\bz_h^-)_t, \\
		\hf\left(\widehat{K\bz_h}\cdot \bz_h^+\right)_t =&\; \hf\widehat{K(\bz_h)_t}\cdot \bz_h^+ +\hf\widehat{K\bz_h}\cdot (\bz_h^+)_t.
	\end{align*}
	Therefore 
	\begin{equation}\label{eq-ener-I}
	\begin{aligned}
		&\frac{d}{dt} \left(\int_{I_j} E(\bz_h) dx -\hf(\widehat{K\bz_h}\cdot\bz_h^-)_{j+\hf}+\hf(\widehat{K\bz_h}\cdot \bz_h^+)_{j-\hf}\right)\\
		 = & -\hf\left(K\bz_h^-\cdot (\bz_h^-)_t -\widehat{K\bz_h}\cdot (\bz_h^-)_t+\widehat{K(\bz_h)_t}\cdot \bz_h^-\right)_{j+\hf}\\
		& + \hf\left(K\bz_h^+\cdot (\bz_h^+)_t -\widehat{K\bz_h}\cdot (\bz_h^+)_t+\widehat{K(\bz_h)_t}\cdot \bz_h^+\right)_{j-\hf}\\
		= & -\hf\left(\cF(\bz_h,(\bz_h)_t) - \hf \left(B[\bz_h]_t\cdot[\bz_h]\right)_t\right)_{j+\hf}
		 + \hf\left(\cF(\bz_h,(\bz_h)_t) + \hf \left(B[\bz_h]_t\cdot[\bz_h]\right)_t\right)_{j-\hf},
		\end{aligned}
	\end{equation}
	where the last equality follows from Lemma \ref{lem-flux}. The local energy conservation property \eqref{localenergy} can be obtained after rearranging terms of \eqref{eq-ener-I}. The global energy conservation property \eqref{eq-energy} follows from summing it over all elements and utilizing the definition of the numerical flux in \eqref{eq-flux} and the periodic boundary condition. 
\end{proof}

The global energy term $\cE_h$ in \eqref{eq-energy} can be further simplified for some special choices of the matrices $A$ and $B$ in the definition of the numerical flux, which is summarized in the following corollary. 

\begin{COR}\label{cor-ec}
	With the decomposition of $K$ in \eqref{K_decomp}, if the matrix $B$ is chosen as $\beta M$ and the last $\lfloor \frac{m}{2}\rfloor$ columns of $MQ^T$ are zero, then the global energy term $\cE_h$ in \eqref{eq-energy} reduces to
	\begin{equation}\label{eq-corec}
	\cE_h =\int_\Omega \Big( S(\bz_h) - \nabla_{\bv} S(\bz_h) \cdot \bv_h \Big) dx + \hf \sum_j \left(\widetilde{A}[\bz_h]\cdot [\bz_h]\right)_{j+\hf},
	\end{equation}
	where $\bv_h \in (V_h)^{\lfloor \frac{m}{2}\rfloor}$ is defined in \eqref{eq-uwv} and 
	\begin{equation} \label{tildeA_decomp}
	\widetilde{A} = 
	\begin{cases}
\displaystyle  Q^T \left(
	\begin{matrix}
	I&0&0\\
	0&1&0\\
	0&0&-I
	\end{matrix}\right)QA, \qquad \text{ if $m$ is odd,}\\
	\\
\displaystyle 	Q^T \left(
	\begin{matrix}
	I&0\\
	0&-I 
	\end{matrix}\right)QA,\qquad\quad\,\, \text{ if $m$ is even.}
\end{cases}
	\end{equation}
	 In particular, for $A$ defined in \eqref{eq-A}, $\widetilde{A}$ is anti-symmetric, hence $\widetilde{A}[\bz_h]\cdot [\bz_h] = 0$, which leads to 
	\begin{equation*}
	\cE_h =\int_\Omega S(\bz_h) - \nabla_{\bv} S(\bz_h) \cdot \bv_hdx .
	\end{equation*}
	 
\end{COR}

\begin{proof}
	Without loss of generality, we assume $m$ is odd in this proof. We first derive the variational form associated with the last $\lfloor \frac{m}{2}\rfloor$ equations. Note the DG scheme \eqref{eq-DG} has the following global form after summing over all the elements: Find $\bz_h \in \bV_h$, such that for all $\pvphi\in \bV_h$, we have
	\begin{equation}\label{eq-DG-global}
	\int_\Omega M (\bz_h)_t \cdot {\pvphi} dx - \int_\Omega
	K \bz_h \cdot {\pvphi}_x dx -\sum_j \left(\widehat{K\bz_h}\cdot{[\pvphi]}\right)_{j+\hf} = \int_\Omega\nabla_\bz S(\bz_h)\cdot\pvphi dx.
	\end{equation}
	We take $\pvphi = Q^T\left(
	\begin{matrix}
	\pmb{0}\\ 0\\ \bxi
	\end{matrix}\right)$ with $\bxi \in (V_h)^{\lfloor \frac{m}{2}\rfloor}$ in \eqref{eq-DG-global}. Under our assumption on $M$, one can obtain
	\begin{equation}\label{eq-cor-1}
	M(\bz_h)_t\cdot \pvphi = - (\bz_h)_t \cdot MQ^T \left(\begin{matrix}
	\pmb{0}\\
	0\\
	\bxi
	\end{matrix}\right) = 0.
	\end{equation}
	Furthermore, denote $Q\bz_h = \left(
	\begin{matrix}
	\bu_h\\w_h\\\bv_h
	\end{matrix}\right)$, and we have
	\begin{equation}\label{eq-cor-2}
		K\bz_h \cdot \pvphi_x =  Q^T \left(
		\begin{matrix}
		0&0&-\Lambda^T\\
		0&0&0\\
		\Lambda&0 & 0
		\end{matrix}\right)\left(
		\begin{matrix}
		\bu_h\\
		w_h\\
		\bv_h
		\end{matrix}\right)\cdot Q^T\left(
		\begin{matrix}
		\pmb{0}\\ 0\\ \bxi_x
		\end{matrix}\right) = \Lambda \bu_h \cdot \bxi_x.
	\end{equation}
	Recall that $B = \beta M$, which leads to $B[\bz_h]_t\cdot [\pvphi]=0$ following the steps in \eqref{eq-cor-1}. Similar computation as in \eqref{eq-cor-2} yields
	\begin{equation*}
	\widehat{ K\bz_h}\cdot [\pvphi] = K\{\bz_h\}\cdot [\pvphi] + \left(A[\bz_h]+B[\bz_h]_t\right)\cdot [\pvphi] = \Lambda\{\bu_h\}\cdot[\bxi] +  A[\bz_h]\cdot [\pvphi].
	\end{equation*}
	Moreover, since
	\begin{equation*}
		\nabla_\bz S(\bz_h)\cdot \pvphi = Q\nabla_\bz S(\bz_h) \cdot \left(\begin{matrix}
		\pmb{0}\\
		0\\
		\bxi
		\end{matrix}\right) = \nabla_{Q\bz} S(\bz_h) \cdot \left(\begin{matrix}
		\pmb{0}\\
		0\\
		\bxi
		\end{matrix}\right)= \nabla_{\bv} S(\bz_h)\cdot \bxi, 
	\end{equation*}
	 the DG scheme \eqref{eq-DG-global} becomes
	\begin{equation}\label{eq-DG-2}
		-\int_\Omega \Lambda \bu_h\cdot \bxi_x dx -\sum_j \left(\Lambda \{\bu_h\} \cdot [\bxi]\right)_{j+\hf} -\sum_j \left(A[\bz_h]  \cdot [\pvphi]\right)_{j+\hf}= \int_\Omega \nabla_{\bv}S(\bz_h)\cdot \bxi dx.
	\end{equation}
	
	We now turn to simplify the energy functional $\cE_h$ defined in \eqref{eq-energy}. One can apply integration by parts to get
	\begin{equation}\label{eq-Kzxz}
	\begin{aligned}
		\int_\Omega K(\bz_h)_x\cdot \bz_h dx
		 =& \int_\Omega \Lambda(\bu_h)_x\cdot \bv_h dx - \int_\Omega \Lambda\bu_h\cdot (\bv_h)_x dx
		 = - 2\int_\Omega \Lambda\bu_h\cdot (\bv_h)_x dx -\sum_j [\Lambda \bu_h \cdot \bv_h]_{j+\hf}.
	\end{aligned}
	\end{equation}
	By taking $\bxi = \bv_h$ in \eqref{eq-DG-2} and combining it with \eqref{eq-Kzxz}, we have
	\begin{equation*}
		\int_\Omega K(\bz_h)_x\cdot \bz_h dx
				= 2\int_\Omega \nabla_{\bv}S(\bz_h)\cdot \bv_h dx
				+\sum_j \Big(2\Lambda\{ \bu_h\}\cdot [\bv_h] - [\Lambda \bu_h \cdot \bv_h]+2A[\bz_h]\cdot[\pvphi_{\bv}]\Big)_{j+\hf},
	\end{equation*}
	where $\pvphi_{\bv} =Q^T\left(\begin{matrix}
	\pmb{0}\\0\\\bv_h
	\end{matrix}\right)$. On the other hand, it can be shown that
	\begin{equation*}
		\sum_j \Big(\left(K\{\bz_h\}+A[\bz_h]\right)\cdot [\bz_h]\Big)_{j+\hf} =\sum_j \Big( \Lambda\{\bu_h\}\cdot [\bv_h] -\Lambda[\bu_h]\cdot \{\bv_h\} +A[\bz_h]\cdot [\bz_h]\Big)_{j+\hf}.
	\end{equation*}
	Therefore, we can rewrite the energy functional $\cE_h$ as 
	\begin{align}\label{eq-energy-proof1}
		\cE_h =& \int_\Omega S(\bz_h) dx - \hf \int_\Omega K(\bz_h)_x\cdot \bz_h dx +
		 	\hf \sum_j\Big(\left(K\{\bz_h\}+A[\bz_h]\right)\cdot [\bz_h]\Big)_{j+\hf} \notag \\
		=&\int_\Omega \left( S(\bz_h) - \nabla_{\bv} S(\bz_h) \cdot \bv_h \right) dx 
		+ \hf \sum_j \left(-{\Lambda\{\bu_h\}}\cdot [\bv_h]+[\Lambda \bu_h\cdot \bv_h ] -\Lambda[\bu_h]\cdot \{\bv_h\}\right)_{j+\hf}  \notag \\
		 &+\sum_j \left(A[\bz_h]\cdot \left(\hf [\bz_h] - [\pvphi_{\bv}]  \right)\right)_{j+\hf}.
	\end{align}
	One can easily verify that 
	\begin{equation}\label{eq-flux-id}
	\sum_j \Big([\Lambda \bu_h\cdot \bv_h ] -
	 {\Lambda\{\bu_h\}}\cdot [\bv_h]-\Lambda[\bu_h]\cdot \{\bv_h\}\Big)_{j+\hf} = 0.
	\end{equation}
	Furthermore, we have
	\begin{equation}\label{eq-tildeA}
	\begin{aligned}
	A[\bz_h]\cdot \left(\hf [\bz_h] - [\pvphi_{\bv}]\right) &= 
	A[\bz_h]\cdot \hf Q^T\left(\begin{matrix}
	[\bu_h]\\
	[w_h]\\
	-[\bv_h]
	\end{matrix}\right)=	\hf Q^T\left(\begin{matrix}
	I&0&0\\
	0&1&0\\
	0&0&-I
	\end{matrix}\right)Q
	A[\bz_h]\cdot [\bz_h].
	\end{aligned}
	\end{equation}
	The equality \eqref{eq-corec} is obtained by substituting \eqref{eq-flux-id} and \eqref{eq-tildeA} into \eqref{eq-energy-proof1}, and using the definition of $\widetilde{A}$ in \eqref{tildeA_decomp}. This completes the proof.
\end{proof}

\section{Application and implementation}\label{sec-examp}
\setcounter{equation}{0}
In this section, we present the DG discretization of a few particular examples of the multi-symplectic HPDEs, including the Hamiltonian wave equation, the BBM equation, the KdV equation, the Schr\"odinger equation and the CH equation. Our main attention is on various numerical fluxes and discrete energy conservation indicated by Theorem \ref{thm-ec}. We detail the choices of numerical fluxes and implementation of our DG method for the Hamiltonian wave equation, the BBM equation and the CH equation. It will be shown that for these models, the proposed methods with appropriate numerical fluxes recover some existing energy conserving DG methods studied in the literature. We have also listed the application to the Schr\"odinger equation and the KdV equation, while these schemes (with central and alternating fluxes) have been studied in \cite{cai2018local}.

For all of these examples of HPDEs, the original equation involves only one unknown $u$, and they can be rewritten in the form of HPDEs by introducing auxiliary variables $\bz$. The DG method \eqref{eq-DG} is presented with the unknown variable $\bz_h$. 
In the numerical implementation of the DG scheme, we want to eliminate all the intermediate unknowns and convert the system back into an equation involving only one unknown $u_h$, which will be explained in this section. The following operators are useful in this procedure to simplify the notations. 
\begin{itemize}
	\item $D_{\alpha}: V_h \to V_h$, such that
	\begin{equation}\label{eq-not-Da}
	\int_{I_j} (D_{\alpha}u_h) \vphi dx = - \int_{I_j} u_h \vphi_x dx + (\widehat{u_h} \vphi^-)_{j+\hf} - (\widehat{u_h} \vphi^+)_{j-\hf}, \qquad \forall \vphi\in V_h,
	\end{equation}
	where $\widehat{u_h} = \{u_h\}+ \alpha [u_h] $. 
	
	\item $L: V_h \to V_h$, such that
	\begin{equation}\label{eq-not-L}
	\int_{I_j} (Lu_h) \vphi dx = [u_h]_{j+\hf}\vphi_{j+\hf}^--[u_h]_{j-\hf}\vphi_{j-\hf}^+, \qquad \forall \vphi\in V_h.
	\end{equation}
	
	\item We also use 
	\begin{equation}\label{eq-not-Pi}
		\Pi:L^2(\Omega)\to V_h
	\end{equation} to represent the standard $L^2$ projection to the piecewise polynomial space $V_h$.
\end{itemize}
With a given set of basis functions, matrices for the local operator $D_\alpha$ and $L$ can be assembled explicitly in each interval $I_j$. The projection $\Pi$ can also be implemented as a computer subroutine.
To implement a DG scheme, for example \eqref{scheme:wave} for the Hamiltonian wave equation, one simply replaces trial functions and operators with corresponding vectors and matrices or computer subroutines. 

\subsection{Hamiltonian wave equation}
	The Hamiltonian wave equation takes the form of
	\begin{equation}\label{eq-wv}
	u_{tt} - u_{xx} = V'(u), 
	\end{equation}
	and includes the linear second order wave, Sine-Gordon and Klein-Gordon equations as special cases.
	By introducing two auxiliary variables $v$ and $w$, it can be formulated as a multi-symplectic system
	\eqref{eq-csv} with $\bz = \left(u\ v\ w\right)^T$ and
	\begin{equation*}
	M = \left(\begin{matrix}
	0&-1&0\\
	1& 0&0\\
	0& 0&0
	\end{matrix}\right),\quad 
	K = \left(\begin{matrix}
	0& 0&1\\
	0& 0&0\\
	-1& 0&0
	\end{matrix}\right),\quad 	S(\bz) = \hf\left(v^2-w^2\right) - V(u),
	\end{equation*}
	which can be expressed as
	\begin{equation} \label{eq-wv2}
	\left\{
	 \begin{aligned}
	  -v_t+w_x=&\;-V'(u), \\
	  u_t = &\;v, \\
	  -u_x = &\;-w.
	 \end{aligned}\right.
	\end{equation}
	By utilizing the fact that $u_x=w$, the energy functional $\cE$ of the continuum equation can be simplified as 
	\begin{align*}
	\cE =\int_\Omega E(\bz) dx & =  \int_\Omega \Big(S(\bz)-\frac{1}{2}K\bz_x\cdot \bz \Big) dx = \int_\Omega \Big(\hf\left(v^2-w^2\right)-V(u) -\hf uw_x+\hf u_xw\Big) dx \\
	& = \int_\Omega \Big(\hf\left(v^2+w^2\right)-V(u)\Big) dx =  \int_\Omega \Big(\hf\left(u_t^2+u_x^2\right)-V(u)\Big) dx.
	\end{align*}

	Applying DG discretization \eqref{eq-wgd} to \eqref{eq-wv2} gives the following scheme: 
	Find $u_h$, $v_h$, $w_h \in V_h$, such that for all  $\vphi_i \in V_h$, $i = 1,2,3$, we have
	\begin{subequations}\label{eq-DG-wv}
	\begin{align}
	-\int_{I_j} (v_h)_t \vphi_1 dx - \int_{I_j} w_h (\vphi_1)_x dx + (\widehat{w_h} \vphi_1^-)_{j+\hf}-(\widehat{w_h} \vphi_1^+)_{j-\hf}  &= - \int_{I_j} V'(u_h)\vphi_1 dx, \\
	\int_{I_j} (u_h)_t \vphi_2 dx +(\widehat{v_h}\vphi_2^-)_{j+\hf}-(\widehat{v_h}\vphi_2^+)_{j-\hf}&= 	\int_{I_j} v_h \vphi_2 dx,\\
	 \int_{I_j} u_h (\vphi_3)_x dx - (\widehat{u_h} \vphi_3^-)_{j+\hf}+(\widehat{u_h} \vphi_3^+)_{j-\hf}&=-\int_{I_j} w_h \vphi_3 dx,
	\end{align}
	\end{subequations}
	where $\widehat{w_h}$, $\widehat{v_h}$, and $\widehat{u_h}$ are the components of the numerical flux $\widehat{K\bz_h}$. Here we limit the choices of the matrices $A$, $B$ in \eqref{eq-flux} to $A = \left(\begin{matrix}
	\alpha_{11}&0&\alpha_{13}\\
	0&0&0\\
	\alpha_{13}&0&\alpha_{33}\\
	\end{matrix}\right)$ and 
	$B = \left(\begin{matrix}
	0&-\beta&0\\
	\beta&0&0\\
	0&0&0\\
	\end{matrix}\right)$. As a result, these numerical fluxes take the form of 	
	\begin{subequations}\label{eq-DG-wv-fl}
	\begin{align}
		\widehat{w_h} =&\; \{w_h\} + \alpha_{11}[u_h] + \alpha_{13}[w_h] - \beta[v_h]_t,\qquad
		\widehat{v_h} = \beta[u_h]_t ,\\
		\widehat{u_h} =&\; \{u_h\}-\alpha_{13}[u_h]  - \alpha_{33}[w_h].
	\end{align}
	\end{subequations}
	
Note that $Q = I$ and $\widetilde{A} = 
\left(
\begin{matrix}
	\alpha_{11}&0&\alpha_{13}\\
	0&0&0\\
-\alpha_{13}&0&-\alpha_{33}\\
\end{matrix}\right)$. Corollary \ref{cor-ec} states the following discrete energy conservation property. 
\begin{PROP}
	The DG scheme \eqref{eq-DG-wv}, with numerical fluxes given in \eqref{eq-DG-wv-fl}, conserves the discrete energy 
	\begin{equation*}
	\cE_h = \int_\Omega \hf\left(v_h^2 + w_h^2 \right) - V(u_h)dx +  \hf\sum_{j}\left({\alpha_{11}}[u_h]^2-\alpha_{33}[w_h]^2\right)_{j+\hf}.
	\end{equation*}
\end{PROP}

	The energy conserving LDG method (with alternating fluxes) for the second order wave equation has been extensively studied in \cite{xing2013energy}, where the optimal error estimate, energy conserving, superconvergence properties are carefully analyzed. It was also demonstrated numerically, that the shape of the solution, after long time integration, is well preserved due to the energy conserving property. The DG method \eqref{eq-DG-wv}, with alternating fluxes $\alpha_{11} = \alpha_{33} = 0$, $\alpha_{13} = \pm\hf$ and $\beta = 0$ in \eqref{eq-DG-wv-fl}, retrieves the same energy conserving LDG method in \cite{xing2013energy}, which means that the proposed semi-discrete method is both multi-symplectic and energy conserving.
	
We now comment on implementation of \eqref{eq-DG-wv} with numerical fluxes \eqref{eq-DG-wv-fl}. With notations in \eqref{eq-not-Da}, \eqref{eq-not-L} and \eqref{eq-not-Pi}, the scheme can be rewritten as follows.
\begin{subequations}\label{eq-wavescheme}
\begin{align}
- (v_h)_t + D_{\alpha_{13}}w_h + \alpha_{11}L u_h- \beta L(v_h)_t&= - \Pi (V'(u_h)),\\
(u_h)_t +\beta L(u_h)_t &= v_h,\\
-\left(D_{-\alpha_{13}}u_h - \alpha_{33} L w_h\right) &= - w_h ,
\end{align}
\end{subequations}
which can be further simplified as
	\begin{equation}\label{scheme:wave}
	(u_h)_{tt}=  (I + \beta L)^{-2}\left(D_{\alpha_{13}}(I+\alpha_{33}L)^{-1}D_{\alpha_{-13}}u_h -\alpha_{11} L u_h  - \Pi (V'(u_h))\right).
	\end{equation}

\subsection{Nonlinear KdV equation}
Next, we consider the nonlinear KdV equation 
\begin{equation}\label{eq-kdv}
u_t + \eta u u_x  + \veps^2 u_{xxx} = 0,
\end{equation}
where $\eta, \veps$ are given parameters. It can be written as a multi-symplectic system \eqref{eq-csv} with $\bz = (\phi\ u\ v\ w)^T$ and
\begin{equation*}
M =  \left(\begin{matrix}
0& \hf&0 &0\\
-\hf& 0&0 &0\\
0& 0&0 &0\\
0& 0&0 &0\\
\end{matrix}\right),\quad 
K = \left(\begin{matrix}
0& 0& 0&1\\
0& 0& -\veps&0\\
0& \veps& 0&0\\
-1& 0& 0& 0
\end{matrix}\right), \quad S(\bz) = \hf v^2 - uw + \frac{\eta}{6} u^3,
\end{equation*}
	which can be expressed as
	\begin{equation} \label{eq-kdv2}
	\left\{
	 \begin{aligned}
	  \hf u_t+w_x=&\;0, \\
	  -\hf \phi_t - \veps v_x =&\; -w+\frac{\eta}{2} u^2, \\
	  \veps u_x =&\; v, \\
	  -\phi_x =&\; -u.
	 \end{aligned}\right.
	\end{equation}
The continuum equation is associated with the energy functional 
\begin{equation} \label{kdvE3}
\cE = \int_\Omega E(\bz) dx =  \int_\Omega \frac{\eta}{6} u^3 - \hf v^2 dx =  \int_\Omega \frac{\eta}{6} u^3 - \frac{1}{2} (\veps u_x)^2 dx.
\end{equation}

Multi-symplectic DG method for the nonlinear KdV equation \eqref{eq-kdv} can be obtained by applying the method \eqref{eq-wgd} to \eqref{eq-kdv2}, 
with the choices $A = \left(\begin{matrix}
0&0&0&\alpha_1\\
0&0&\veps \alpha_2&0\\
0&\veps \alpha_2&0&0\\
\alpha_1&0&0&0\\
\end{matrix}\right)$ and $B = 0$. To be more specific, we look for $u_h, v_h, w_h \in V_h$, such that for all $\vphi_i \in V_h$, $i = 1, 2,3,4$, we have
\begin{subequations}\label{eq-dg-kdv}
	\begin{align}
		\hf \int_{I_j} (u_h)_t \vphi_1 dx - \int_{I_j}w_h (\vphi_1)_x dx + (\widehat{w_h}\vphi_1^-)_{j + \hf} - (\widehat{w_h}\vphi_1^+)_{j - \hf} &= 0,\label{eq-kdv-1}\\
		-\hf \int_{I_j} (\phi_h)_t \vphi_2 dx +\veps\left( \int_{I_j} v_h (\vphi_2)_x dx - (\widehat{v_h} \vphi_2^-)_{j+\hf} + (\widehat{v_h} \vphi_2^+)_{j_-\hf} \right) &= \int_{I_j} (-w_h + \frac{\eta}{2}u_h^2)\vphi_2 dx,\label{eq-kdv-2}\\
		\veps\left( - \int_{I_j} u_h (\vphi_3)_x dx + (\widehat{u_h}\vphi_3^-)_{j+\hf} - (\widehat{u_h}\vphi_3^+)_{j-\hf}\right) &= \int_{I_j} v_h \vphi_3 dx, \label{eq-kdv-3}\\
		\int_{I_j} \phi_h (\vphi_4)_x dx - (\widehat{\phi_h}\vphi_4^-)_{j+\hf} +(\widehat{\phi_h}\vphi_4^+)_{j-\hf} &= -\int_{I_j} u_h \vphi_4 dx,\label{eq-kdv-4}
	\end{align}
\end{subequations}
where 
\begin{equation}\label{eq-dg-kdv-flux}
	\widehat{w_h} = \{w_h\} + \alpha_1 [w_h], \quad \widehat{v_h} = \{v_h\} - \alpha_2 [v_h], \quad
	\widehat{u_h} = \{u_h\} + \alpha_2 [u_h], \quad 
	\widehat{\phi_h} = \{\phi_h\} - \alpha_1 [\phi_h]. 
\end{equation}
If we choose $\alpha_1=0$, and $\alpha_2 = 0$ (leading to central flux) or $\alpha_2 =  \pm \hf$ (leading to alternating flux), the proposed methods reduce to the same ones that
was studied in \cite{tang2017discontinuous}. 

With $Q = I$ and $\tilde{A} = 
\left(
\begin{matrix}
0&0&0&\alpha_1\\
0&0&\veps\alpha_2&0\\
0&-\veps\alpha_2&0&0\\
-\alpha_1&0&0&0\\
\end{matrix}\right)$, 
we have the following discrete energy conservation from Corollary \ref{cor-ec}.
\begin{PROP}
	For any $\alpha_1$ and $\alpha_2$, the DG scheme \eqref{eq-dg-kdv} with numerical fluxes \eqref{eq-dg-kdv-flux} conserves the discrete energy 
	\begin{equation*}
	\cE_h 
	= \int_\Omega \left( \frac{\eta}{6}u_h^3-\hf v_h^2 \right) dx.
	\end{equation*}
\end{PROP}

As for the implementation, we look into the DG scheme in the operator form. 
\begin{align*}
\hf \left(u_h\right)_t + D_{\alpha_{1}}w_h &= 0,\\
-\hf (\phi_h)_t - \veps D_{-\alpha_{2}}v_h &= -w_h +  \frac{\eta}{2}\Pi u_h^2,\\
\veps D_{\alpha_2}u_h &= v_h,\\
-D_{-\alpha_1}\phi_h &= -u_h.
\end{align*}
When $
\alpha_1 = 0$, one can eliminate all auxiliary variables to obtain 
\begin{equation*}
	(u_h)_t = -\veps^2 D_{0}D_{-\alpha_2}D_{\alpha_2} u_h -\frac{\eta}{2} D_0 \Pi u_h^2.
\end{equation*}

	The ``energy conserving'' DG method, which conserves the $\int_\Omega u_h^2dx$ (the momentum in the multi-symplectic community), have been recently proposed and studied in \cite{bona2013conservative}. Suboptimal error estimate was provided. Later, the Hamiltonian preserving DG method, which conserves the energy \eqref{kdvE3}, was presented in \cite{liu2016hamiltonian}, and it was shown in the paper that the Hamiltonian preserving DG methods have a slightly improved long time behavior over those in \cite{bona2013conservative}. The DG method \eqref{eq-dg-kdv}, with central fluxes, retrieves the same Hamiltonian conserving DG method in \cite{liu2016hamiltonian}, which means that semi-discrete method is both multi-symplectic and energy conserving. 
	To see connections between the scheme  \eqref{eq-dg-kdv} and that in \cite{liu2016hamiltonian}, we apply $\hf \partial_t$ to \eqref{eq-kdv-4}, take $\vphi_4 = \vphi_1$ and then combine it with \eqref{eq-kdv-1} to get
	\begin{equation}\label{eq-kdv-1-1}
	\int_{I_j} (u_h)_t \vphi_{1} dx - \int_{I_j}(w_h-\hf (\phi_h)_t) (\vphi_1)_x dx + (\{w_h -\hf (\phi_h)_t\}\vphi_1^-)_{j+\hf} - (\{w_h -\hf (\phi_h)_t\}\vphi_1^+)_{j-\hf} = 0.
	\end{equation}
	Let $q_h = w_h -\hf (\phi_h)_t$ in \eqref{eq-kdv-1-1} and \eqref{eq-kdv-2}. Then \eqref{eq-kdv-1-1}, \eqref{eq-kdv-2} and \eqref{eq-kdv-3} recover the scheme in 
	\cite{liu2016hamiltonian}. 
	
\subsection{BBM equation}
The BBM equation 
\begin{equation*}
u_t-\sigma u_{xxt}+uu_x=0,
\end{equation*}
with $\sigma$ being a given positive parameter, models long shallow water waves of small amplitude and is widely studied.
It can be reformulated as (by introducing $V(u)=-u^3/6$)
\begin{equation}\label{eq-bbm}
u_t - \sigma u_{xxt}  = V'(u)_x, 
\end{equation}
which has the multi-symplectic form \eqref{eq-csv} with $\bz = \left(\phi\ u\ v\ w\ p\right)^T$, and
\begin{equation*}
\begin{aligned}
M &=  \left(\begin{matrix}
0& -\hf&0 &0&0\\
\hf& 0&-\frac{\sigma}{2} &0&0\\
0& \frac{\sigma}{2}&0 &0&0\\
0& 0&0 &0&0\\
0& 0&0 &0&0\\
\end{matrix}\right),\quad 
K = \left(\begin{matrix}
0& 0& 0& 0&-1\\
0& 0& 0& -\frac{\sigma}{2}&0\\
0& 0& 0&0&0\\
0& \frac{\sigma}{2}& 0& 0&0\\
1& 0& 0& 0&0
\end{matrix}\right),\quad 	
S(\bz) = up-V(u)+\frac{\sigma}{2}v w.
\end{aligned}
\end{equation*}
	Then the BBM equation \eqref{eq-bbm} can be expressed as
	\begin{equation*}
	\left\{ \begin{aligned}
	  -\hf u_t-p_x=&\;0, \\
	  \hf \phi_t - \frac{\sigma}{2} v_t - \frac{\sigma}{2} w_x =&\; p-V'(u), \\
	  \frac{\sigma}{2} u_t =&\; \frac{\sigma}{2} w, \\
	  \frac{\sigma}{2} u_x =&\; \frac{\sigma}{2} v, \\
	  \phi_x =&\; u.
	 \end{aligned}\right.
	\end{equation*}
After simplification, the corresponding energy functional is 
\begin{equation*}
\cE = \int_\Omega E(\bz) dx =  -\int_\Omega V(u) dx. 
\end{equation*}

	Applying DG discretization \eqref{eq-wgd} to \eqref{eq-bbm} yields the following scheme: 
Find $\phi_h,u_h,v_h,w_h,p_h \in V_h$, such that for all $\vphi_i \in V_h$, $i = 1,2,3,4,5$, we have 
\begin{subequations}\label{eq-DG-BBM}
\begin{align}
	-\hf\int_{I_j}( u_h)_t \vphi_1 dx+ \int_{I_j} p_h (\vphi_1)_x dx - (\widehat{p_h}\vphi_1^-)_{j+\hf} + (\widehat{p_h}\vphi_1^+)_{j-\hf} = &\; 0,\\
	\hspace{-1cm}\hf\int_{I_j}( \phi_h - \sigma v_h)_t \vphi_2 dx+ \frac{\sigma}{2}\bigg(\int_{I_j} w_h (\vphi_2)_x dx - (\widehat{w_h}\vphi_2^-)_{j+\hf} + (\widehat{w_h}\vphi_2^+)_{j-\hf}\bigg)
	 =&\; \int_{I_j}(p_h - V'(u_h))\vphi_2 dx,\\
	\frac{\sigma}{2}\int_{I_j}(u_h)_t \vphi_3 dx = &\;	\frac{\sigma}{2}\int_{I_j}w_h \vphi_3 dx ,\\
	\frac{\sigma}{2} \left(-\int_{I_j} u_h (\vphi_4)_x dx + (\widehat{u_h}\vphi_4^-)_{j+\hf} - (\widehat{u_h}\vphi_4^+)_{j-\hf}\right) = &\;	\frac{\sigma}{2}\int_{I_j}v_h \vphi_4 dx ,\\
	-\int_{I_j} \phi_h (\vphi_5)_x dx + (\widehat{\phi_h}\vphi_5^-)_{j+\hf} - (\widehat{\phi_h}\vphi_5^+)_{j-\hf} = 	&\; \int_{I_j}u_h \vphi_5 dx .
\end{align}
\end{subequations}
Here we limit the choices of the matrices $A$, $B$ in the numerical flux $\widehat{K\bz_h}$ defined in \eqref{eq-flux} to be $
A = \left(\begin{matrix}
0& \alpha_0& 0& 0&\alpha_1\\
\alpha_0&0& 0& \frac{\sigma}{2}\alpha_2&0\\
0& 0& 0&0&0\\
0& \frac{\sigma}{2}\alpha_2& 0& 0&0\\
\alpha_1& 0 & 0& 0&0
\end{matrix}\right)$ and $B = 0$. As a result, these numerical fluxes in \eqref{eq-DG-BBM} take the form of 	
\begin{subequations}\label{eq-DG-BBM-fl}
\begin{align}
\widehat{p_h} =&\; \{p_h\} - \alpha_0 [u_h] - \alpha_1[p_h],\quad 
&\;\widehat{w_h} = &\;\{w_h\} - \frac{2}{\sigma}\alpha_0[\phi_h]- \alpha_2 [w_h],\\
\widehat{u_h} =&\; \{u_h\} + \alpha_2 [u_h],\quad
&\;\widehat{\phi_h} =&\; \{\phi_h\} + \alpha_1[\phi_h].
\end{align}
\end{subequations}

Note that $Q = I$ and $\tilde{A} = \left(\begin{matrix}
0& \alpha_0& 0& 0&\alpha_1\\
\alpha_0&0& 0& \frac{\sigma}{2}\alpha_2&0\\
0& 0& 0&0&0\\
0& -\frac{\sigma}{2}\alpha_2& 0& 0&0\\
-\alpha_1& 0 & 0& 0&0
\end{matrix}\right)$. One can apply Corollary \ref{cor-ec} to show the proposed DG method conserves a modified energy, given in Proposition \ref{prop-BBM}.
\begin{PROP} \label{prop-BBM}
	The DG scheme \eqref{eq-DG-BBM}, with numerical fluxes given in \eqref{eq-DG-BBM-fl}, conserves the discrete energy
	\begin{equation*}
	\cE_h = -\frac{1}{6}\int_\Omega u_h^3 dx + \frac{\alpha_0}{2}\sum_j \left([u_h+\phi_h]^2-[u_h-\phi_h]^2\right)_{j+\hf}.
	\end{equation*}
\end{PROP}

For implementation, the scheme \eqref{eq-DG-BBM} and \eqref{eq-DG-BBM-fl} can be rewritten as the following system
\begin{subequations} 
	\begin{align}
	- \hf ( {u_h})_t - {D_{-\alpha_1}} {p_h} +\alpha_0 L u_h&= 0,\\
	( \hf {\phi_h}-\frac{\sigma}{2}{v_h})_t - \frac{\sigma }{2} {D_{-\alpha_2}} {w_h} +\alpha_0 L \phi_h &=  {p_h}-  \Pi(V'(u)),\\
	\frac{\sigma}{2}( {u_h})_t &= \frac{\sigma}{2} w_h,\\
	\frac{ \sigma}{2}{D_{\alpha_2}}  {u_h}&=\frac{\sigma}{2}{v_h} ,\\
	{D_{\alpha_1} }	 {\phi_h}&= {u_h} \label{eq-CH1-scheme-5}.
	\end{align}
\end{subequations}
One can eliminate $v_h$, $w_h$ and $p_h$ to obtain
\begin{align*}
	\hf\left( D_{-\alpha_{1}}\phi_h +\left(I - 2\sigma D_{-\alpha_{1}}D_0\right) u_h\right)_t =&\; -\alpha_0 D_{-\alpha_{1}}L\phi_h + \alpha_0 L u_h -D_{-\alpha_{1}}\Pi(V'(u_h)),\\
	u_h=&\; {D_{\alpha_1} }	 {\phi_h}.
\end{align*}
Note although $\alpha_2$ appears in the weak formulation, it eventually gets cancelled and does not play a role in the numerical scheme. 

{Case 1} (Central flux): Suppose $\alpha_0 = \alpha_1 = 0$, we can eliminate $\phi_h$ to get 
\begin{equation*}
( {u_h})_t = -(I-\sigma  {D_0}^2)^{-1} {D_0} \Pi(V'(u_h)) .
\end{equation*}
Numerically, we see that this scheme is $(k+1)$th order accurate for even $k$ on uniform mesh, but is only $k$-th order for odd $k$ or on nonuniform mesh.
 
{Case 2} (Generalized alternating flux): Suppose $\alpha_0 = 0$ and $\alpha_1 \neq 0$. We solve for $\phi_h$ from $u_h$ in \eqref{eq-CH1-scheme-5} for further reducing the scheme. The method can then be formally written as 
\begin{equation*}
( {u_h})_t = -\left(\frac{I +  {D_{-\alpha_{1}}} {D_{\alpha_1}^{-1}}}{2} - \sigma  {D_{-\alpha_1}}  {D_0}\right)^{-1}  {D_{-\alpha_1}}\Pi(V'(u_h)).
\end{equation*}
However, the kernel of $ {D_{\alpha_1}}$ is one dimension and  $ {D_{\alpha_1}}^{-1}$ is only well-defined on $\{v_h \in V_h: \int_\Omega v_h dx = 0\}$. Therefore, when implementing the numerical scheme, we do time marching for $u_h - \int_\Omega u_h dx$. We drop the last cell average and perform the matrix inversion on a smaller vector in $\mathbb{R}^{N(k+1)-1}$, and finally retrieve the last cell average using the zero average condition. To have the other matrix inverse $\left(\frac{I +  {D_{-\alpha_{1}}} {D_{\alpha_1}}^{-1}}{2} - \sigma  {D_{-\alpha_1}}  {D_0}\right)^{-1}$  well-defined, we need the polynomial degree $k$ to be even and the number of spatial mesh cells $N$ to be odd. 
Although the implementation of this case can be quite involved, numerically it still suffers order reduction as in Case 1. 

{Case 3}: We then consider $\alpha_0 \neq 0$ and $\alpha_1 = 0$. Then we can obtain
\begin{equation*}
( {u_h})_t = -(I-\sigma  {D_0}^2)^{-1}  \left( {D_0} \Pi(V'(u_h))+\alpha_0 \left(D_0LD_0^{-1} - L\right)u_h\right).
\end{equation*}
This scheme seems to retrieve optimal convergence rate for $u_h$ numerically. 

\subsection{CH equation}
The nonlinear CH equation, a bi-Hamiltonian model for waves in the shallow water, takes the form of 
\begin{equation*}
u_t - u_{xxt}+ 3uu_x-2u_xu_{xx} - uu_{xxx} = 0,
\end{equation*}
and corresponds to the multi-symplectic system \eqref{eq-csv} with $\bz = \left(u\ \phi\ \rho \ v\ w\right)^T$, and
\begin{equation*}
M =  \left(\begin{matrix}
0& \hf&-\hf &0&0\\
-\hf& 0&0 &0&0\\
\hf& 0&0 &0&0\\
0& 0&0 &0&0\\
0& 0&0 &0&0\\
\end{matrix}\right),\quad 
K = \left(\begin{matrix}
0& 0& 0&-1&0\\
0& 0& 0& 0&1\\
0& 0& 0&0&0\\
1& 0& 0& 0&0\\
0& -1& 0& 0&0
\end{matrix}\right),\quad 
S(\bz) = -wu - \frac{u^3}{2} - \frac{u\rho^2}{2}+\rho v.
\end{equation*}	
The CH equation can then be expressed as 
\begin{equation*}
	\left\{
\begin{aligned}
\hf( {\phi} -  {\rho})_t - v_x  =&\; -\left(w+\frac{3}{2} {u}^2 +\hf {\rho^2}\right),\\
-\hf  {u}_t + {w_x} =&\; 0,\\
\hf  {u}_t =&\; -  {u}{\rho}+  {v},\\
 {u_x} =&\;  {\rho},\\
-\phi_x =&\;-  u.
\end{aligned}\right.
\end{equation*}
The associated continuous energy functional is given by
\begin{equation*}
\cE=\int_\Omega E(\bz) dx = -\hf\int_\Omega u(u^2 + \rho^2) dx = -\hf\int_\Omega u(u^2 + u_x^2) dx.
\end{equation*}

The variational form of the DG scheme \eqref{eq-wgd} for the multi-symplectic CH equation is given as follows:
Find $u_h$, $\phi_h$,  $w_h$, $v_h$,  $\rho_h \in V_h$, such that for all test functions $\vphi_{i}\in V_h$, $i = 1,2,3,4,5$, we have
\begin{subequations}\label{eq-DG-CH1}
\begin{align}
 \hf \int_{I_j} \left(\phi_h - \rho_h\right)_t\vphi_{1}dx +  \int_{I_j}  v_h(\vphi_1)_xdx - (\widehat{v_h}\vphi_1^-)_{j+\hf} + (\widehat{v_h}\vphi_1^+)_{j-\hf}  
  =&-\int_{I_j}  \left(w_h+\frac{3}{2}u_h^2 + \frac{1}{2}\rho_h^2\right) \vphi_{1} dx,\\
-\hf \int_{I_j}  (u_h)_t \vphi_2dx - \int_{I_j}  w_h(\vphi_2)_xdx +  \left(\widehat{w_h}\vphi_2^+\right)_{j+\hf} - \left(\widehat{w_h}\vphi_2^-\right)_{j-\hf} =&\;0,\\
\hf \int_{I_j}  (u_h)_t \vphi_{3} dx= \int_{I_j} \left(-u_h\rho_h + v_h\right)&\;\vphi_{3} dx,\\
- \int_{I_j}  u_h (\vphi_4)_x dx +  \left(\widehat{u_h}\vphi_4^+\right)_{j+\hf}-\left(\widehat{u_h}\vphi_4^-\right)_{j-\hf}  =&\;\int_{I_j}  \rho_h\vphi_{4}dx,\label{eq-CH14}\\
\int_{I_j}  \phi_h (\vphi_5)_x dx -\left(\phi_h\vphi_5^+\right)_{j+\hf}+ \left(\phi_h\vphi_5^-\right)_{j-\hf} = &\;-\int_{I_j}  u_h\vphi_{5}dx,\label{eq-CH15}
\end{align}
\end{subequations}
Here we limit the choices of the matrices $A$, $B$ in the numerical flux $\widehat{K\bz_h}$ defined in \eqref{eq-flux} to be 
$A = \left(\begin{matrix}
0& \alpha_0& 0&0&0\\
\alpha_0& 0& 0& 0&0\\
0& 0& 0&0&0\\
0& 0& 0& 0&0\\
0& 0& 0& 0&0
\end{matrix}\right)$ and $B = 0$. As a result, these numerical fluxes in \eqref{eq-DG-CH1} take the form of 
\begin{align}\label{eq-DG-CH1-fl}
\widehat{v_h} = \{v_h\} - \alpha_0[\phi_h],\quad\widehat{w_h} = \{w_h\} + \alpha_0[u_h],\quad
\widehat{u_h} = \{u_h\},\quad
\widehat{\phi_h} = \{\phi_h\}.
\end{align}
One can apply Corollary \ref{cor-ec} to show the proposed DG method conserves a modified energy, as explained below.
\begin{PROP}
	The DG scheme \eqref{eq-DG-CH1}, with numerical fluxes given in \eqref{eq-DG-CH1-fl}, conserves the discrete energy
	\begin{equation*}
	\cE(\bz_h) = -\hf \int_\Omega u_h\left(u_h^2 + \rho_h^2\right) dx + \frac{\alpha_0}{2}\sum_j \left([u_h+\phi_h]^2-[u_h-\phi_h]^2\right)_{j+\hf}. 
	\end{equation*}
\end{PROP}

For implementation, the DG scheme can be written in the following operator form.
	\begin{align*}
	\hf( {\phi_h} -  {\rho_h})_t - {D_{0}} {v_h} + \alpha_0 L \phi_h &= -\Pi\left({w_h}+\frac{3}{2} {u_h}^2 +\hf {\rho_h^2}\right),\\
	-\hf ( {u_h})_t + {D_{0}}  {w_h} + \alpha_0 L u_h&= 0,\\
	\hf ( {u_h})_t &= \Pi\left(-  {u_h}{\rho_h}+  {v_h}\right),\\
	{D_{0}}  {u_h} &=  {\rho_h},\\
	-{D_{0}}  {\phi_h} &=-  {u_h}.
	\end{align*}
After simplification, we have one equation left to update $u_h$
\begin{align*}
\left(u_h\right)_t =&\;  \left(I - D_0^2\right)^{-1}\left(D_{0} \left(D_0\Pi\left(u_h D_0 u_h\right)-\Pi\left(\frac{3}{2}u_h^2 + \hf\left(D_0 u_h\right)^2\right)\right) - \alpha_0\left(D_{0}LD_0^{-1} - L\right)u_h\right).
\end{align*}
%

\subsection{Nonlinear Schr\"odinger equation}
Consider the nonlinear Schr\"odinger equation
\begin{equation}\label{eq-nls}
i u_t + u_{xx} + \alpha|u|^2 u = 0, 
\end{equation}
where $i$ is the imaginary unit and $\alpha$ is a given positive parameter. Let $u = p + iq$. Then \eqref{eq-nls} can be written as the multi-symplectic system \eqref{eq-csv} with $\bz = (p\ q\ v\ w)^T$, and
\begin{equation*}
M =  \left(\begin{matrix}
0& 1&0 &0\\
-1& 0&0 &0\\
0& 0&0 &0\\
0& 0&0 &0\\
\end{matrix}\right),\quad 
K = \left(\begin{matrix}
0& 0&-1&0\\
0& 0& 0&-1\\
1& 0& 0&0\\
0& 1& 0& 0
\end{matrix}\right),\quad 	S(\bz) = \frac{1}{2}\left(v^2+w^2+\frac{\alpha}{2}(p^2+q^2)^2\right).
\end{equation*}
This continuum equation preserves the energy functional 
\begin{equation*}
\cE=\int_\Omega E(\bz) dx =\int_\Omega \frac{\alpha}{2}(p^2+q^2)^2-\hf\left(v^2+w^2\right) dx = \int_\Omega \frac{\alpha}{2}(p^2+q^2)^2-\hf\left(p_x^2+q_x^2\right) dx.
\end{equation*}

We omit the detailed schemes in this section to save space. With alternating fluxes, the DG method \eqref{eq-wgd} applied to \eqref{eq-nls}, is the same as the original LDG method for the NLS equation introduced in \cite{xushuNLS}. Later, the same DG method, with central or alternating fluxes \eqref{eq-wgd}, has been studied in \cite{cai2018local} and \cite{tang2017discontinuous} for its
multi-symplectic property. The energy conservation has also been investigated in these papers, which can also be obtained through Corollary \ref{cor-ec} as below.
\begin{PROP}
	With central or alternating fluxes, the DG scheme \eqref{eq-wgd} for the nonlinear Schr\"odinger equation \eqref{eq-nls} conserves the discrete energy 
	\begin{equation*}
	\begin{aligned}
	\cE_h 
	=&\int_\Omega\frac{\alpha}{2}(p_h^2+q_h^2)^2-\hf\left(v_h^2+w_h^2\right) dx.
	\end{aligned}
	\end{equation*}
\end{PROP}

As has been pointed out in \cite{cai2018local}, this DG scheme also conserves the total charge 
$\cC_h = \int_\Omega p_h^2 + q_h^2  dx$. 
This conservation property relies on special structures of \eqref{eq-nls}, and is hence not covered under the current framework.

\subsection{BBM--KdV equation}
\begin{equation}\label{eq-BBM-KdV}
u_t - \sigma u_{xxt} =  V'(u)_x +  \nu u_{xxx}
\end{equation} 
has the multi-symplectic form \eqref{eq-csv} with $\bz = \left(u\ \theta \ \phi \ w\ \rho\ v \right)^T$, and
\begin{equation*}
\begin{aligned}
M &=  \left(\begin{matrix}
0& \frac{\sigma}{2} &-\hf &0&0&0\\
-\frac{\sigma}{2}&0&0&0&0&0\\
\hf&0&0&0&0&0\\
0&0&0&0&0&0\\
0&0&0&0&0&0\\
0&0&0&0&0&0
\end{matrix}\right),\quad 
K = \left(\begin{matrix}
0& 0& 0& 0&\frac{\sigma}{2}&\nu\\
0& 0& 0& 0&0&0\\
0& 0& 0& -1&0&0\\
0& 0& 1& 0&0&0\\
-\frac{\sigma}{2}& 0& 0& 0&0&0\\
-\nu & 0& 0& 0&0&0
\end{matrix}\right),
\end{aligned}
\end{equation*}
\begin{equation*}
S(\bz) = uw - V(u) -\frac{\nu}{2}v^2-\frac{\sigma}{2}\theta\rho.
\end{equation*}
After simplification, the associated energy functional is given by
\begin{equation*}
\cE=\int_\Omega E(\bz) dx = \int_\Omega  -V(u) + \frac{\nu}{2}v^2 dx.
\end{equation*}
To avoid complications in implementation, we consider the DG method with central fluxes for \eqref{eq-BBM-KdV}. The scheme actually retrieves that for KdV equation \eqref{eq-kdv} with $\sigma = 0, V(u) =  -\frac{\eta}{6}u^3$ and that for BBM equation \eqref{eq-bbm} with $\nu = 0$. Furthermore, we have the following energy conservation property. 
\begin{PROP}
	With central fluxes, the DG scheme \eqref{eq-wgd} for the BBM--KdV equation \eqref{eq-BBM-KdV} preserves the discrete energy
	\begin{equation*}
	\cE_h = \int_\Omega- V(u_h) + \frac{\nu}{2}v_h^2 dx. 
	\end{equation*}	
\end{PROP}

\section{Numerical tests}\label{sec-num}
\setcounter{equation}{0}

In this section, we provide some numerical results to demonstrate the behavior of the proposed multi-symplectic and energy conserving DG methods for the Hamiltonian wave equation, the BBM equation and the CH equation. We refer to \cite{cai2018local} for the performance of the DG methods for the KdV and Schr\"odinger equations, when central or alternating fluxes are used. In the accuracy tests, both uniform and nonuniform meshes are considered. For nonuniform meshes in all tests, the mesh size ratios are set as $2:1:2:1:\cdots$, i.e., we have $\Delta x_{2j-1}= 2\Delta x_{2j}$, $j = 1,2,\cdots \frac{N}{2}$ with $N$ being the number of grid points. For simplicity, periodic conditions are used for all tests. Various RK temporal discretizations are used in the tests, with details provided in each test.  

\subsection{Hamiltonian wave equation}
The DG method with alternating numerical fluxes for the second order wave equation has been studied in \cite{xing2013energy}, and here we provide some additional numerical results when other choices of numerical fluxes are considered. 
\begin{examp}[Accuracy test]\label{examp-wave}
	In this example, we examine the accuracy of the DG scheme \eqref{eq-DG-wv} for the wave equation \eqref{eq-wv} with different choices of the numerical fluxes \eqref{eq-DG-wv-fl}. Consider the setup of the problem as
	\begin{subequations}
		\begin{align}
			u_{tt} &= u_{xx},\qquad x\in (0,2\pi), \qquad t> 0,\label{eq-wv-acc}\\
			u(x,0) &= e^{\sin(x)}, \qquad u_t(x,0) = \cos(x) e^{\sin(x)}.
		\end{align}
	\end{subequations}
	The exact solution to the problem is $u(x,t) = e^{\sin(x+t)}$. For DG scheme with $P^k$ elements, we apply $(k+1)$th order RK method to \eqref{eq-wavescheme}
	 for time discretization. The time step is taken as $\Delta t = 0.01\Delta x$. To reveal the convergence rate, we carefully choose the initial condition $u_h$, such that $w_h$ also has the optimal convergence rate. 
	
	We start with examining the accuracy with central fluxes $(\alpha_{11}, \alpha_{13}, \alpha_{33},\beta) = (0,0,0,0)$.  
	From Table \ref{tab-wave-central}, one can see that, for uniform meshes and for both $u_h$ and $w_h$, the scheme is $k$th order accurate if $k$ is odd and $(k+1)$th order if $k$ is even. However, for nonuniform meshes, the scheme degenerates to $k$th order for all cases. These observations are consistent with our understanding when central fluxes are used. 
	
	Then we study the effect of numerical fluxes on accuracy. Only $P^1$ and $P^2$ schemes are presented for simplicity. Tests in Table \ref{tab-wave-000} indicate that the accuracy of $u_h$ can be improved if any of the following hold: $\alpha_{11}\geq 0$, $\alpha_{13}\neq 0$ or $\beta\neq 0$; while optimal rates of $w_h$ are retrieved when $\alpha_{33}< 0$ or $\alpha_{13}\neq0$. 
	
 \begin{table}[h!]
 	\small
	\centering
	\begin{tabular}{c| c | c | c | c | c  | c | c | c | c }
		\hline
		& & \multicolumn{4}{c|}{Uniform mesh}&\multicolumn{4}{c}{Nonuniform mesh}\\
		\hline
		$k$& $N$  &$\|u-u_h\|_{L^2}$ & order &$\|w-w_h\|_{L^2}$  &order&$\|u-u_h\|_{L^2}$ & order &$\|w-w_h\|_{L^2}$  &order\\
		\hline
		\multirow{3}{*}{$ 0$}
		&40&1.0151E-01& -&1.3108E-01& -&1.3225E-01& -&7.0897E-01&- \\
		&80&5.0774E-02& 1.00&6.5221E-02& 1.01&7.8750E-02& 0.75&7.3469E-01&-0.05 \\
		&160&2.5383E-02& 1.00&3.2589E-02& 1.00&5.4303E-02& 0.54&7.5264E-01&-0.03 \\
		\hline
		\multirow{3}{*}{$ 1$}
		&40&2.6510E-02&-&6.9983E-02& -&2.9029E-02& -&8.4319E-02& - \\
		&80&1.3403E-02& 0.98&3.5313E-02& 0.99&1.4693E-02& 0.98&4.2335E-02& 0.99  \\
		&160&6.7004E-03& 1.00&1.7740E-02& 0.99&7.4419E-03& 0.98&2.1432E-02& 0.98  \\
		\hline
		\multirow{3}{*}{$ 2$}
		&40&6.0027E-05& -&4.8934E-04& -&6.4861E-04& -&2.4136E-03& -  \\
		&80&7.6747E-06& 2.97&6.8270E-05& 2.84&1.5607E-04& 2.06&7.3858E-04& 1.71  \\
		&160&9.5034E-07& 3.01&8.2966E-06& 3.04&3.9928E-05& 1.97&1.7212E-04& 2.10  \\
	    \hline
	    \multirow{3}{*}{$ 3$}
	    &40&8.9060E-06& -&1.0554E-04& -&1.5157E-05& -&1.8731E-04& - \\
	    &80&1.0982E-06& 3.02&1.2701E-05& 3.05&1.7993E-06& 3.07&2.5038E-05& 2.90 \\
	    &160&1.3591E-07& 3.01&1.7115E-06& 2.89&2.2965E-07& 2.97&3.2699E-06& 2.94 \\
	    \hline
	    \multirow{3}{*}{$4$}
        &40&2.5719E-08& -&6.7153E-07& -&3.7699E-07& -&5.8904E-06& - \\
        &80&8.7391E-10& 4.88&2.2994E-08& 4.87&2.3248E-08& 4.02&3.8392E-07& 3.94 \\
        &160&2.6916E-11& 5.02&2.8840E-10& 6.32&1.5393E-09& 3.92&2.4844E-08& 3.95 \\
        \hline
	\end{tabular}
	\caption{Accuracy test for Example \ref{examp-wave} with numerical flux $(\alpha_{11},\alpha_{13},\alpha_{33},\beta)=(0,0,0,0)$. }\label{tab-wave-central}
\end{table}
\begin{table}[h!]
	\small
	\centering
	\begin{tabular}{c| c | c | c | c | c  | c | c | c }
		\hline
		 & \multicolumn{4}{c|}{$k = 1$, uniform mesh}&\multicolumn{4}{c}{$k =2$, nonuniform mesh}\\
		\hline
		 $N$  &$\|u-u_h\|_{L^2}$ & order &$\|w-w_h\|_{L^2}$  &order&$\|u-u_h\|_{L^2}$ & order &$\|w-w_h\|_{L^2}$  &order\\
		\hline
		\multicolumn{9}{c}{$(\alpha_{11},\alpha_{13},\alpha_{33},\beta) = (1,0,0,0)$}\\
		\hline
		40&5.2146E-03& -&1.1208E-01& - &1.9070E-04& -&3.8593E-03& -  \\
		80&1.1765E-03& 2.15&6.2546E-02& 0.84&4.4557E-05& 2.10&9.3371E-04& 2.05  \\
		160&5.9084E-04& 0.99&3.0261E-02& 1.05  &4.5986E-06& 3.28&2.9120E-04& 1.68 \\
		320&9.3838E-05& 2.65&1.6099E-02& 0.91&8.2454E-07& 2.48&7.0660E-05& 2.04 \\
		\hline
		\multicolumn{9}{c}{$(\alpha_{11},\alpha_{13},\alpha_{33},\beta) = (0,0,-1,0)$}\\
		\hline
		40&2.5705E-02& -&7.4606E-03& -
		&1.2537E-03& -&3.4249E-04& -  \\
		80&1.3034E-02& 0.98&1.9608E-03& 1.93
		&3.0751E-04& 2.03&4.5086E-05& 2.93 \\
		160&6.6025E-03& 0.98&4.9729E-04& 1.98
		&7.6151E-05& 2.01&6.7623E-06& 2.74 \\
		320&3.3245E-03& 0.99&1.2526E-04& 1.99
		&1.8940E-05& 2.01&7.5680E-07& 3.16 \\
		\hline
		\multicolumn{9}{c}{$(\alpha_{11},\alpha_{13},\alpha_{33},\beta) = (1,0,-1,0)$}\\
		\hline
		40&5.3574E-03& -&4.7856E-03& -
		&4.2658E-04& -&4.7855E-04& - \\
		80&1.0954E-03& 2.29&1.1935E-03& 2.00
		&7.2114E-05& 2.56&6.4808E-05& 2.88 \\
		160&3.1848E-04& 1.78&3.0054E-04& 1.99
		&8.3568E-06& 3.11&1.0014E-05& 2.69 \\
		320&6.5782E-05& 2.28&7.5085E-05& 2.00
		&8.0188E-07& 3.38&7.5943E-07& 3.72 \\
		\hline
		\multicolumn{9}{c}{$(\alpha_{11},\alpha_{13},\alpha_{33},\beta) = (0,\frac{1}{8},0,0)$}\\
		\hline
		40&1.3728E-02& -&2.0183E-02& -
		&2.6276E-04& -&6.8421E-04& -\\
		80&3.3673E-03& 2.03&6.8188E-03& 1.57
		&3.5665E-05& 2.88&1.1415E-04& 2.58 \\
		160&8.9476E-04& 1.91&1.2067E-03& 2.50
		&5.2386E-06& 2.77&1.6382E-05& 2.80 \\
		320&2.2337E-04& 2.00&3.1962E-04& 1.92
		&6.8104E-07& 2.94&1.4917E-06& 3.46 \\
		\hline
		\multicolumn{9}{c}{$(\alpha_{11},\alpha_{13},\alpha_{33},\beta) = (0,0,0,1)$}\\
		\hline
		40&3.4741E-03& -&1.2906E-01& -
		&2.0486E-04& -&7.3941E-03& - \\
		80&8.7692E-04& 1.99&6.5102E-02& 0.99
		&2.5851E-05& 2.99&1.7905E-03& 2.05 \\
		160&2.1881E-04& 2.00&3.2592E-02& 1.00
		&3.5567E-06& 2.86&4.7744E-04& 1.91 \\
		320&5.4648E-05& 2.00&1.6300E-02& 1.00
		&4.1391E-07& 3.10&1.1651E-04& 2.03 \\
		\hline
		\multicolumn{9}{c}{$(\alpha_{11},\alpha_{13},\alpha_{33},\beta) = (0,0,0,-1)$}\\
		\hline
		40&3.4730E-03& -&1.2902E-01& -
		&2.0323E-04& &7.3240E-03& - \\
		80&8.7688E-04& 1.99&6.5099E-02& 0.99
		&2.5744E-05& 2.98&1.7804E-03& 2.04 \\
		160&2.1881E-04& 2.00&3.2592E-02& 1.00
		&3.5506E-06& 2.86&4.7626E-04& 1.90 \\
		320&5.4648E-05& 2.00&1.6300E-02& 1.00
		&4.1351E-07& 3.10&1.1636E-04& 2.03 \\
		\hline
	\end{tabular}
	\caption{Accuracy test for Example \ref{examp-wave} with different numerical fluxes.}\label{tab-wave-000}
\end{table}
\end{examp}
\begin{examp}[Error and energy change]
	In this test, we focus on the numerical error and energy change ($\Delta \cE_h(t):=\cE_h(t) -\cE_h(0)$) after long time simulation. We consider the problem \eqref{eq-wv-acc} together with the initial condition $u(x,0) = \sin(\cos(x))$ and $u_t(x,0) = 0$. The exact solution of the problem is $u(x,t) = \hf \left(\sin(\cos(x+t))+\sin(\cos(x-t))\right)$. $P^3$ elements with $100$ mesh cells are used for numerical simulations. We use the fifth order RK method for time marching and the time step is set as $\Delta t = 0.01\Delta x$ to reduce the temporal error. The final time is set as $T = 200\pi$, which corresponds to $100$ period. The plots for $L^2$ error and the energy change are given in Figure \ref{fig-wv-longtime}. With central fluxes, the scheme has numerical error around $10^{-7}$, and the error is around $10^{-8}$ for other numerical fluxes. It can also be seen that the energy change remains at a small magnitude of around $10^{-13}\sim 10^{-12}$ for all these tests, and the minor energy change seems to be due to time discretization which is not energy conserving. 
	\begin{figure}[h!]
		\centering
		\subfigure[$ (0,0,0,0)$.]{\includegraphics[width=0.23\textwidth]{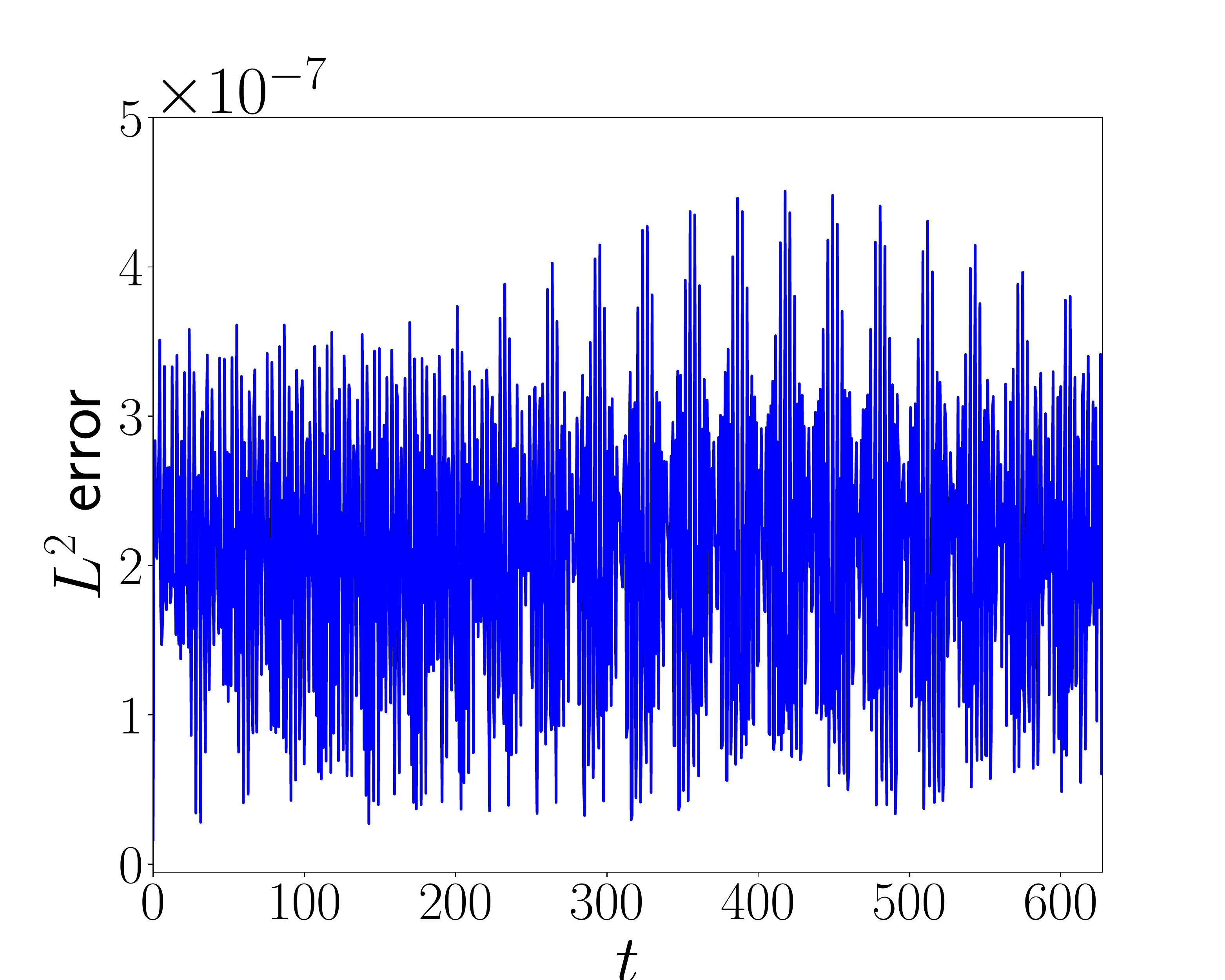}}
		\subfigure[$ (0,\hf,0,0)$.]{\includegraphics[width=0.23\textwidth]{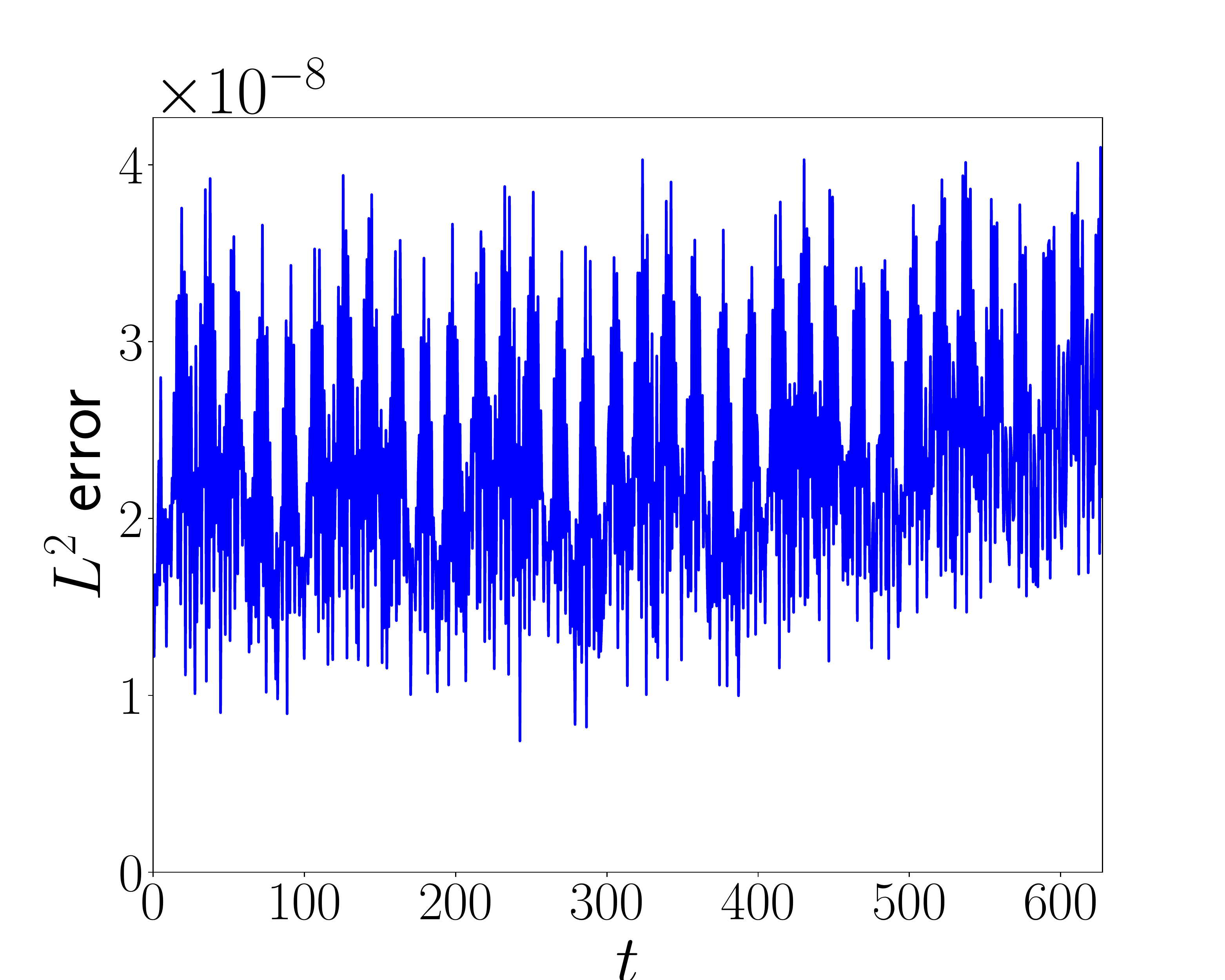}}
		\subfigure[$ (1,0,-1,0)$.]{\includegraphics[width=0.23\textwidth]{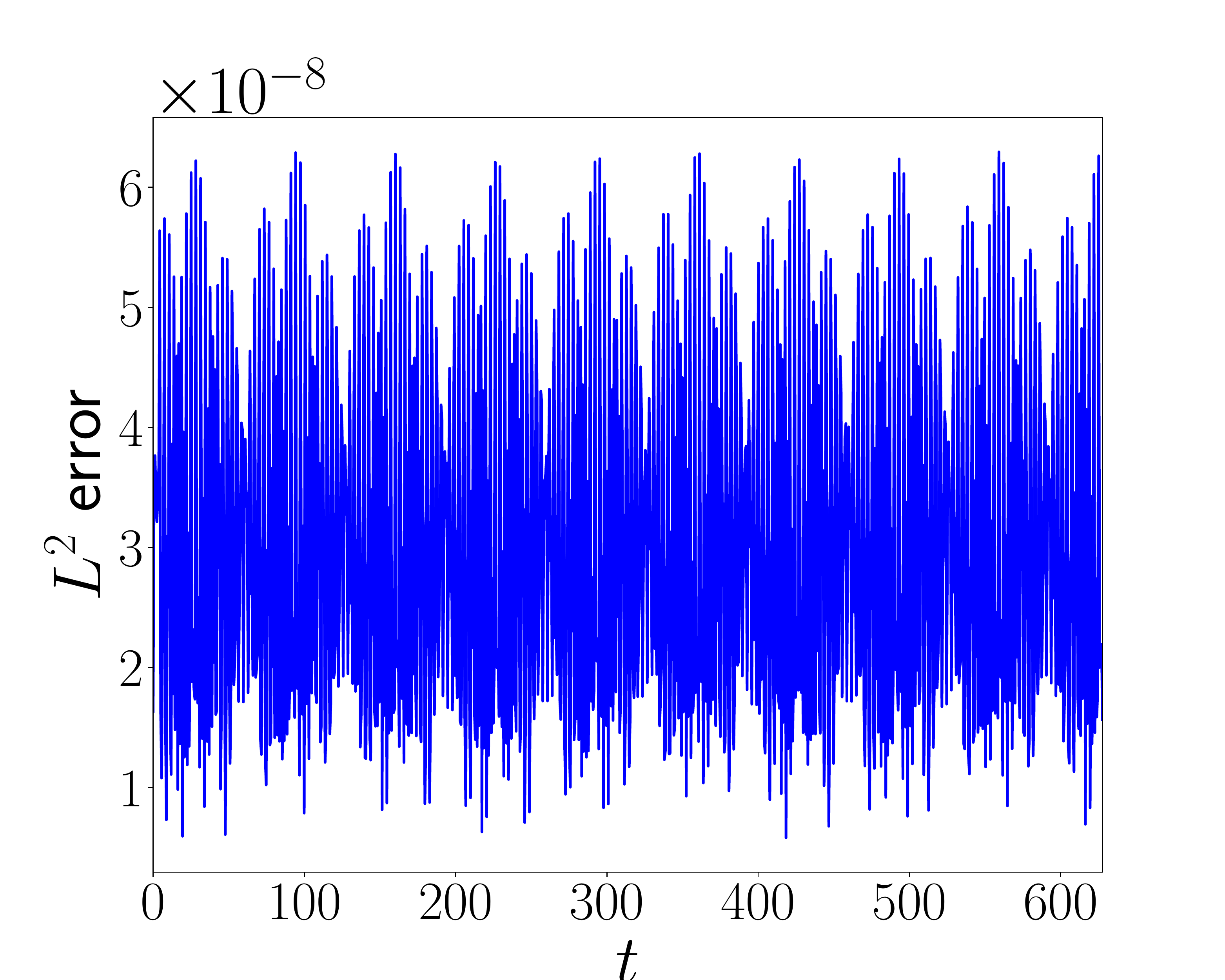}}
		\subfigure[$ (0,0,0,1)$.]{\includegraphics[width=0.23\textwidth]{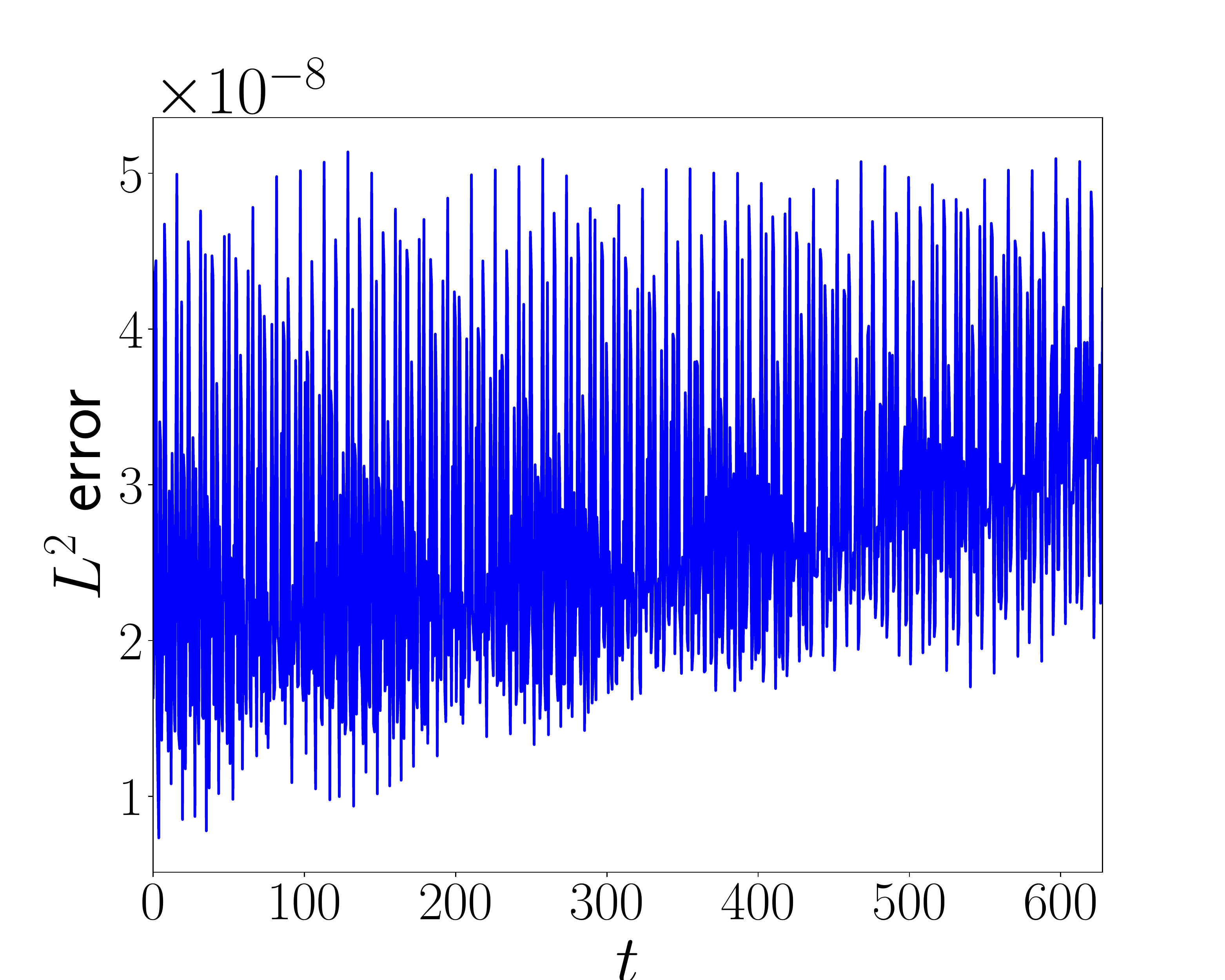}}\\
		\subfigure[$ (0,0,0,0)$.]{\includegraphics[width=0.23\textwidth]{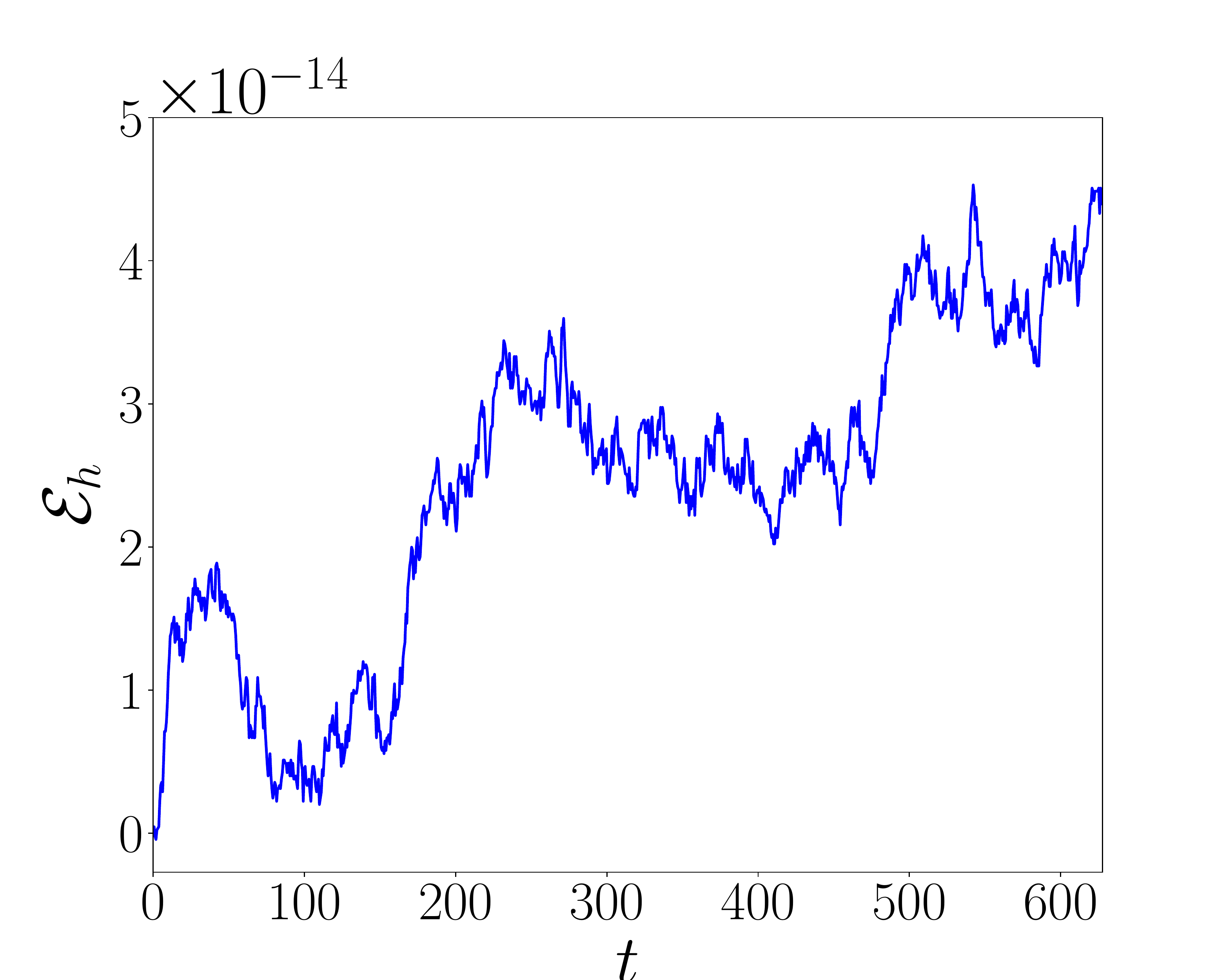}}
		\subfigure[$ (0,\hf,0,0)$.]{\includegraphics[width=0.23\textwidth]{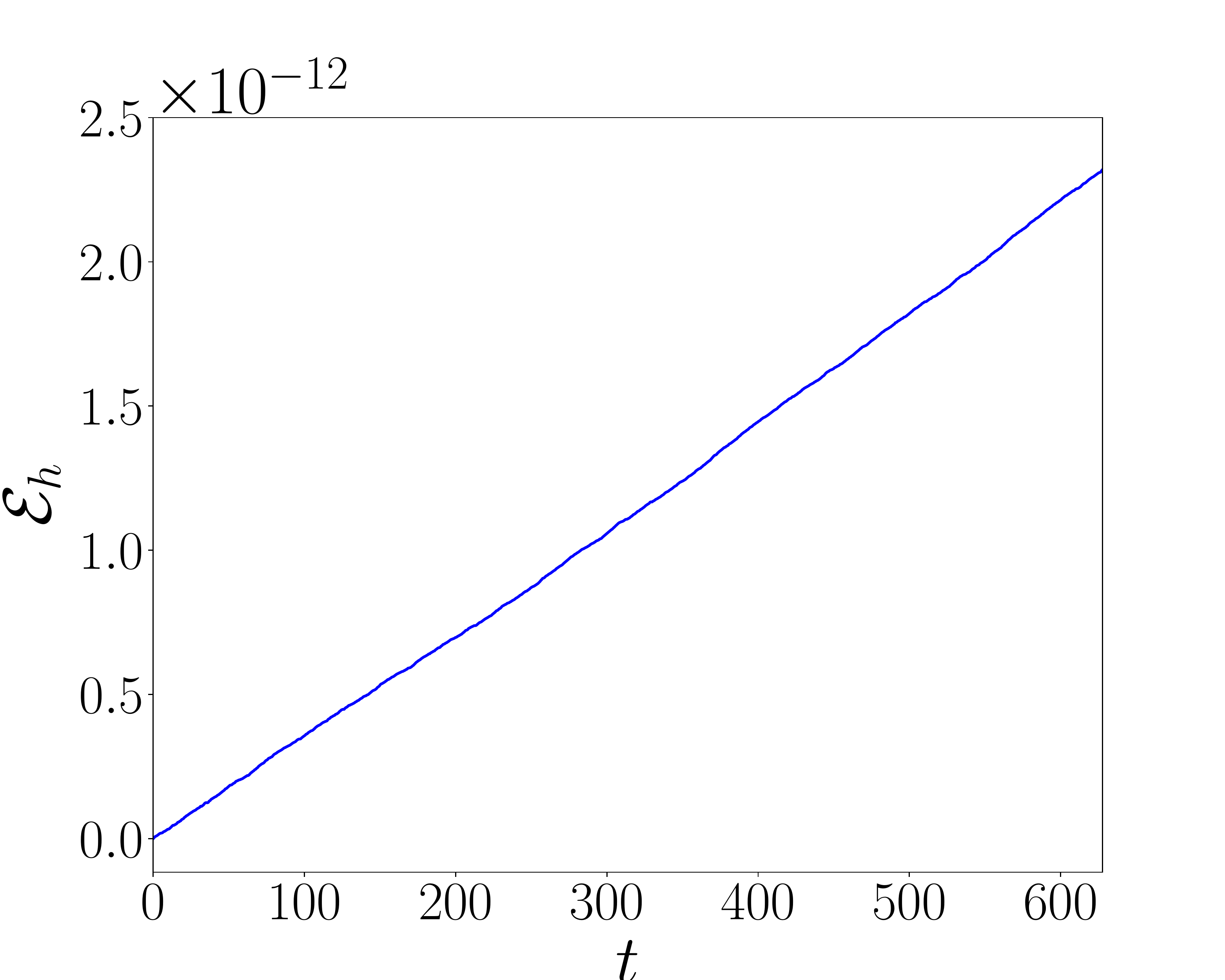}}
		\subfigure[$ (1,0,-1,0)$.]{\includegraphics[width=0.23\textwidth]{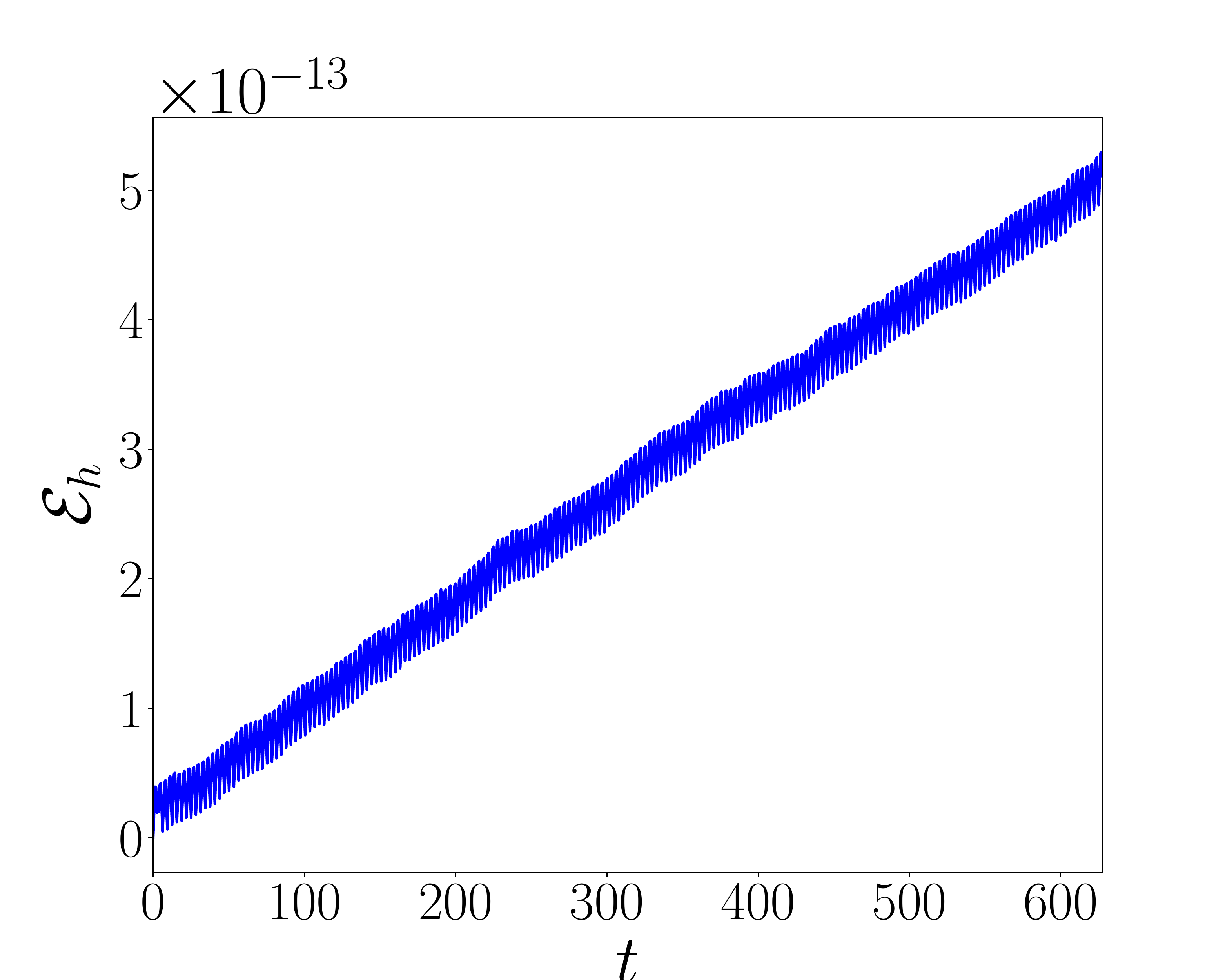}}
		\subfigure[$ (0,0,0,1)$.]{\includegraphics[width=0.23\textwidth]{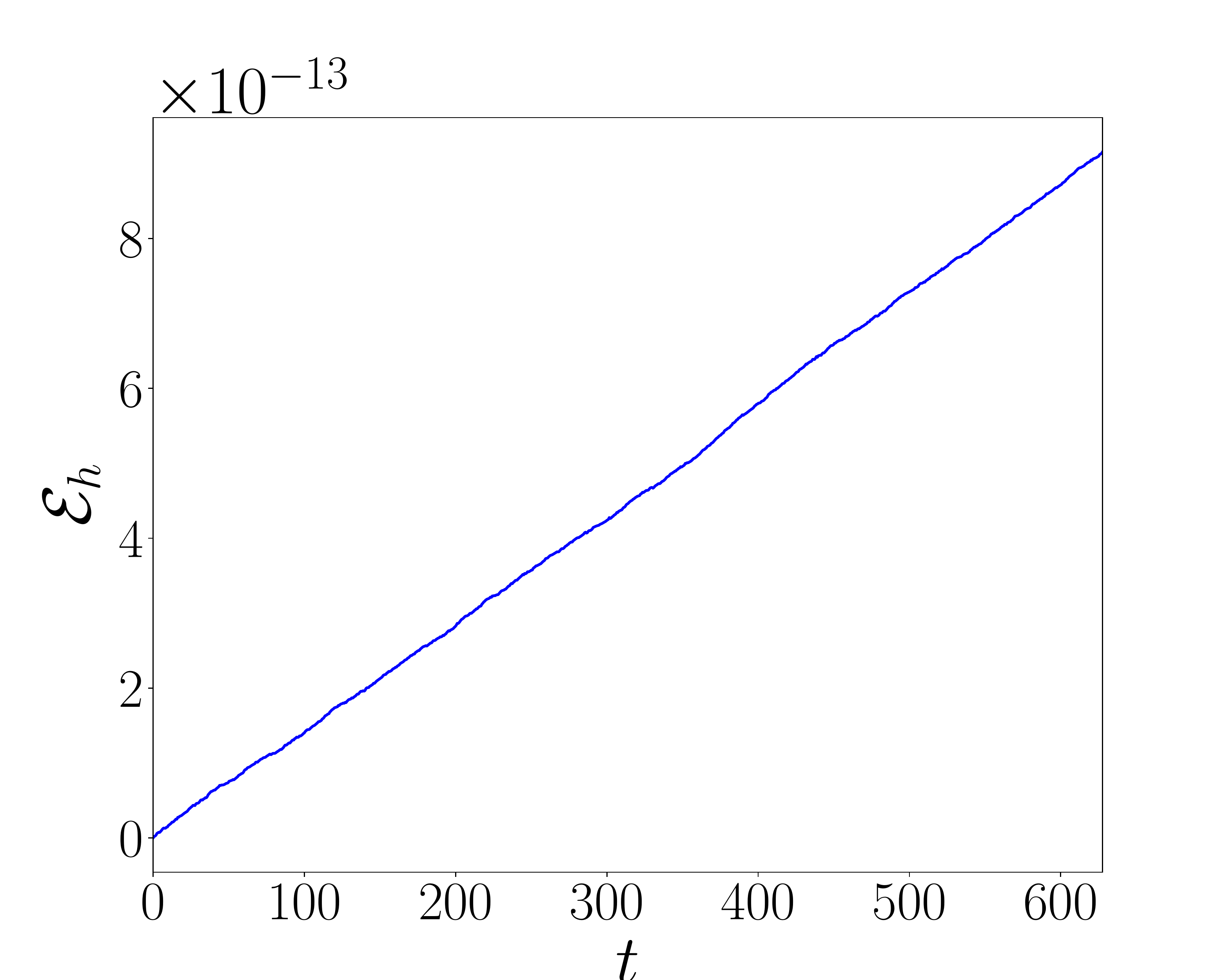}}
		\caption{Numerical error and the energy change for the wave equation with different numerical fluxes. Sub-captions correspond to flux parameters $(\alpha_{11},\alpha_{13},\alpha_{33},\beta)$. The first row: $L^2$ error, the second row: $\Delta \cE_h$. }\label{fig-wv-longtime}
	\end{figure}
\end{examp}

\subsection{BBM equation}\label{examp-BBM-accu}
\begin{examp}[Cnoidal wave]
	The cnoidal wave solution to \eqref{eq-bbm} with $V(u) = -\frac{u^3}{6}$ has the form \begin{equation*}
		u_{\mathrm{c}}(x,t;c,x_0,m) = \frac{3mc}{2m-1}\mathrm{cn}^2\left(\frac{x-ct-x_0}{\sqrt{4(2m-1)\sigma}};m\right),
	\end{equation*}
	where $\mathrm{cn}$ is the Jacobi elliptic function; $c$, $x_0$ and $m$ are given parameters. The solution is periodic on $(0,  {4K(m)}\sqrt{\sigma(2m-1)})$, with $K(\cdot)$ the complete elliptic integral.
	In this numerical example, we consider the setup as in \cite{li2019optimal}.
	\begin{equation*}
	m = 0.9, \quad c = \frac{2m-1}{3m},  \quad \sigma = 10^{-2},  \quad  \Omega =(0,  {4K(m)}\sqrt{\sigma(2m-1)})\approx (0,0.92237).
	\end{equation*}
	
	We start with the accuracy test. Again, for $P^k$ elements, explicit RK method of order $k+1$ is used for time integration. The time step is taken as $\Delta t = 0.5\Delta x$ and we compute to $T = 1$. From Table \ref{tab-BBM-accu}, we observe a similar convergence pattern as that in Table \ref{tab-wave-central} if $\alpha_0 = 0$ (which includes central fluxes). While if we set $\alpha_0 = \frac{1}{4}$ and $\alpha_1 = \alpha_2 = 0$ for numerical fluxes, it achieves the $(k+1)$th order convergence rate for $u$. Other choices of $\alpha_0>0$ seem to retrieve the optimal rate as well. We remark that a large $\alpha_0$ may result in a more restricted time step, and a small $\alpha_0$ may require a more refined mesh to observe the actual convergence rate. When $\alpha_0<0$ is used, the method seems to be unstable numerically.

\begin{table}[h!]
	\centering
	\small
	\begin{tabular}{ c|c | c | c | c | c |c | c  | c |c|c}
		\hline
		&& & \multicolumn{4}{c|}{Uniform mesh}&\multicolumn{4}{c}{Nonuniform mesh}\\
		\hline
		&&&\multicolumn{2}{c|}{$u_h$}&\multicolumn{2}{c|}{$D_0 u_h$}&\multicolumn{2}{c|}{$u_h$}&\multicolumn{2}{c}{$D_0 u_h$}\\
		\hline
		$(\alpha_0,\alpha_1)$&$k$&$N$  &$L^2$ error &order&$L^2$ error &order &$L^2$ error &order&$L^2$ error &order\\
		\hline
		\multirow{8}{*}{$(0,0)$}
			&&40&5.8144E-03& -&1.5710E-01& -&6.4770E-03& -&1.7620E-01& -\\
			&&80&2.8995E-03& 1.00&7.9530E-02& 0.98&3.2078E-03& 1.01&9.1124E-02& 0.95\\
			&$1$&160&1.4514E-03& 1.00&4.0006E-02& 0.99&1.6062E-03& 1.00&4.6065E-02& 0.98 \\
			&&320&7.2588E-04& 1.00&2.0016E-02& 1.00&8.0555E-04& 1.00&2.3327E-02& 0.98 \\
			\cline{2-11}
			&&40&1.2050E-05& -&7.9337E-04& -&1.1517E-04& -&6.9499E-03& - \\
			&&80&1.4591E-06& 3.05&9.7364E-05& 3.03&2.4906E-05& 2.21&1.8067E-03& 1.94 \\
			&$2$&160&1.8812E-07& 2.96&1.2117E-05& 3.01&5.9702E-06& 2.06&4.2151E-04& 2.10\\
			&&320&2.3464E-08& 3.00&1.5153E-06& 3.00&1.5150E-06& 1.98&1.0782E-04& 1.97\\
		\hline
		\multirow{4}{*}{$(0,\hf)$}
		&&41
		&9.8036E-06& -
		&7.1663E-04& -&2.2656E-04& -&9.0495E-03& -\\
		&&81
		&1.5634E-06& 2.65
		&9.6452E-05& 2.89&4.9424E-05& 2.20&2.2454E-03& 2.01\\
		&$2$&161
		&1.9919E-07& 2.97
		&1.2075E-05& 3.00&1.6780E-05& 1.56&5.6368E-04& 1.99\\
		&&321
		&2.5117E-08& 2.99
		&1.5121E-06& 3.00&4.2224E-06& 1.99&1.3827E-04& 2.03\\
		\hline
			\multirow{10}{*}{$(\frac{1}{4},0)$}
			&&40&2.2113E-03& -&1.3945E-01& -&3.3113E-03& -&3.1932E-01& -\\
			&&80&6.7017E-04& 1.72&7.4211E-02& 0.91&9.4065E-04& 1.82&1.5979E-01& 1.00\\
			&$1$&160&1.8814E-04& 1.83&3.8404E-02& 0.95&2.5637E-04& 1.88&7.9702E-02& 1.00\\
			&&320&4.9789E-05& 1.92&1.9481E-02& 0.98&6.6741E-05& 1.94&3.9867E-02& 1.00\\
			&&640&1.2732E-05& 1.97&9.7875E-03& 0.99&1.6952E-05& 1.98&1.9916E-02& 1.00\\
			\cline{2-11}
			&&40&3.7048E-05& -&1.1866E-02& -&7.5531E-05& -&1.5760E-02& -\\
			&&80&4.5232E-06& 3.03&2.8650E-03& 2.05&1.3017E-05& 2.72&4.3587E-03& 2.01\\
			&$2$&160&5.5741E-07& 3.02&7.0477E-04& 2.02&1.9624E-06& 2.73&1.0743E-03& 2.02\\
			&&320&6.9194E-08& 3.01&1.7480E-04& 2.01&2.7868E-07& 2.82&2.6867E-04& 2.00\\
			&&640&8.6198E-09& 3.00&4.3531E-05& 2.01&3.7487E-08& 2.89&6.7173E-05& 2.00\\
		\hline
	\end{tabular}
	\caption{Accuracy test of BBM equation with the cnoidal wave. }\label{tab-BBM-accu}
\end{table}
\end{examp}

\begin{examp}[Energy conserving property]
	We use the same test problem as in the previous example, to compare the performance of $\mathcal{E}_{3}$-conserving scheme (which preserves $\mathcal{E}_{3} =\frac{1}{6}\int_\Omega u^3 dx$) in this paper and the $\mathcal{E}_{2}$-conserving scheme \cite{li2019optimal} (which preserves $\mathcal{E}_{2} =\int_\Omega (u^2+\sigma u_x^2) dx$) for long time simulation. Here $\mathcal{E}_{3}$ and $\mathcal{E}_{2}$ denote the third and second invariant of the BBM equation, respectively. 
	
	$P^2$ polynomials with $10$ mesh cells ($\Delta x \approx 9.2237 \times 10^{-2}$) are used in the test. To reduce the temporal error, we apply the fifth order RK method with a small time step $\Delta t = 0.01\Delta x$ for time discretization.  The $L^2$ error, $\Delta \mathcal{E}_{h,3}$ and $\Delta \mathcal{E}_{h,2}$ are plotted in Figure \ref{fig-BBM-energy}, and the snapshots at $T = 200, 1000, 3000, 5000$ are given in Figure \ref{fig-BBM-profiles}. By comparing the Figure \ref{subfig-BBM-E3-L2err} and Figure \ref{subfig-BBM-E2-L2err}, one can tell the $L^2$ error of the $\mathcal{E}_{3}$-conserving scheme grows slower than that of the $\mathcal{E}_{2}$-conserving scheme. According to Figure \ref{fig-BBM-profiles}, although both methods preserve the profile of the traveling wave well after very long time, the $\mathcal{E}_{3}$-conserving scheme commits slightly smaller phase error compared with the $\mathcal{E}_{2}$-conserving scheme. 
	\begin{figure}[h!]
		\centering
		\subfigure[$L^2$ error.\label{subfig-BBM-E3-L2err}]{		\includegraphics[width=0.3\textwidth]{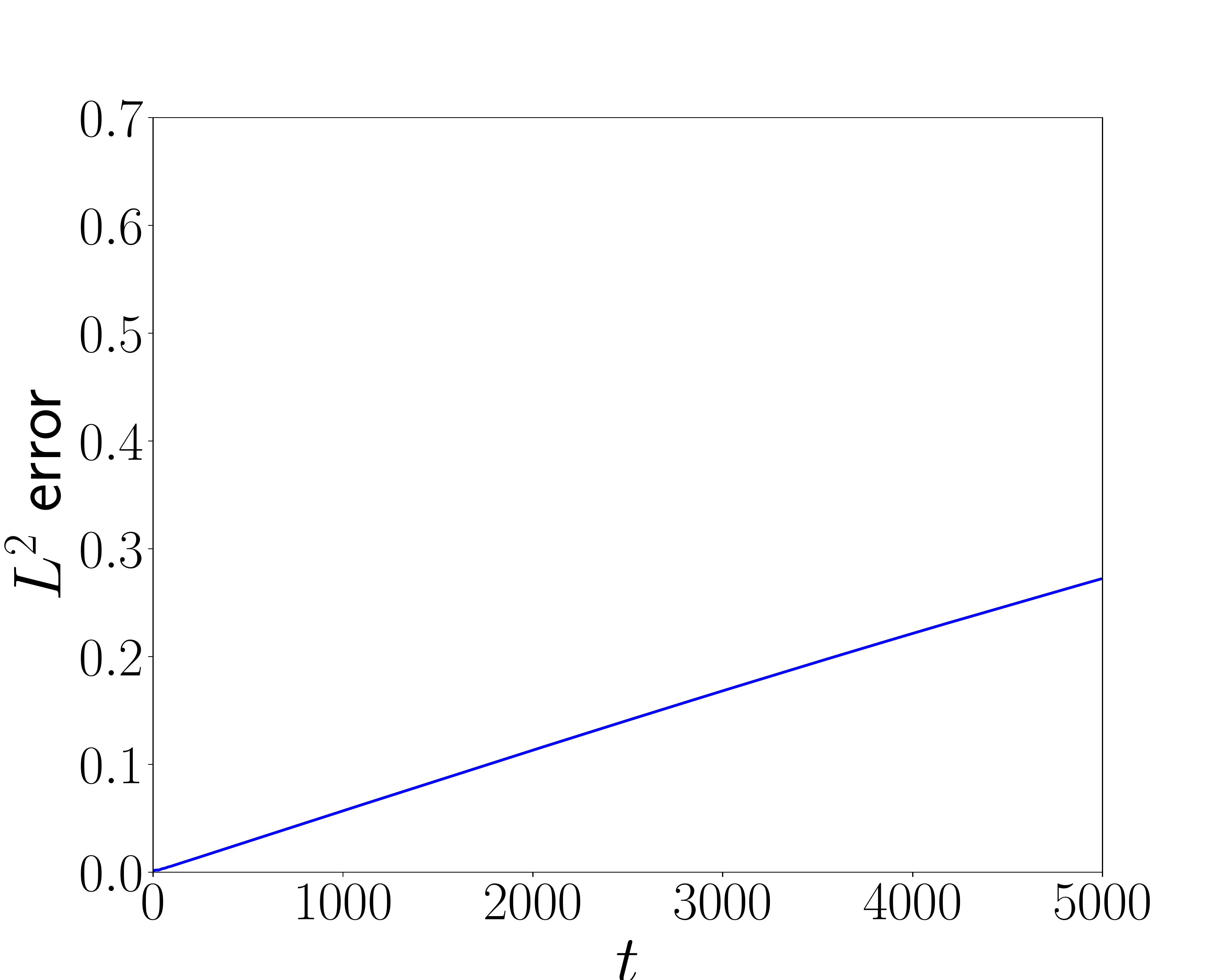}}
		\subfigure[$\mathcal{E}_{h,2}$.]{\includegraphics[width=0.3\textwidth]{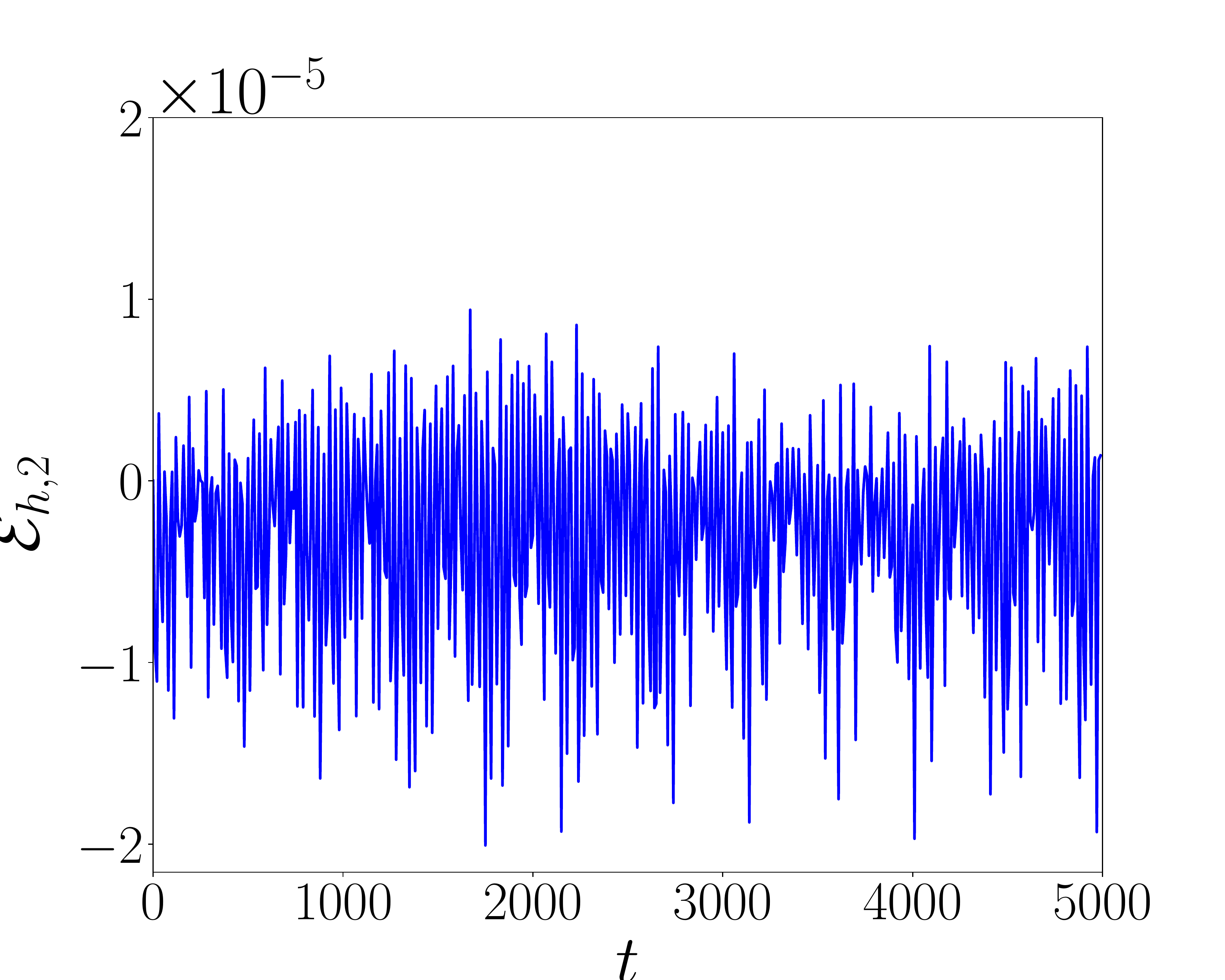}}
		\subfigure[$\mathcal{E}_{h,3}$.]{	\includegraphics[width=0.3\textwidth]{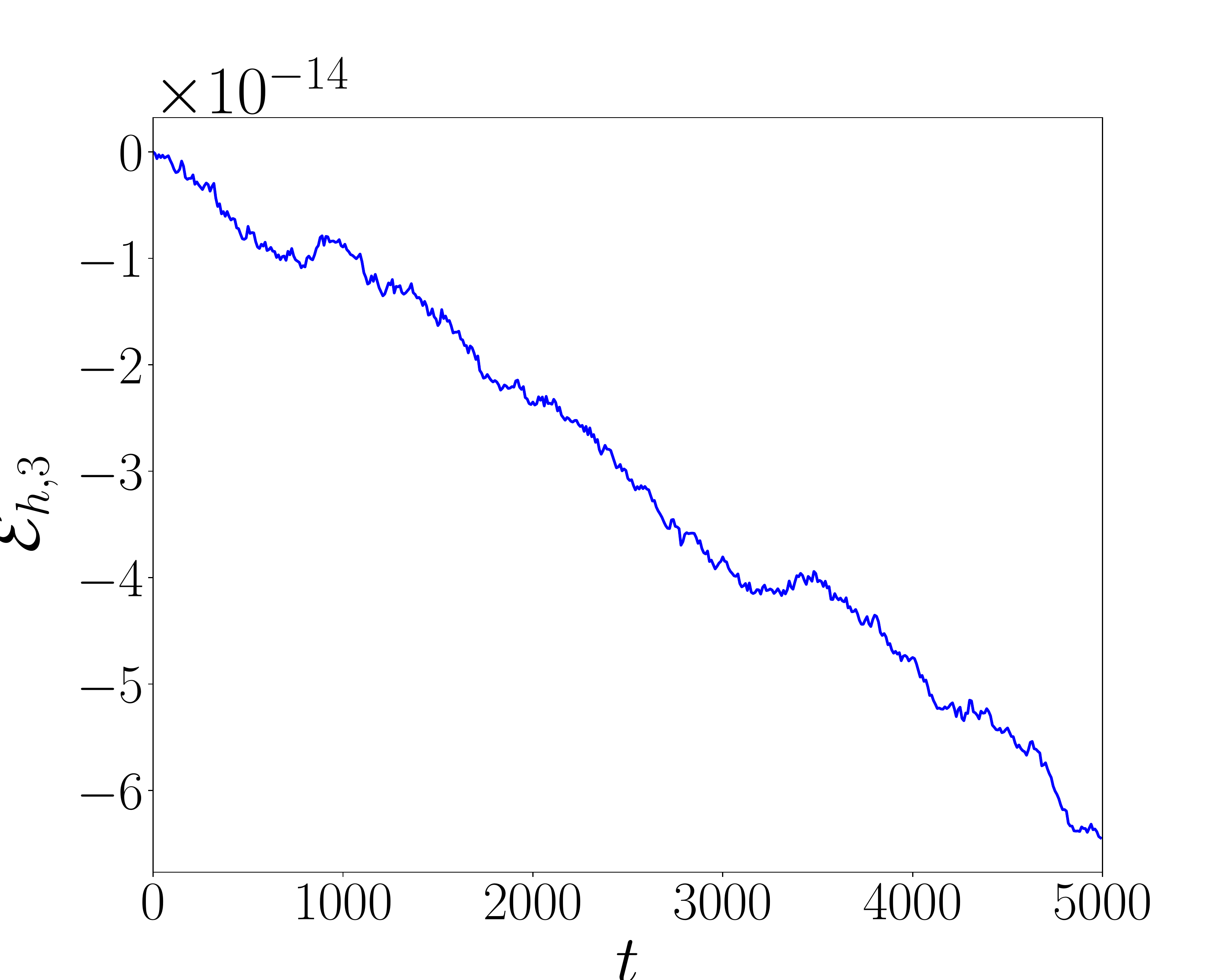}}
		\subfigure[$L^2$ error.\label{subfig-BBM-E2-L2err}]{		\includegraphics[width=0.3\textwidth]{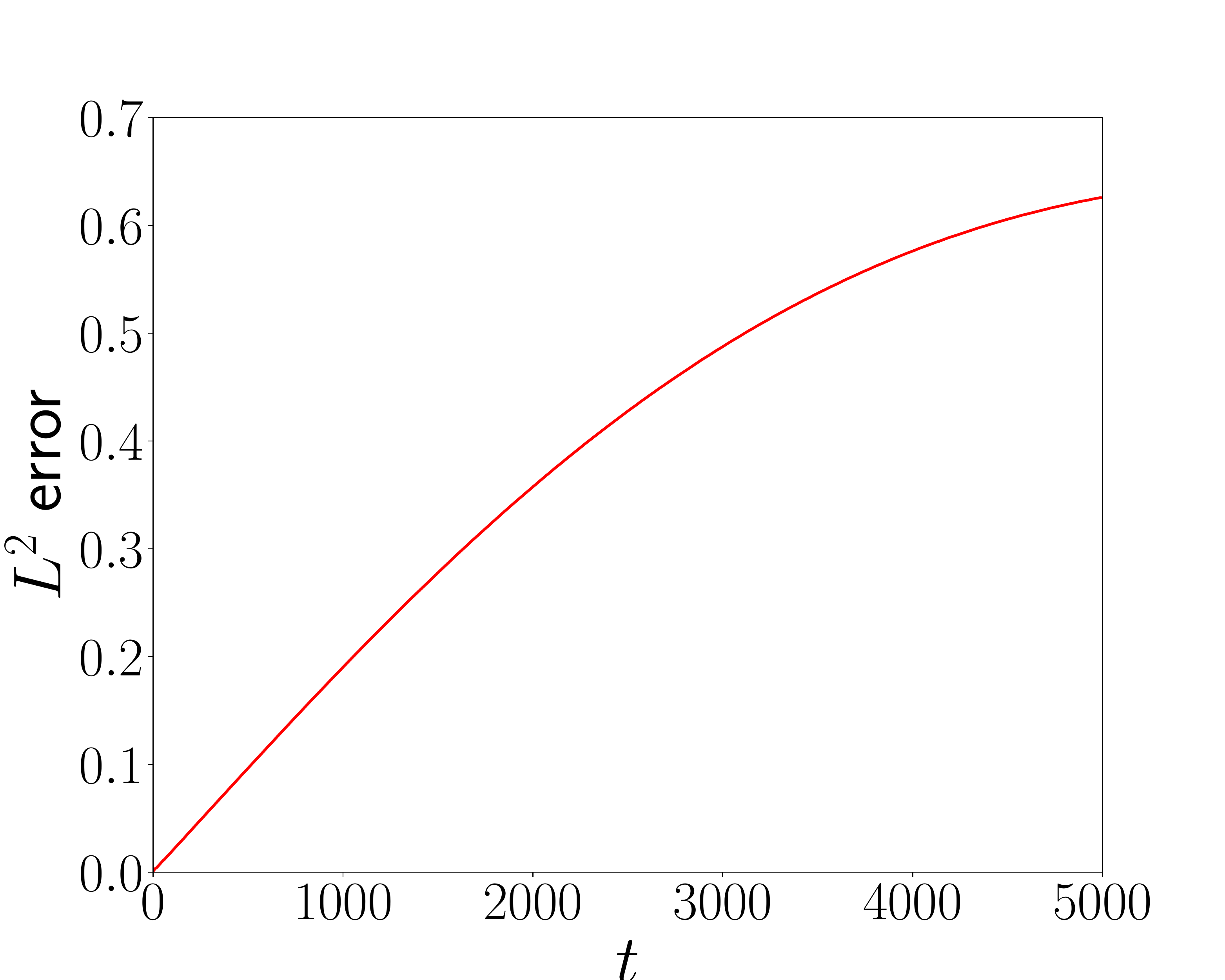}			}
		\subfigure[ $\mathcal{E}_{h,2}$.]{\includegraphics[width=0.3\textwidth]{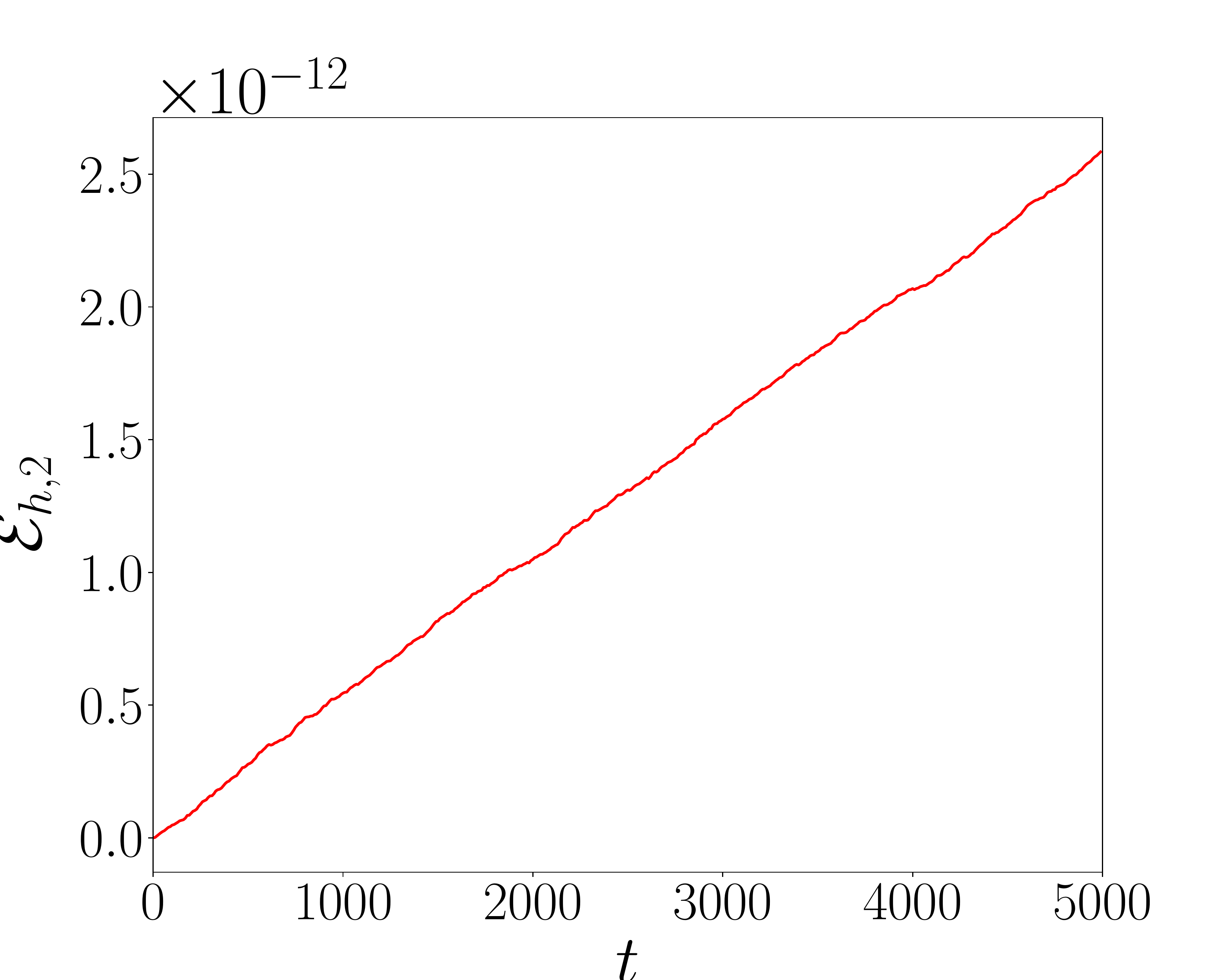}}
		\subfigure[ $\mathcal{E}_{h,3}$.]{\includegraphics[width=0.3\textwidth]{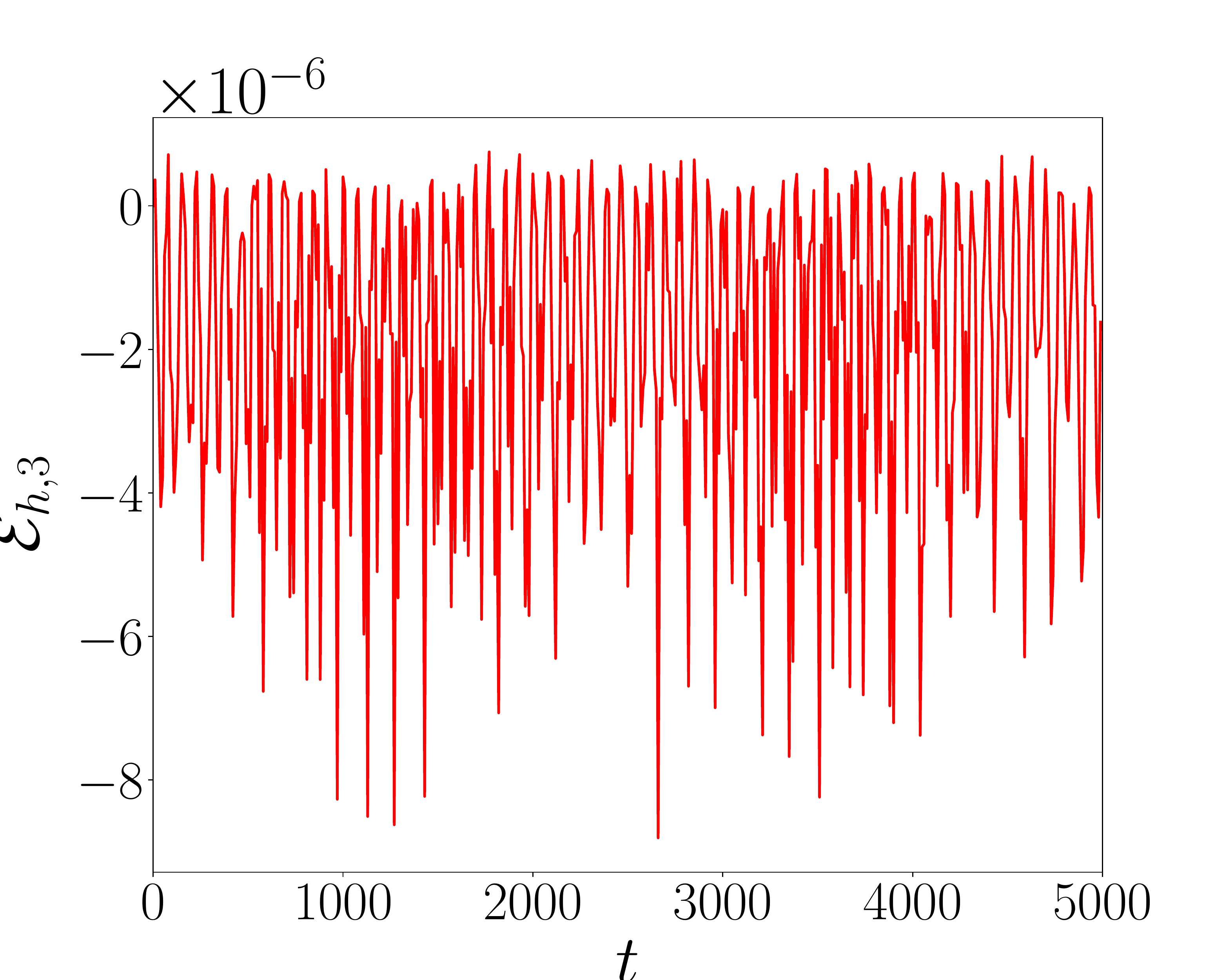}			}
		\caption{Numerical error and the energy change for the BBM equation with different numerical fluxes.  The first row: $\mathcal{E}_{3}$-conserving scheme. The second row: $\mathcal{E}_{2}$-conserving scheme. }\label{fig-BBM-energy}
	\end{figure}
\begin{figure}[h!]
	\centering
	\subfigure[$T = 200$.]{\includegraphics[width=0.24\textwidth]{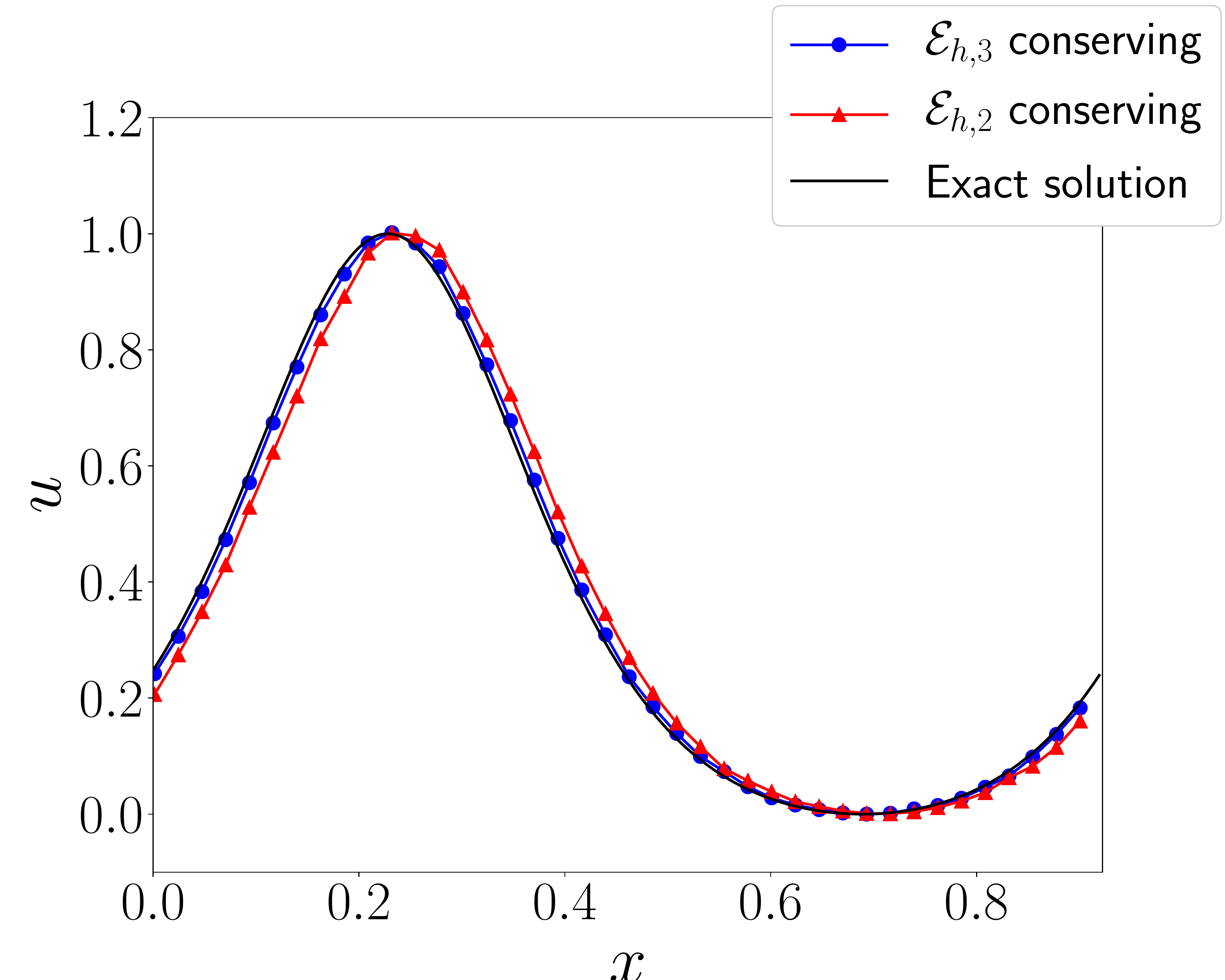}}
	\subfigure[$T = 1000$.]{\includegraphics[width=0.24\textwidth]{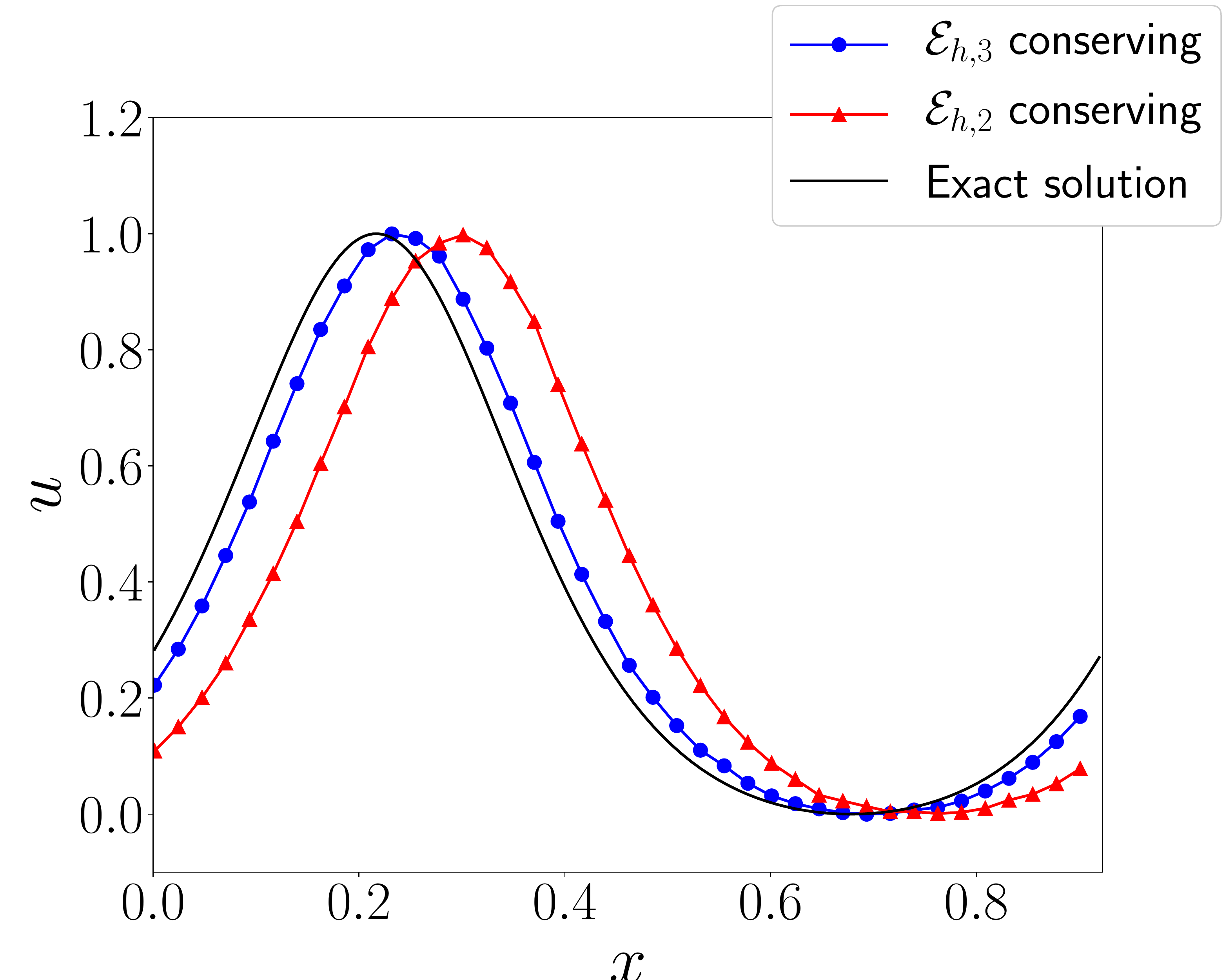}}
	\subfigure[$T = 3000$.]{\includegraphics[width=0.24\textwidth]{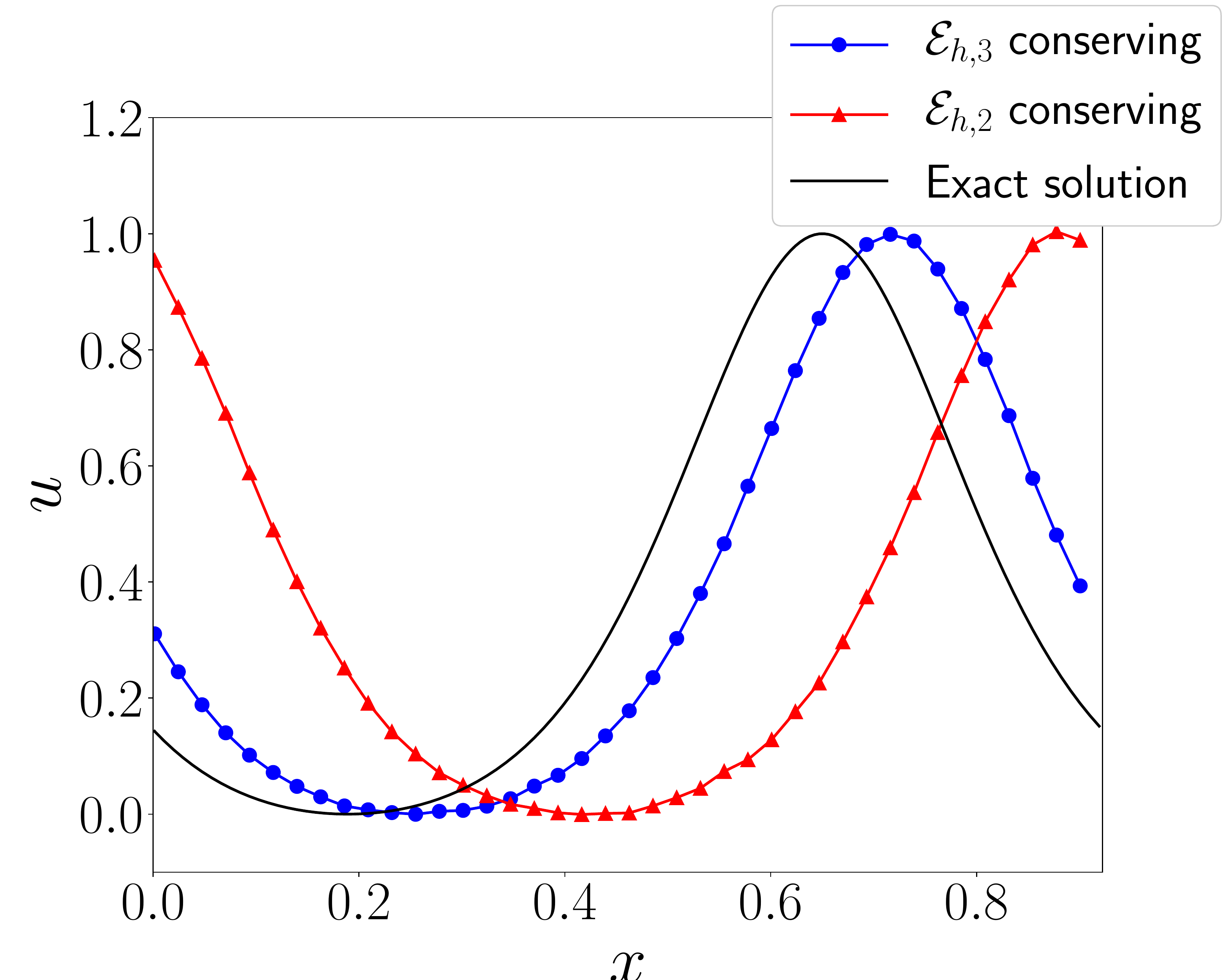}}
	\subfigure[$T = 5000$.]{\includegraphics[width=0.24\textwidth]{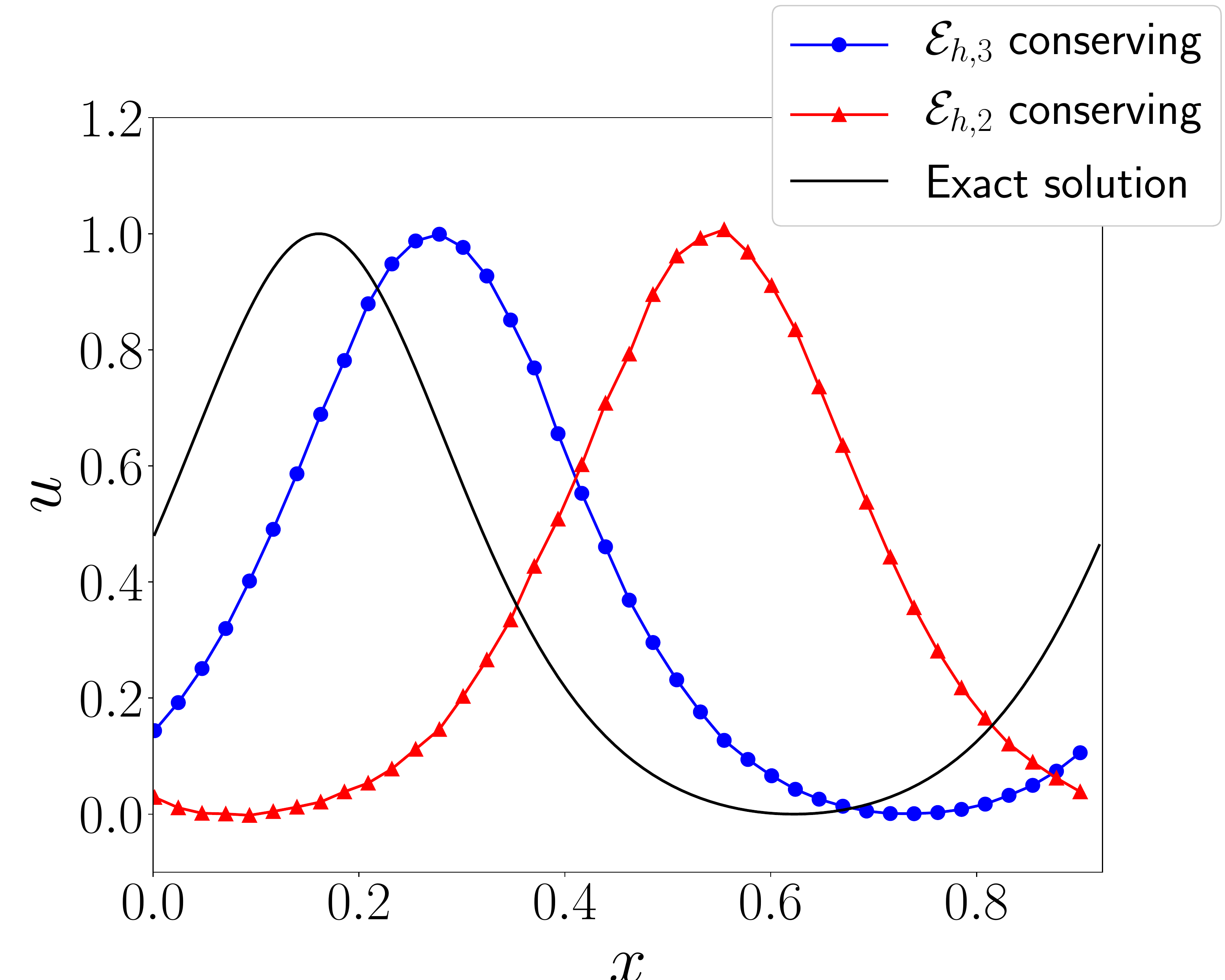}}
	\caption{Solution profiles of the cnoidal waves of the BBM equation at $T = 200, 1000,3000,5000$.}\label{fig-BBM-profiles}
\end{figure}
\end{examp}

\begin{examp}[Solitary wave and soliton interaction]
	Solitary waves of the BBM equation are given as 
	\begin{equation*}
		u_{\mathrm{s}}(x,t; c, x_0) = 3c\ \mathrm{sech}^2\left(\hf \sqrt{\frac{1}{\sigma}}(x-x_0-ct)\right).
	\end{equation*}
	In this example, we perform numerical simulations for the  soliton and multi-soliton interaction of the BBM equation. The initial data are given in Table \ref{tab-BBM-init}. We fix $\sigma = \left(\frac{11}{100}\right)^2$ and use $P^4$ elements with central fluxes for spatial discretization. For time integration, the fifth order RK method with $\Delta t = 0.05 \Delta x$ is used. Settings of the numerical tests and the corresponding results
	are documented in Table \ref{tab-BBM-init}.
	
		\begin{table}[h!]
		\centering
		\small
		\begin{tabular}{c|c|c|c|c}
			\hline
			&Initial condition $u(x,0)$&$T$ & $\Omega$&Results \\
			\hline
			Single soliton&$ u_{\mathrm{s}}(x,0,\frac{1}{5},-2)$& $20$&$(-5,5)$& Figure \ref{fig-BBM-1}\\
			\hline
			Two-soliton& $u_{\mathrm{s}}(x,0,\frac{3}{4}-12) + u_{\mathrm{s}}(x,0,\frac{1}{4},-6)$&$30$&$(-15,15)$&Figure \ref{fig-BBM-2}\\
			\hline
			Four-soliton& $\begin{aligned}&u_{\mathrm{s}}(x,0,\frac{1}{4},-1) + u_{\mathrm{s}}(x,0,\frac{1}{2},-3)\\
			+&u_{\mathrm{s}}(x,0,\frac{3}{4},-5) + u_{\mathrm{s}}(x,0,\frac{5}{4},-13)\end{aligned}$&$20$&$(-15,15)$&Figure \ref{fig-BBM-3}\\
			\hline
		\end{tabular}
		\caption{Setups for BBM equation simulation. $T$: the final time. $\Omega$: the spatial domain. }\label{tab-BBM-init}
	\end{table}
	\begin{figure}[h!]
		\centering
		\subfigure[Solution profile.]{\includegraphics[width=0.38\textwidth]{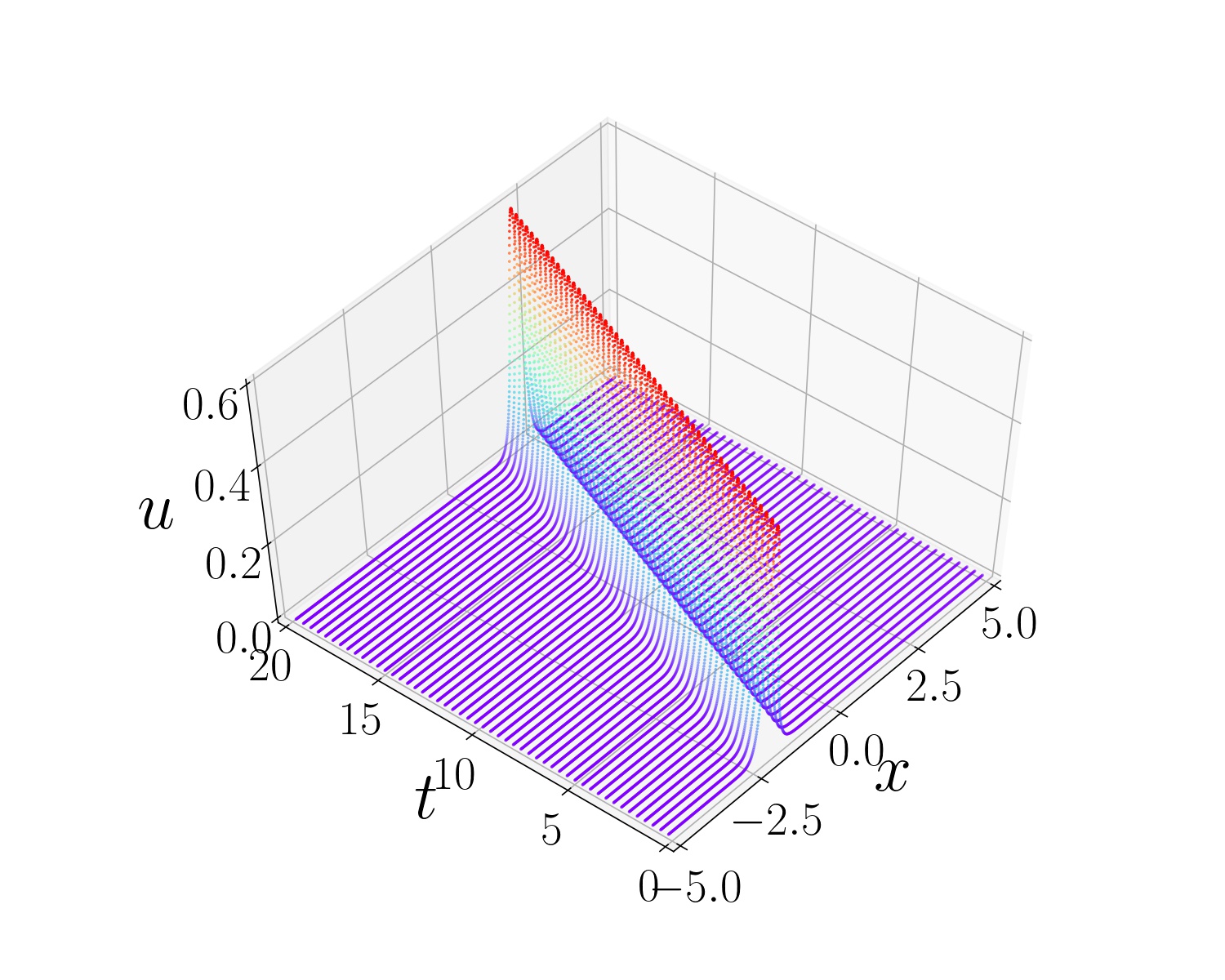}}
		\subfigure[Contour plot.]{\includegraphics[width=0.3\textwidth]{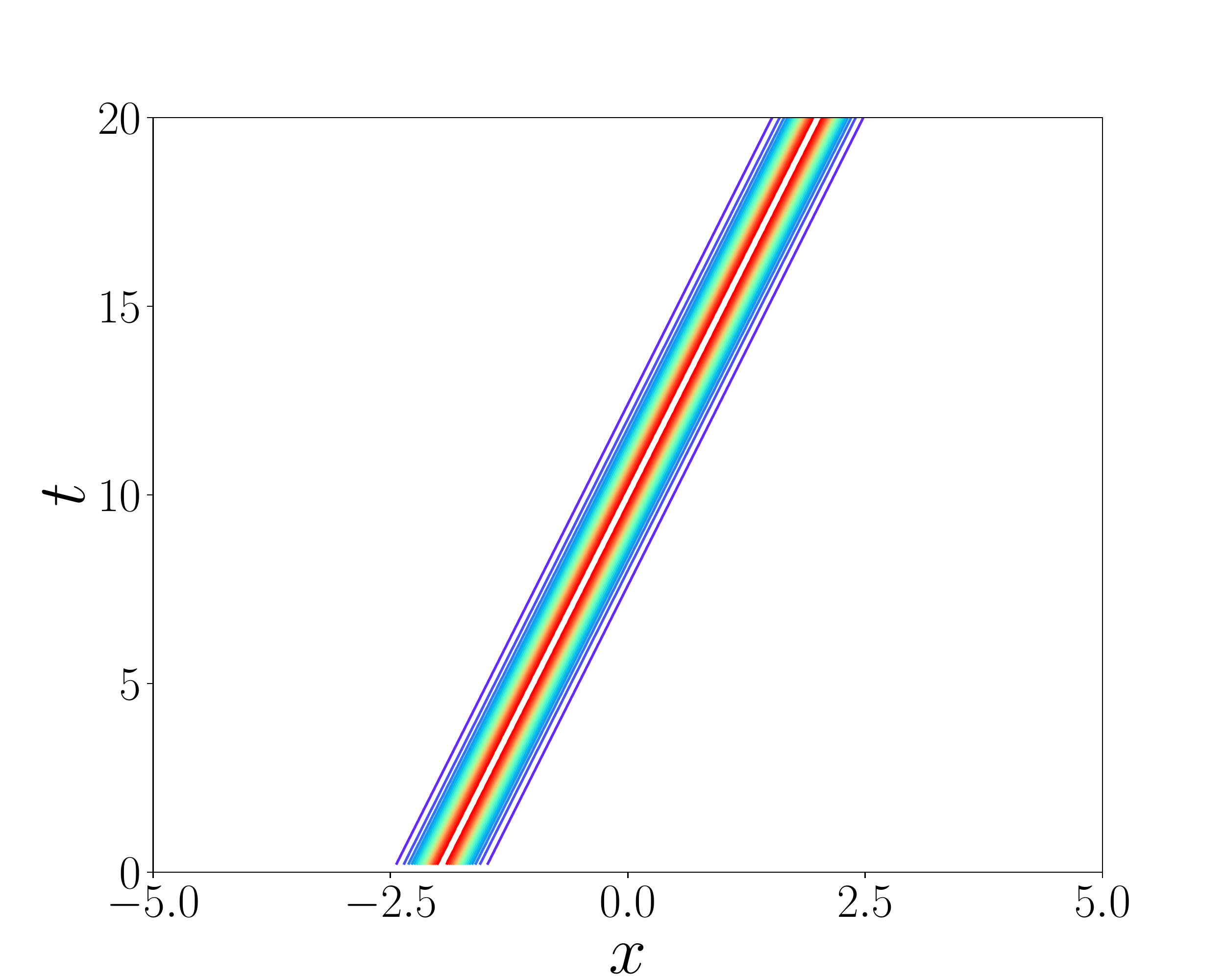}}
		\subfigure[Solution at $T = 20$.]{\includegraphics[width=0.3\textwidth]{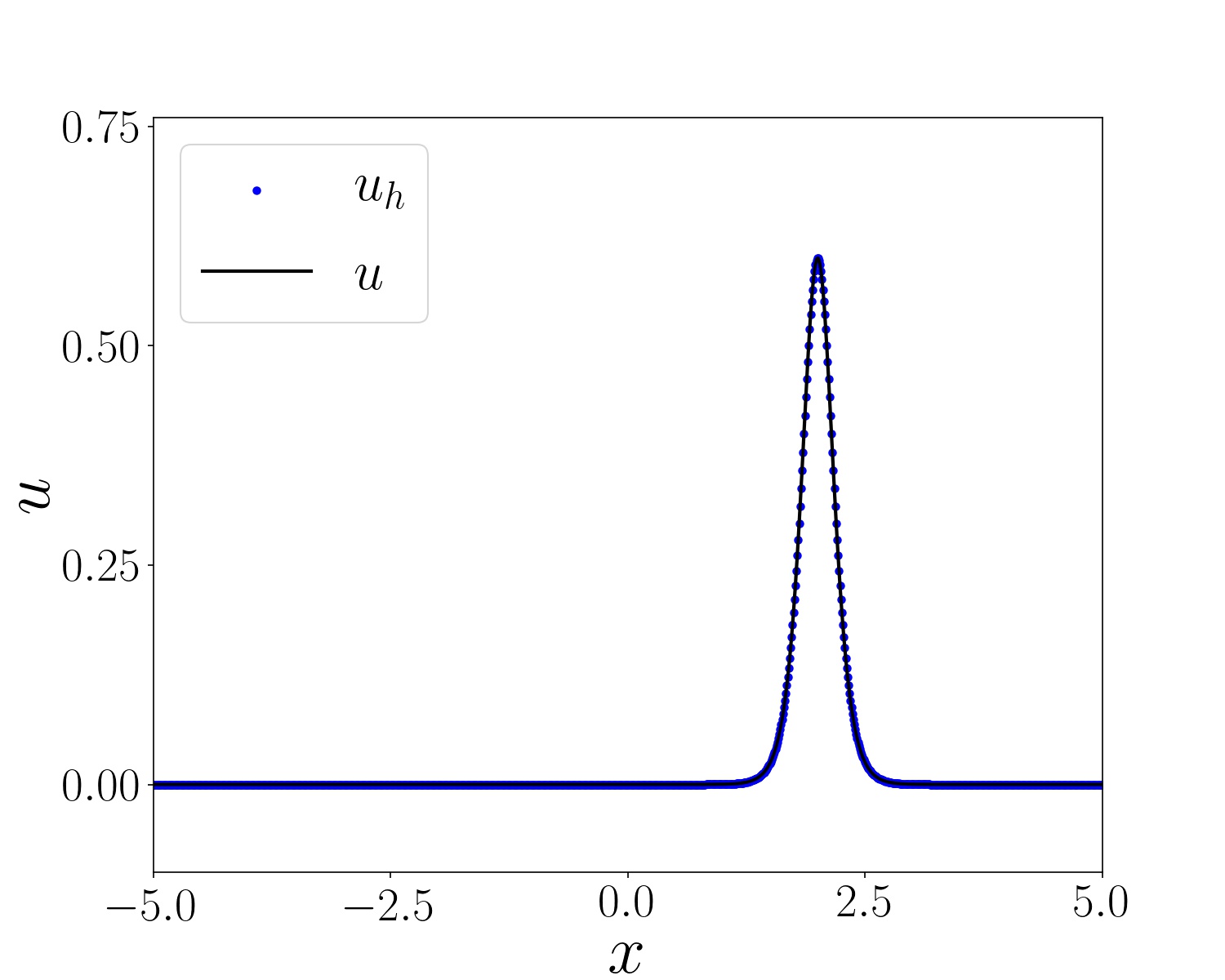}}
		\caption{Single traveling soliton of the BBM equation. }\label{fig-BBM-1}
	\end{figure}
	\begin{figure}[h!]
		\centering
		\subfigure[Solution profile.]{\includegraphics[width=0.38\textwidth]{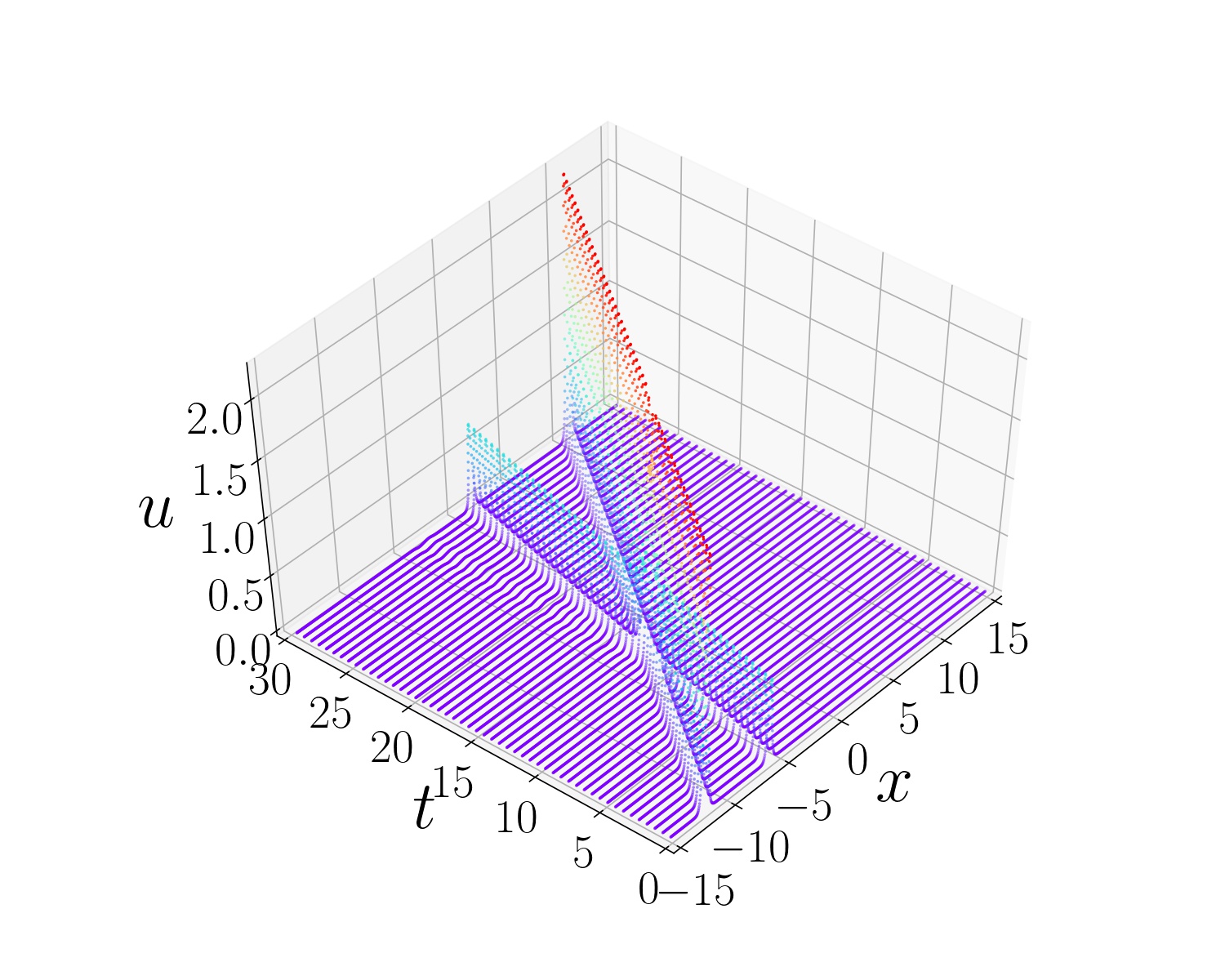}}
		\subfigure[Contour plot.]{\includegraphics[width=0.3\textwidth]{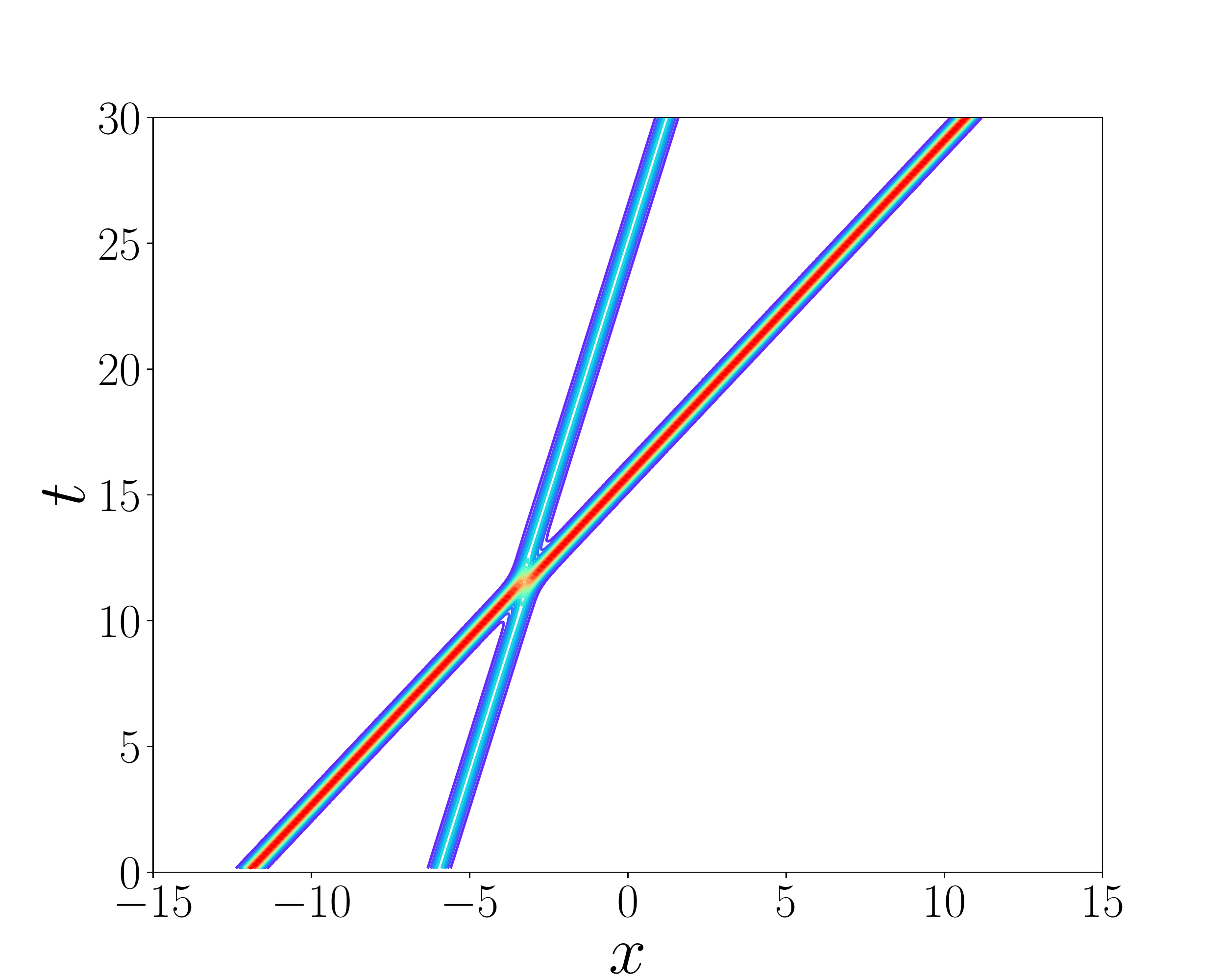}}
		\subfigure[Solution at $T = 30$.]{\includegraphics[width=0.3\textwidth]{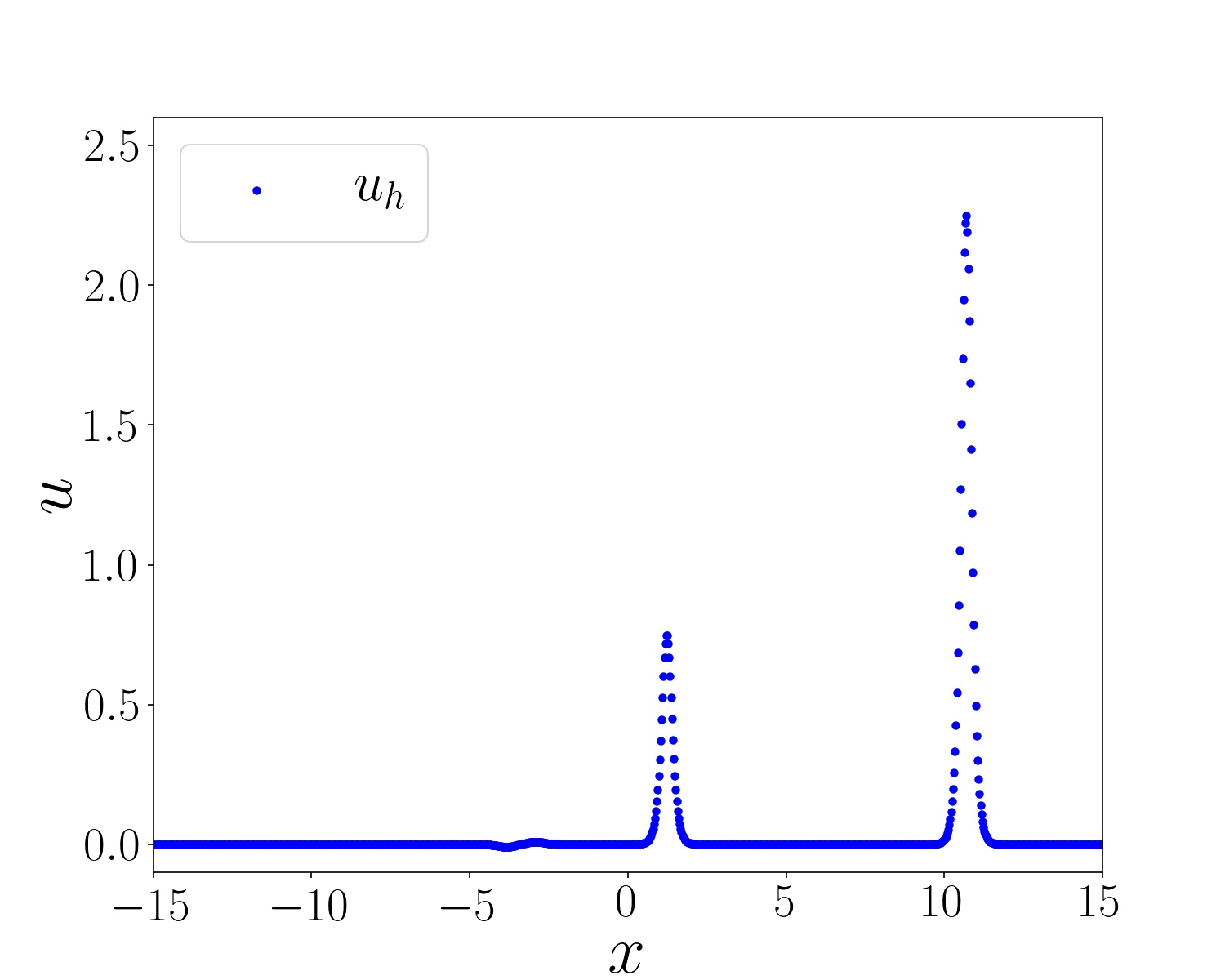}}
		\caption{Two-soliton interaction of the BBM equation.}\label{fig-BBM-2}
	\end{figure}
	\begin{figure}[h!]
	\centering
	\subfigure[Solution profile.]{\includegraphics[width=0.38\textwidth]{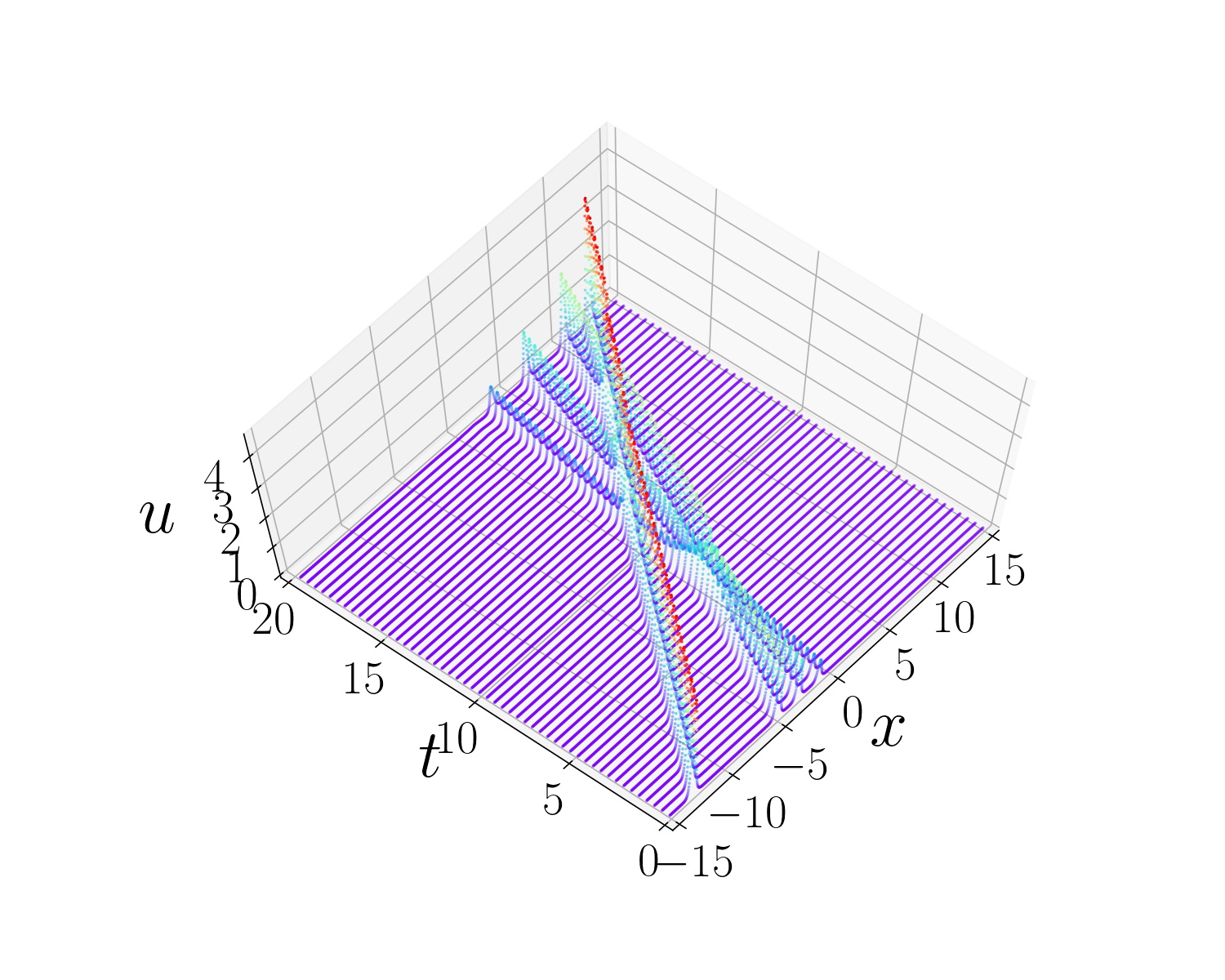}}
	\subfigure[Contour plot.]{\includegraphics[width=0.3\textwidth]{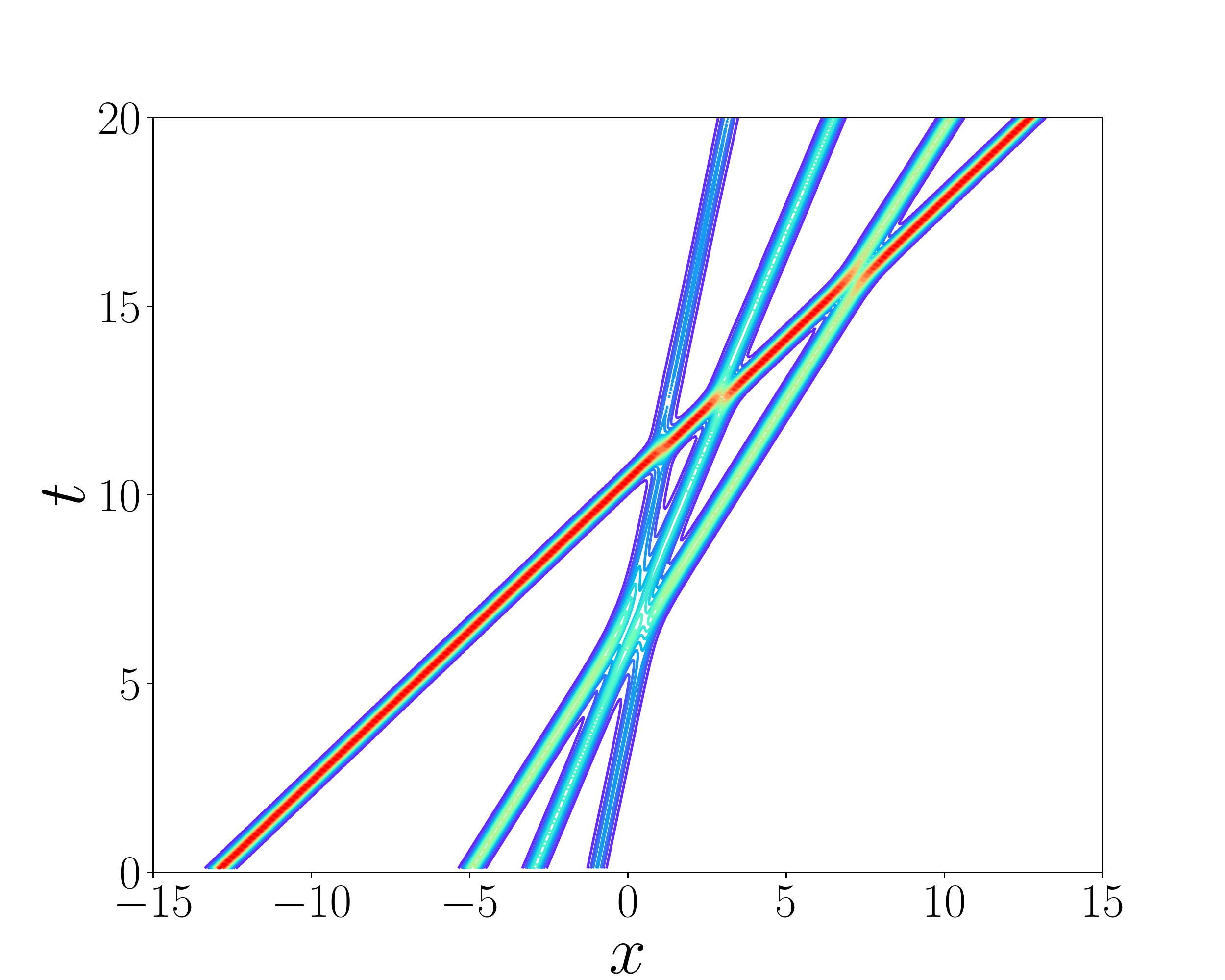}}
	\subfigure[Solution at $T = 20$.]{\includegraphics[width=0.3\textwidth]{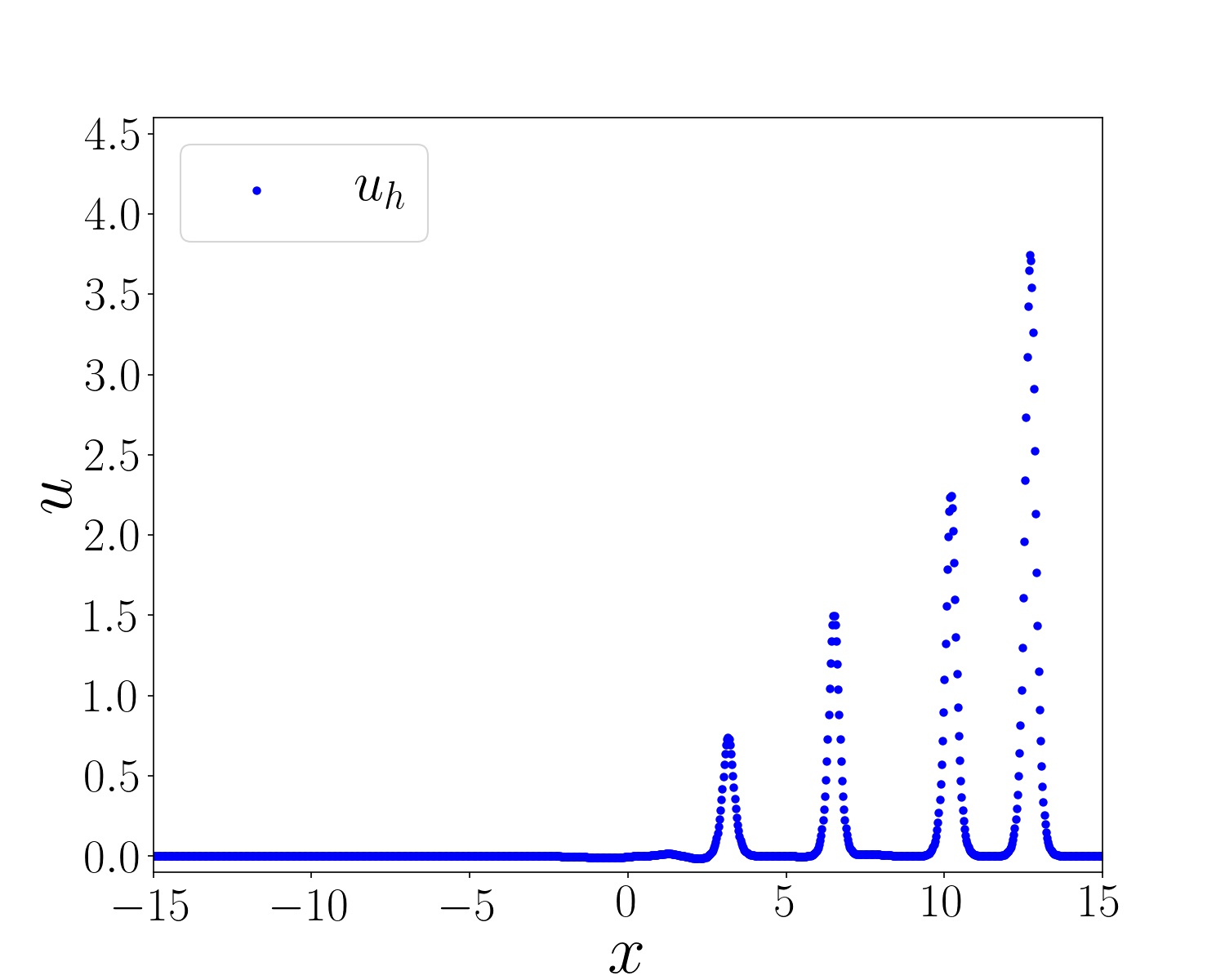}}
	\caption{Four-soliton interaction of the BBM equation.}\label{fig-BBM-3}
	\end{figure}
\end{examp}

\subsection{CH equation}
\begin{examp}[Accuracy test]
	We start with the accuracy test with the fabricated solution. The exact solution is set as $u(x,t) = \sin(x+t)$ and the source term is computed accordingly. We compute to $T = 1$ with the time step $\Delta t = 0.01 \Delta x$ using $(k+1)$th order RK scheme. Numerical fluxes with $\alpha_0 = 0$ and $\alpha_0 = 3$ are tested under this setting. The numerical error is documented in Table \ref{tab-CH-accu}, where the same even-odd phenomenon is observed when $\alpha_0 = 0$ and the optimal convergence rate for $u_h$ seems to be retrieved when $\alpha_0$ is chosen to be $3$.
\begin{table}[h!]
	\centering
	\small
	\begin{tabular}{ c|c | c | c | c | c |c | c  | c |c|c}
		\hline
		&& & \multicolumn{4}{c|}{Uniform mesh}&\multicolumn{4}{c}{Nonuniform mesh}\\
		\hline
		&&&\multicolumn{2}{c|}{$u_h$}&\multicolumn{2}{c|}{$D_0 u_h$}&\multicolumn{2}{c|}{$u_h$}&\multicolumn{2}{c}{$D_0 u_h$}\\
		\hline
		$\alpha_0$&$k$&$N$  &$L^2$ error &order&$L^2$ error &order &$L^2$ error &order&$L^2$ error &order\\
		\hline
		\multirow{8}{*}{$0$}
		&&40&3.3396E-02& -&2.7168E-01& -&3.7844E-02& -&2.9383E-01& -\\
		&&80&1.6468E-02& 1.02&1.5546E-01& 0.81&1.8275E-02& 1.05&1.6763E-01& 0.81\\
		&$1$&160&8.1877E-03& 1.01&8.1319E-02& 0.93&9.1096E-03& 1.00&8.6809E-02& 0.95 \\
		&&320&4.0881E-03& 1.00&4.0392E-02& 1.01&4.5485E-03& 1.00&4.4772E-02& 0.96\\
		\cline{2-11}
		&&40&2.1591E-05& -
			&7.9954E-05& -&5.1452E-04& -&5.6247E-03& -\\
		&&80&2.7672E-06& 2.96
			&9.9825E-06& 3.00&1.2242E-04& 2.07&1.8689E-03& 1.59\\
		&$2$&160&3.4540E-07& 3.00
			&1.2304E-06& 3.02&2.9967E-05& 2.03&6.2347E-04& 1.58\\
		&&320&4.3188E-08& 3.00
			&1.4942E-07& 3.04&7.3536E-06& 2.03&1.7305E-04& 1.85\\
		\hline
		\multirow{8}{*}{$3$}
		&&40& 9.1796E-03& -	  &1.1015E-01& -&1.0887E-02& -&1.2349E-01& -\\
		&&80& 2.5478E-03& 1.85&4.3960E-02& 1.33&2.9429E-03& 1.89&5.0469E-02& 1.29\\
		&$1$&160& 6.6140E-04& 1.95&2.0392E-02& 1.11&7.7927E-04& 1.92&2.3859E-02& 1.08\\
		&&320& 1.6783E-04& 1.98&1.0088E-02& 1.02&1.9728E-04& 1.98&1.1902E-02& 1.00\\
		\cline{2-11}
		&&40&2.1598E-05& -	 &8.0121E-05& -&2.0347E-04& -&4.9410E-03& -\\
		&&80&2.7683E-06& 2.96&9.8237E-06& 3.03&3.2837E-05& 2.63&9.8174E-04& 2.33\\
		&$2$&160&3.4542E-07& 3.00&1.2217E-06& 3.01&4.4606E-06& 2.88&1.7877E-04& 2.46\\
		&&320&4.3193E-08& 3.00&1.5006E-07& 3.03&5.8167E-07& 2.94&3.5117E-05& 2.35\\
		\hline
	\end{tabular}
	\caption{Accuracy test of CH equation with fabricated solutions. }\label{tab-CH-accu}
\end{table}
\end{examp}
\begin{examp}[Traveling peakon and peakon-interaction]
	The periodic peakon solution to CH equation on $(0,x_r)$ is given by 
	\begin{equation}
		u_{\mathrm{p}}(x,t;x_r,c,x_0) = \frac{c}{\cosh(\frac{x_r}{2})}\cosh\left( -\left(x-x_0-ct\right) + x_r \lfloor\frac{x-x_0-ct}{x_r} + \hf\rfloor\right).
	\end{equation}
	In the following numerical tests, we simulate the single traveling peakon, two-peakon interaction, three-peakon interaction and peakon-antipeakon interaction using multi-symplectic DG scheme with central fluxes. Except for the single traveling peakon, other test problems are taken from \cite{liu2016invariant}.  $P^4$ elements with $\Delta x = 0.075$ are used to resolve the solution. The spatial domain is set as $\Omega = (0,30)$. To ensure stability, we use the third order strong-stability-preserving RK method with superviscosity stabilization after each time stage. See \cite{sun2019enforcing} for further studies on this stabilization approach. The time step is set as $\Delta t = 0.01\Delta x$. Other settings of the tests are documented in Table \ref{tab-CH}. The profiles of numerical solutions are consistent with those in \cite{liu2016invariant}.
	
	\begin{table}[h!]
		\centering
		\small
		\begin{tabular}{c|c|c|c}
			\hline
			&Initial condition $u(x,0)$&$T$ & Results \\
			\hline
			Single peakon&$u_{\mathrm{p}}(x,0;30,1,-10)$& $20$& Figure \ref{fig-CH-peakon}\\
			\hline
			Two-peakon& $u_{\mathrm{p}}(x,0;30,2,-5) + u_{\mathrm{p}}(x,0;30,1,5)$&$18$&Figure \ref{fig-CH-2peakon}\\
			\hline
			Three-peakon& $u_{\mathrm{p}}(x,0;30,2,-5)+u_{\mathrm{p}}(x,0;30,1,-3) +u_{\mathrm{p}}(x,0;30,0.8,-1)$&$6$&Figure \ref{fig-CH-3peakon}\\
			\hline
			Peakon-antipeakon&$ u_{\mathrm{p}}(x,0;30,1,-2)+u_{\mathrm{p}}(x,0;30,-1,2)$&$10$&Figure \ref{fig-CH-antipeakon}\\
			\hline
		\end{tabular}
	\caption{Settings for CH equation simulation. $T$: the final time of the simulation. }\label{tab-CH}
	\end{table}
\begin{figure}[h!]
	\centering
	\subfigure[Solution profile.]{\includegraphics[width=0.38\textwidth]{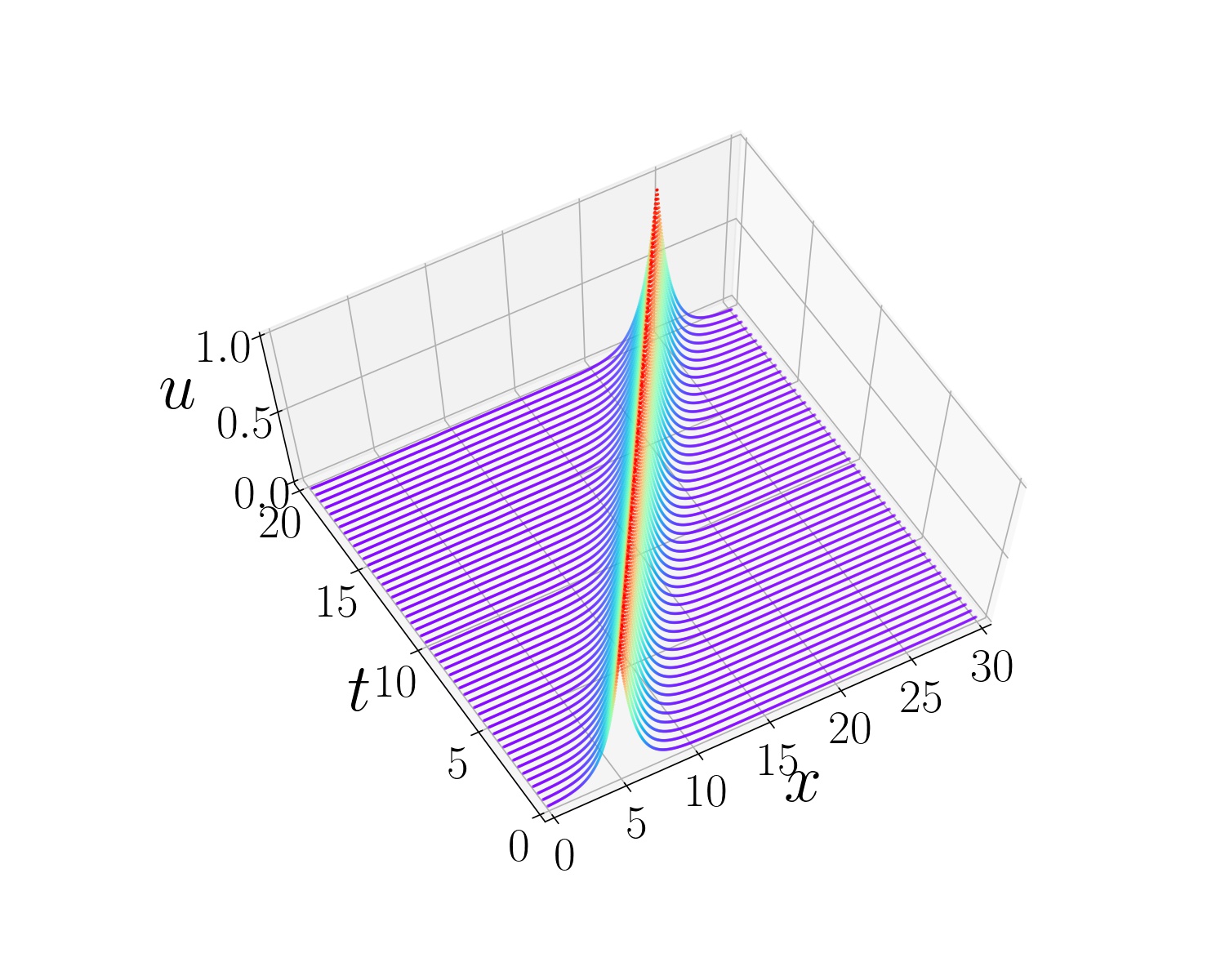}}
	\subfigure[Contour plot.]{\includegraphics[width=0.3\textwidth]{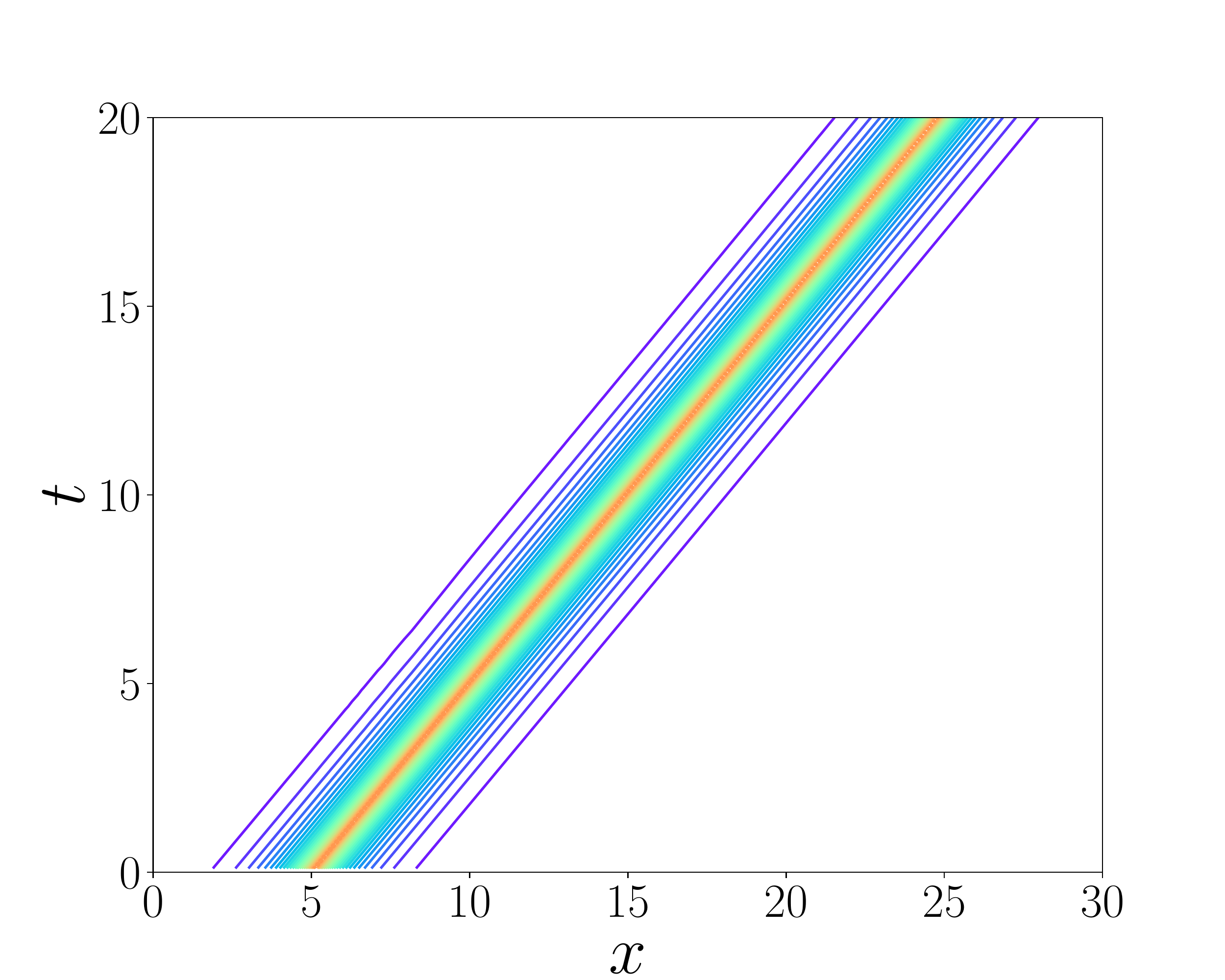}}
	\subfigure[Solution at $T = 20$.]{\includegraphics[width=0.3\textwidth]{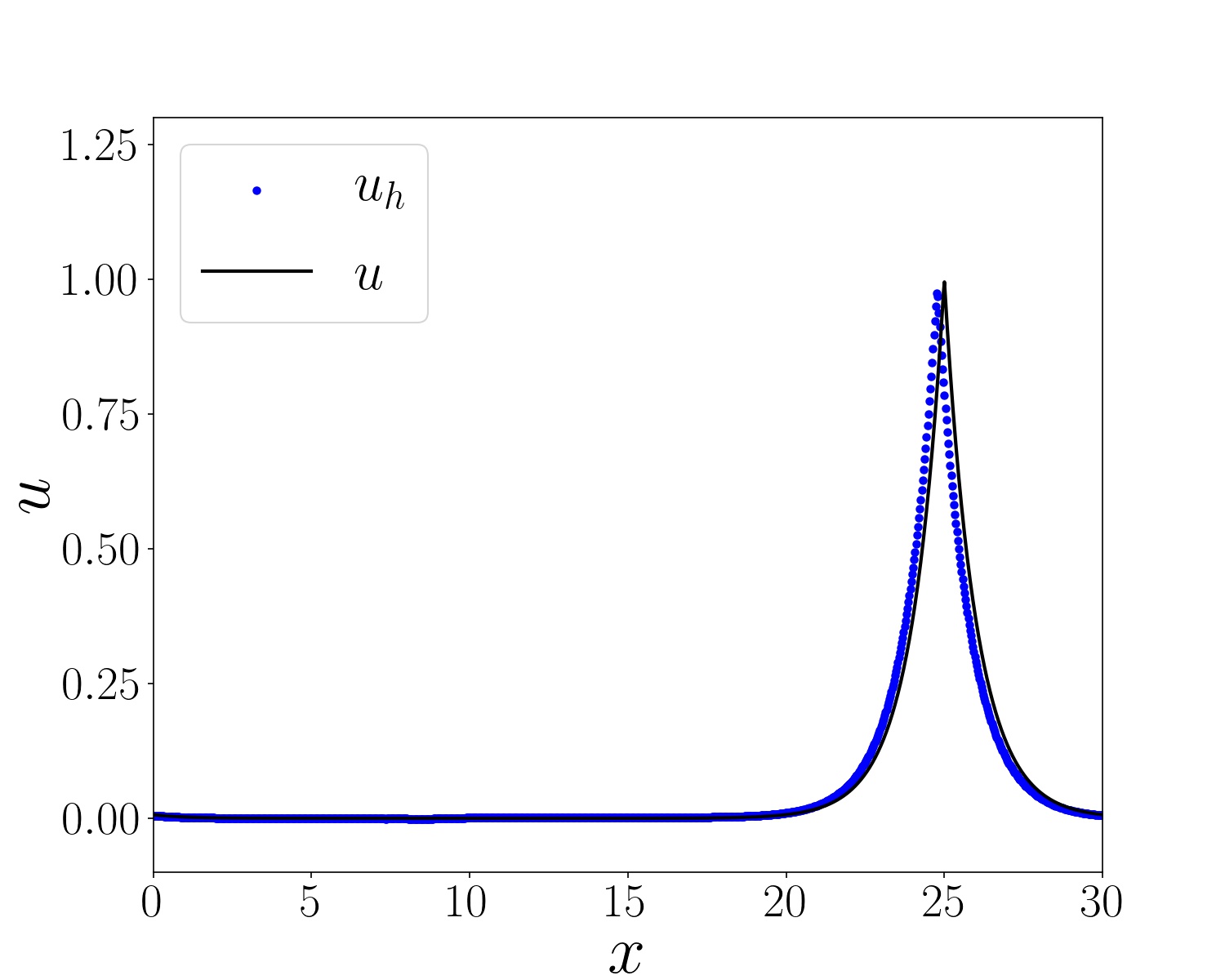}}
	\caption{Single traveling peakon of the CH equation. }\label{fig-CH-peakon}
\end{figure}
\begin{figure}[h!]
	\centering
	\subfigure[Solution profile.]{\includegraphics[width=0.38\textwidth]{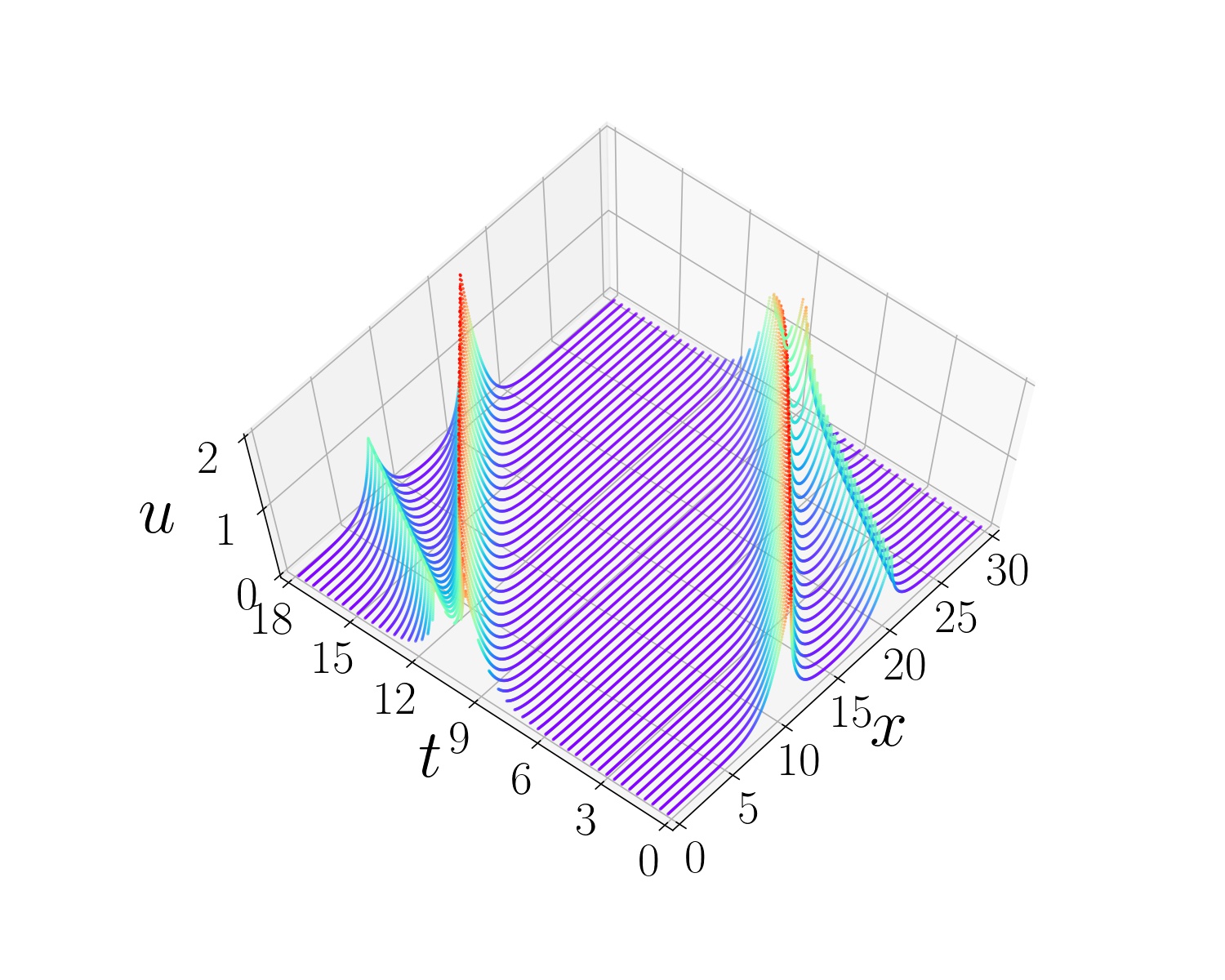}}
	\subfigure[Contour plot.]{\includegraphics[width=0.3\textwidth]{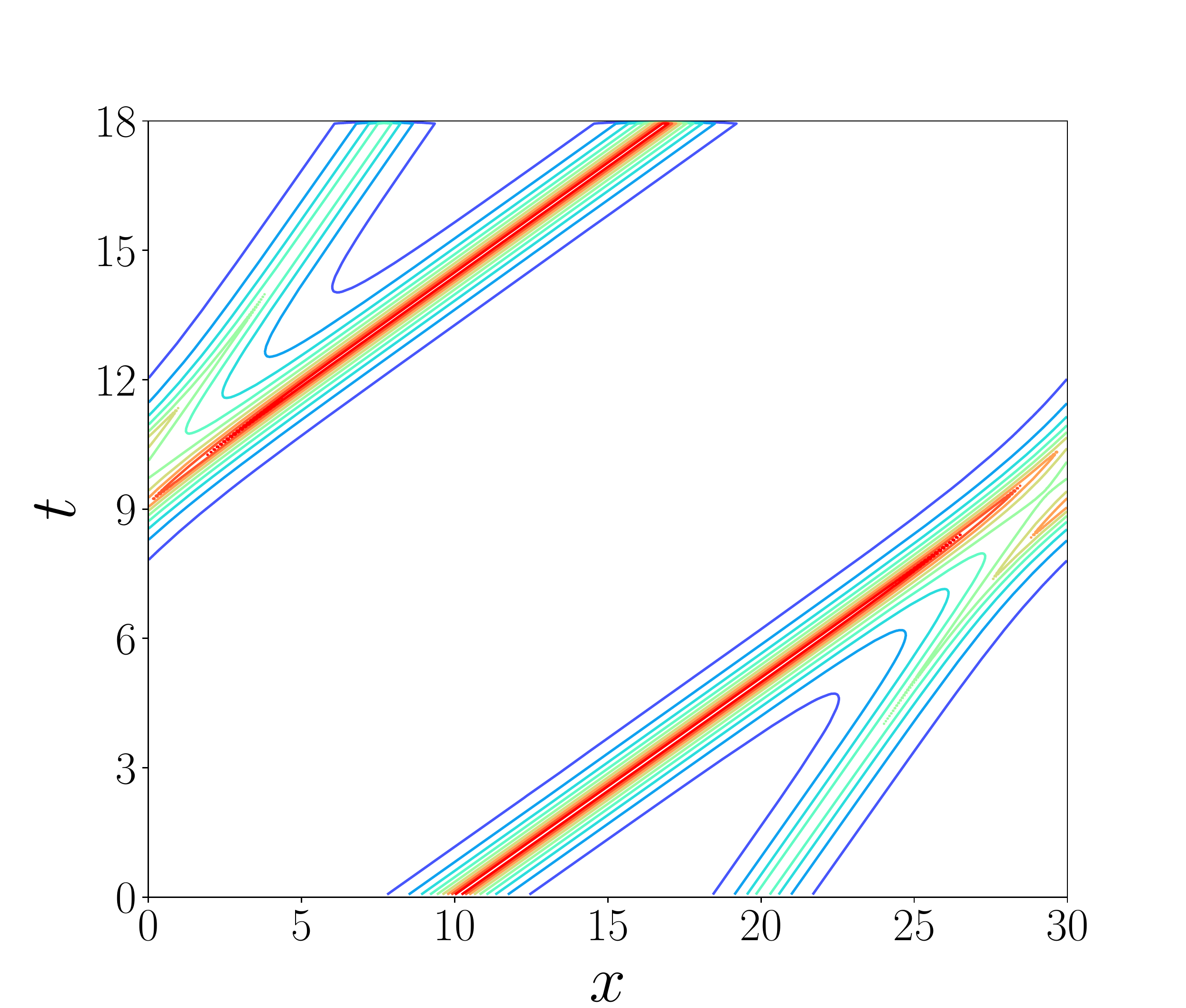}}
	\subfigure[Solution profile at $T = 18$]{\includegraphics[width=0.3\textwidth]{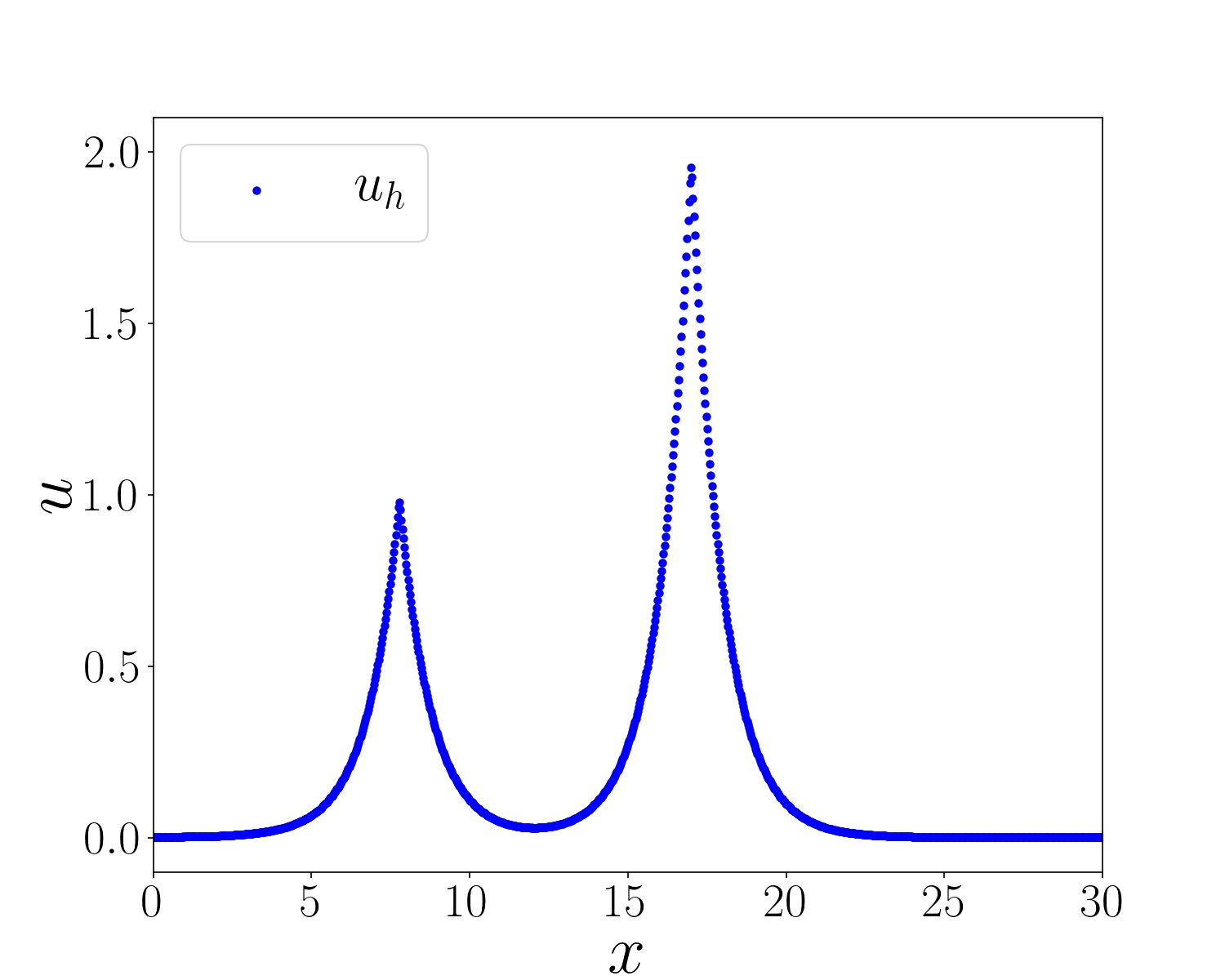}}
	\caption{Two-peakon interaction of the CH equation.}\label{fig-CH-2peakon}
\end{figure}
\begin{figure}[h!]
	\centering
	\subfigure[Solution profile.]{\includegraphics[width=0.38\textwidth]{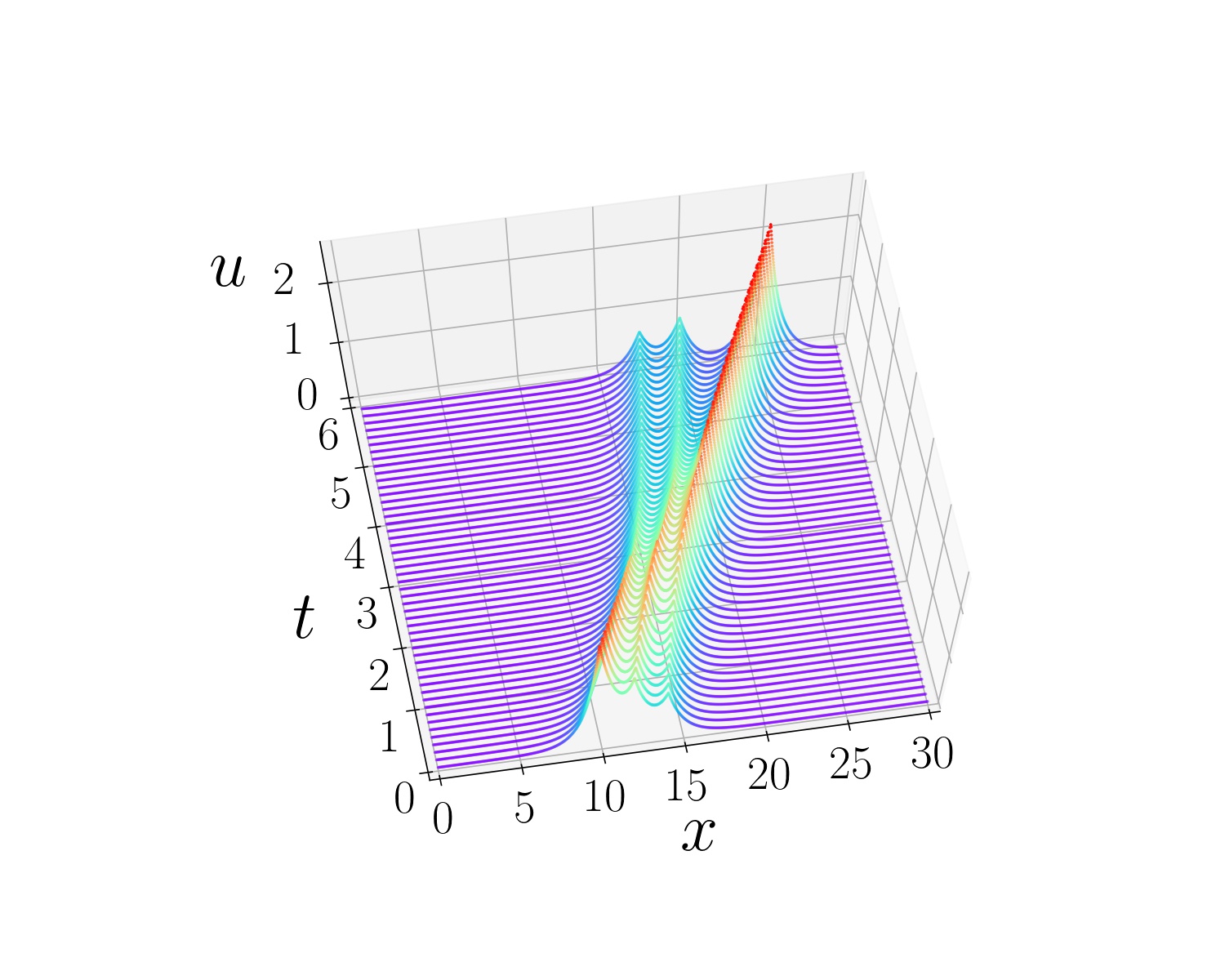}}
	\subfigure[Contour plot.]{\includegraphics[width=0.3\textwidth]{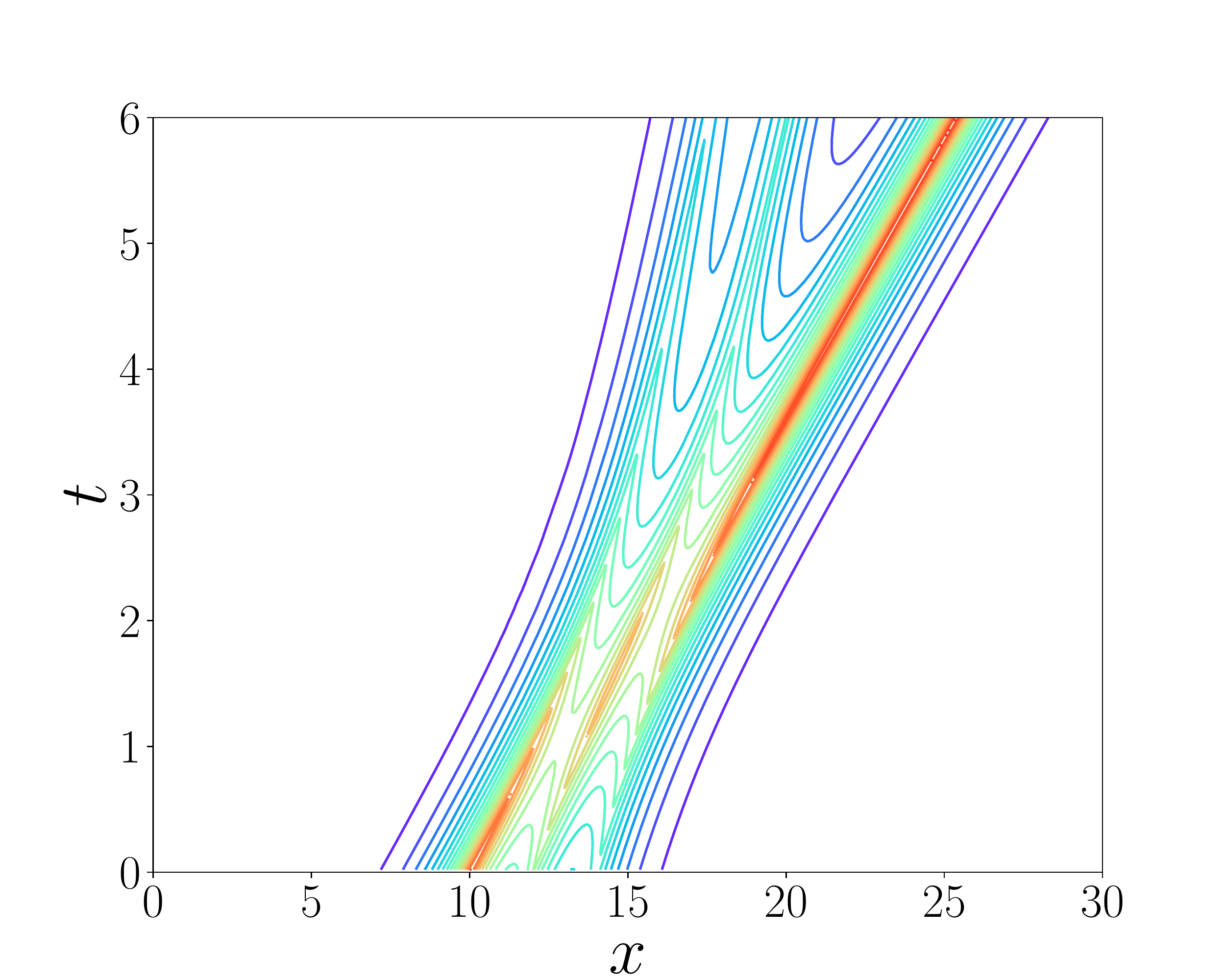}}
	\subfigure[Solution profile at $T = 6$.]{\includegraphics[width=0.3\textwidth]{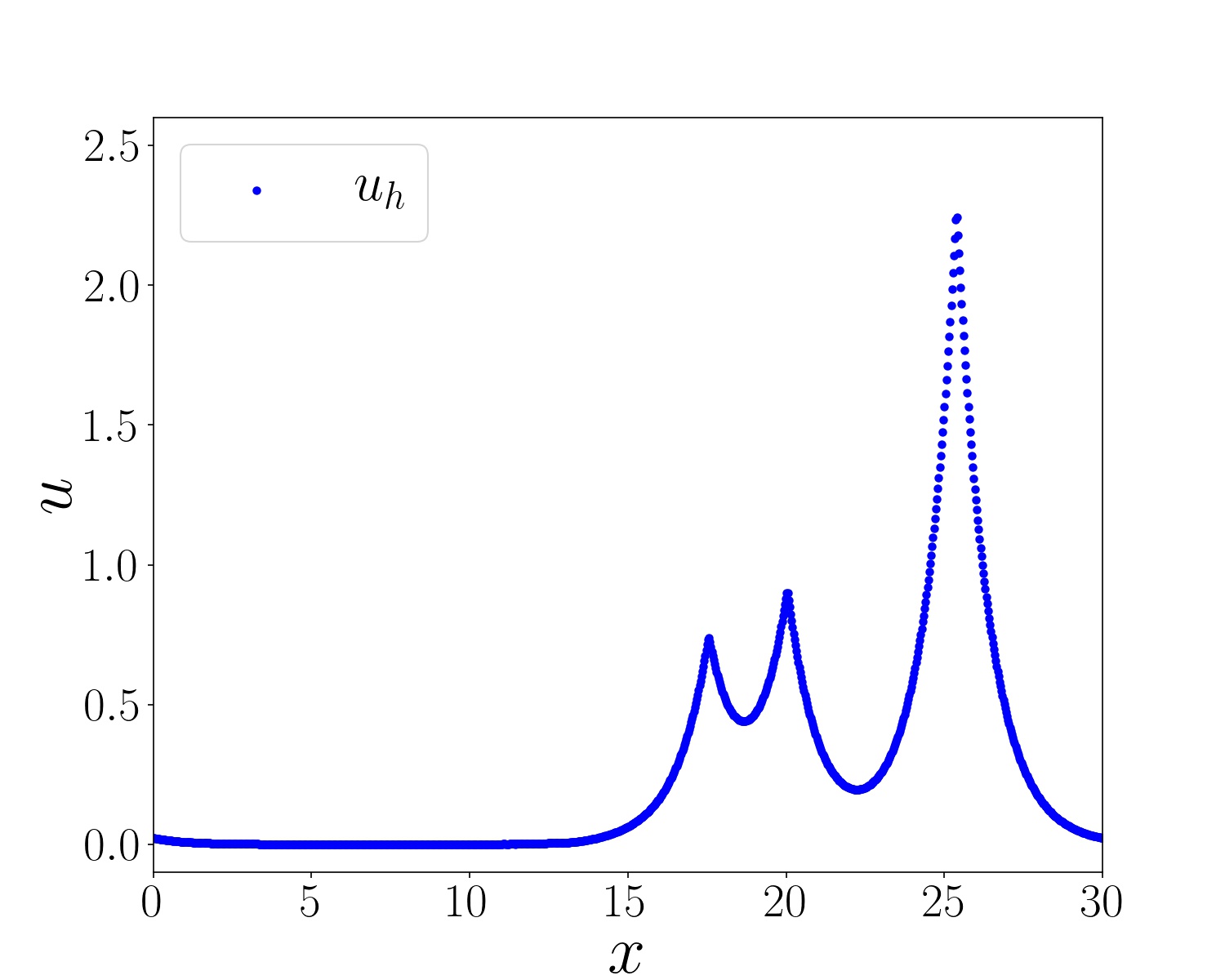}}
	\caption{Three-peakon interaction of the CH equation.}\label{fig-CH-3peakon}
\end{figure}
\begin{figure}[h!]
	\centering
	\subfigure[Solution profile.]{\includegraphics[width=0.38\textwidth]{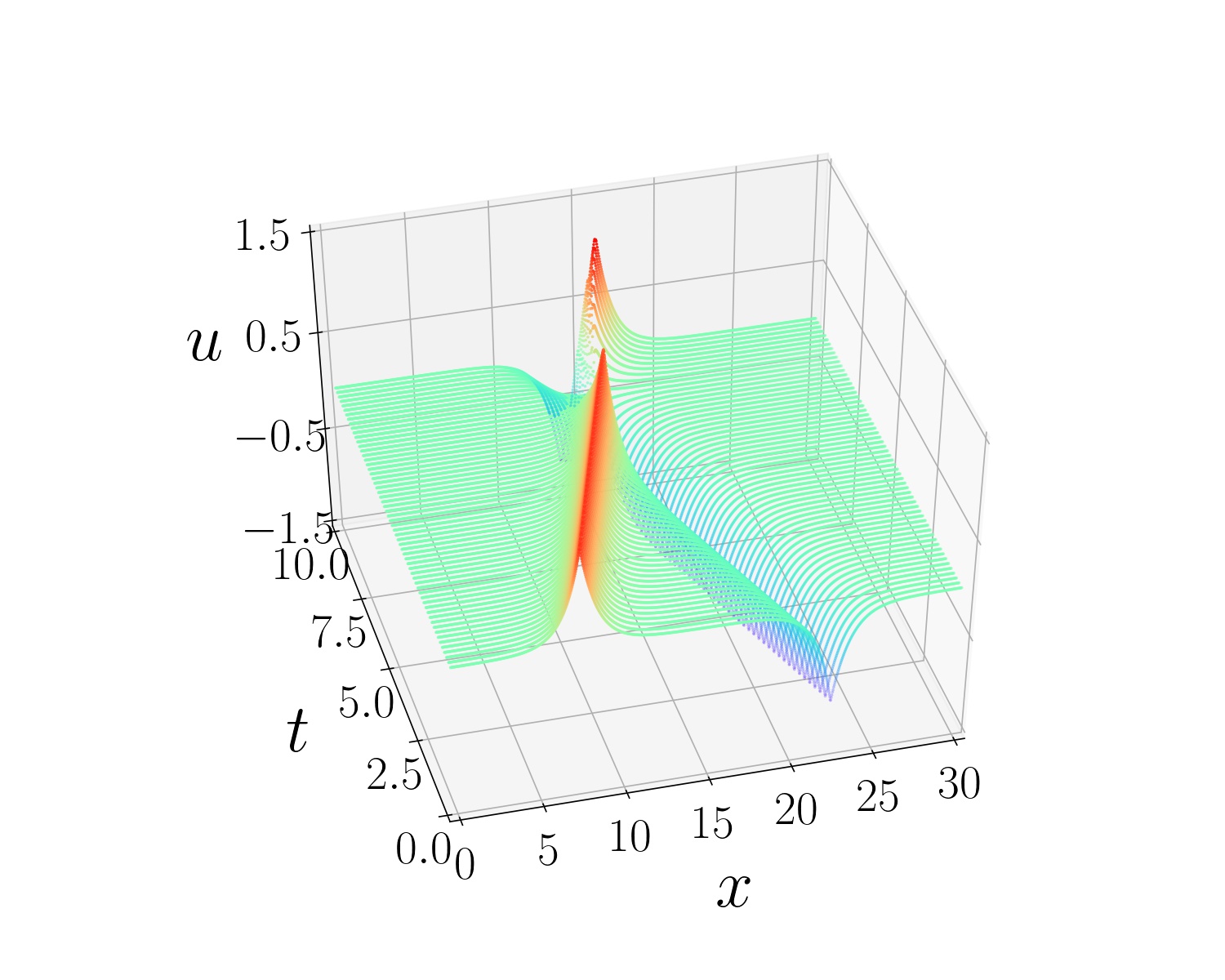}}
	\subfigure[Contour plot.]{\includegraphics[width=0.3\textwidth]{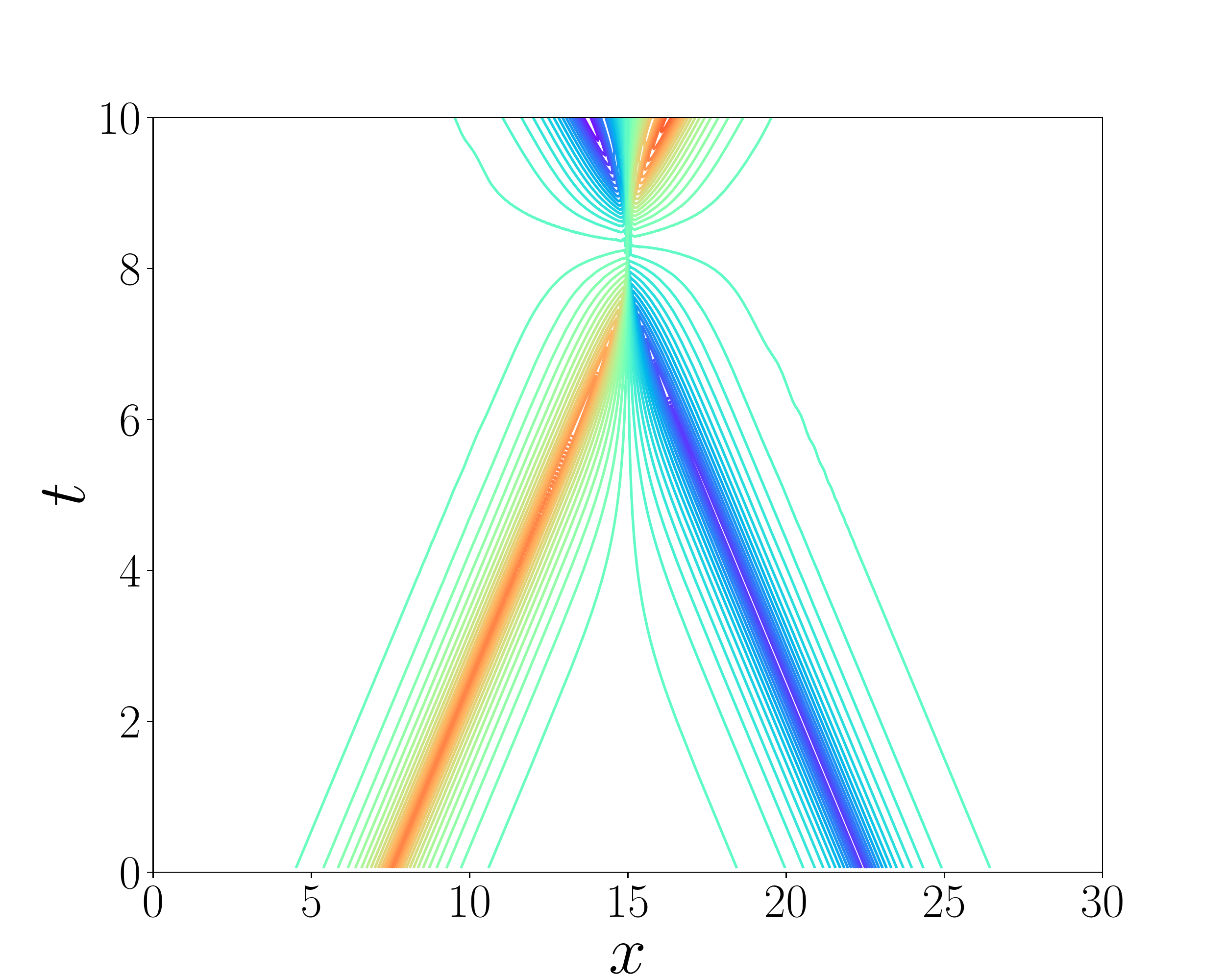}}
	\subfigure[Solution profile at $T = 10$.]{\includegraphics[width=0.3\textwidth]{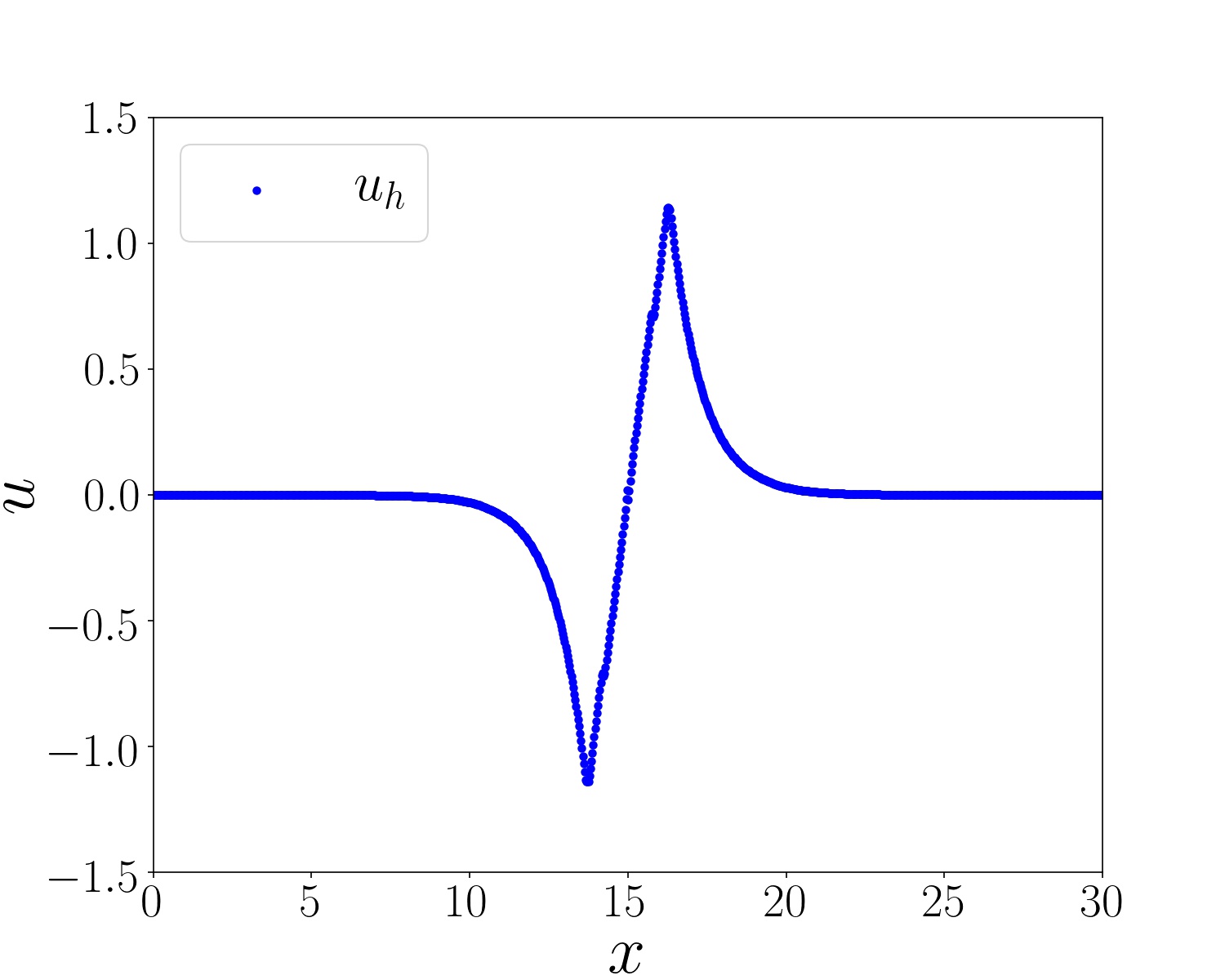}}
	\caption{Peakon-antipeakon interaction of the CH equation.}\label{fig-CH-antipeakon}
\end{figure}
\end{examp}

\section{Conclusions}\label{sec-concl}
\setcounter{equation}{0}

In this paper, we study the semi-discrete DG discretization of multi-symplectic systems, along with its application to various multi-symplectic HPDEs. With a general class of numerical fluxes, the semi-discrete DG schemes are shown to preserve the multi-symplectic structure and the local energy conservation simultaneously. For fully discrete schemes, one of the properties can also be achieved with suitable time integrators chosen. For applications, we particularly consider the wave equation, the BBM equation and the CH equation. The preserved energy functionals and their implementation methods are discussed. Numerically, we observe that different choices of numerical fluxes do have an influence on the accuracy of the schemes, and optimal convergence rate could be achieved with carefully chosen numerical fluxes. 
The DG scheme for each equation preserves its corresponding energy functional well in long time simulation with accurate time discretization, and enjoy the benefit of energy conserving methods like an improved long time behavior. Numerical simulations of multi-wave interactions are also provided to illustrate the performance of the scheme. A general class of numerical fluxes has been discussed in the paper, and we would like to provide the following remark as a guideline on the choice of numerical fluxes.
\begin{REM}[Choices of numerical fluxes]\label{rem-choices}
	Although in principle any numerical fluxes of the form \eqref{eq-flux} preserve the multi-symplectic structure and the energy conservation, in practice, we recommend fluxes that are easy to implement and retrieve optimal convergence rates. Central fluxes, with $A = B = 0$, usually lead to a scheme that is easy to implement. While it may suffer order degeneration for odd order polynomials or on nonuniform meshes. We recommend to use alternating fluxes if the implementation is not an issue, since the convergence rates are usually optimal for both $u_h$ and its discrete derivative. While alternating fluxes may not work as desired for the BBM and CH equations. In this situation, we suggest to tune $A$ such that the numerical flux for $u_h$ is no longer central, and then supplement extra fluxes to ensure symmetry of $A$. With this approach, the accuracy can be improved and the implementation is usually not that complicated. 
\end{REM}


\begin{thebibliography}{10}
	
	\bibitem{ascher2004multisymplectic}
	U.~M. Ascher and R.~I. McLachlan.
	\newblock Multisymplectic box schemes and the {K}orteweg--de {V}ries equation.
	\newblock {\em Applied Numerical Mathematics}, 48(3-4):255--269, 2004.
	
	\bibitem{bona2013conservative}
	J.~Bona, H.~Chen, O.~Karakashian, and Y.~Xing.
	\newblock Conservative, discontinuous {G}alerkin--methods for the generalized
	{K}orteweg--de {V}ries equation.
	\newblock {\em Mathematics of Computation}, 82(283):1401--1432, 2013.
	
	\bibitem{bridges2001multi}
	T.~J. Bridges and S.~Reich.
	\newblock Multi-symplectic integrators: numerical schemes for {H}amiltonian
	{P}{D}{E}s that conserve symplecticity.
	\newblock {\em Physics Letters A}, 284(4-5):184--193, 2001.
	
	\bibitem{bridges2001multiZK}
	T.~J. Bridges and S.~Reich.
	\newblock Multi-symplectic spectral discretizations for the
	{Z}akharov--{K}uznetsov and shallow water equations.
	\newblock {\em Physica D: Nonlinear Phenomena}, 152:491--504, 2001.
	
	\bibitem{cai2019two}
	J.~Cai and J.~Shen.
	\newblock Two classes of linearly implicit local energy-preserving approach for
	general multi-symplectic {H}amiltonian {P}{D}{E}s.
	\newblock {\em Journal of Computational Physics}, 401:108975, 2020.
	
	\bibitem{cai2018local}
	W.~Cai, Y.~Sun, Y.~Wang, and H.~Zhang.
	\newblock Local discontinuous {G}alerkin methods based on the multisymplectic
	formulation for two kinds of {H}amiltonian {PDEs}.
	\newblock {\em International Journal of Computer Mathematics}, 95(1):114--143,
	2018.
	
	\bibitem{celledoni2012preserving}
	E.~Celledoni, V.~Grimm, R.~I. McLachlan, D.~McLaren, D.~O’Neale, B.~Owren,
	and G.~Quispel.
	\newblock Preserving energy resp. dissipation in numerical {P}{D}{E}s using the
	“average vector field” method.
	\newblock {\em Journal of Computational Physics}, 231(20):6770--6789, 2012.
	
	\bibitem{chen2001multi}
	J.-B. Chen and M.-Z. Qin.
	\newblock Multi-symplectic {F}ourier pseudospectral method for the nonlinear
	{S}chr{\"o}dinger equation.
	\newblock {\em Electronic Transactions on Numerical Analysis}, 12:193--204,
	2001.
	
	\bibitem{chen2017entropy}
	T.~Chen and C.-W. Shu.
	\newblock Entropy stable high order discontinuous {G}alerkin methods with
	suitable quadrature rules for hyperbolic conservation laws.
	\newblock {\em Journal of Computational Physics}, 345:427--461, 2017.
	
	\bibitem{CSX2014}
	C.-S. Chou, C.-W. Shu, and Y.~Xing.
	\newblock Optimal energy conserving local discontinuous {G}alerkin methods for
	second-order wave equation in heterogeneous media.
	\newblock {\em Journal of Computational Physics}, 272:88--107, 2014.
	
	\bibitem{rkdg4}
	B.~Cockburn, S.~Hou, and C.-W. Shu.
	\newblock The {R}unge--{K}utta local projection discontinuous {G}alerkin finite
	element method for conservation laws. {I}{V}. the multidimensional case.
	\newblock {\em Mathematics of Computation}, 54(190):545--581, 1990.
	
	\bibitem{rkdg3}
	B.~Cockburn, S.-Y. Lin, and C.-W. Shu.
	\newblock {T}{V}{B} {R}unge--{K}utta local projection discontinuous {G}alerkin
	finite element method for conservation laws {I}{I}{I}: one-dimensional
	systems.
	\newblock {\em Journal of Computational Physics}, 84(1):90--113, 1989.
	
	\bibitem{rkdg2}
	B.~Cockburn and C.-W. Shu.
	\newblock {T}{V}{B} {R}unge--{K}utta local projection discontinuous {G}alerkin
	finite element method for conservation laws. {I}{I}. general framework.
	\newblock {\em Mathematics of Computation}, 52(186):411--435, 1989.
	
	\bibitem{rkdg1}
	B.~Cockburn and C.-W. Shu.
	\newblock The {R}unge--{K}utta local projection $
	{P}^1$-discontinuous-{G}alerkin finite element method for scalar conservation
	laws.
	\newblock {\em ESAIM: Mathematical Modelling and Numerical Analysis},
	25(3):337--361, 1991.
	
	\bibitem{CS1998SINUM}
	B.~Cockburn and C.-W. Shu.
	\newblock The local discontinuous {G}alerkin method for time-dependent
	convection-diffusion systems.
	\newblock {\em SIAM Journal on Numerical Analysis}, 35:2440--2463, 1998.
	
	\bibitem{rkdg5}
	B.~Cockburn and C.-W. Shu.
	\newblock The {R}unge--{K}utta discontinuous {G}alerkin method for conservation
	laws {V}: multidimensional systems.
	\newblock {\em Journal of Computational Physics}, 141(2):199--224, 1998.
	
	\bibitem{cohen2008multi}
	D.~Cohen, B.~Owren, and X.~Raynaud.
	\newblock Multi-symplectic integration of the {C}amassa--{H}olm equation.
	\newblock {\em Journal of Computational Physics}, 227(11):5492--5512, 2008.
	
	\bibitem{guermond2010asymptotic}
	J.-L. Guermond and G.~Kanschat.
	\newblock Asymptotic analysis of upwind discontinuous {G}alerkin approximation
	of the radiative transport equation in the diffusive limit.
	\newblock {\em SIAM Journal on Numerical Analysis}, 48(1):53--78, 2010.
	
	\bibitem{hairer2006geometric}
	E.~Hairer, C.~Lubich, and G.~Wanner.
	\newblock {\em Geometric numerical integration: structure-preserving algorithms
		for ordinary differential equations}, volume~31.
	\newblock Springer Science \& Business Media, 2006.
	
	\bibitem{hong2006multi}
	J.~Hong and C.~Li.
	\newblock Multi-symplectic {R}unge--{K}utta methods for nonlinear {D}irac
	equations.
	\newblock {\em Journal of Computational Physics}, 211(2):448--472, 2006.
	
	\bibitem{hong2006multiPRK}
	J.~Hong, H.~Liu, and G.~Sun.
	\newblock The multi-symplecticity of partitioned {R}unge-{K}utta methods for
	{H}amiltonian {P}{D}{E}s.
	\newblock {\em Mathematics of Computation}, 75(253):167--181, 2006.
	
	\bibitem{HLY2014}
	Y.~Huang, H.~Liu, and N.~Yi.
	\newblock A conservative discontinuous {G}alerkin method for the
	{Degasperis-Procesi} equation.
	\newblock {\em Methods and Applications of Analysis}, 21:67--90, 2014.
	
	\bibitem{leimkuhler2004simulating}
	B.~Leimkuhler and S.~Reich.
	\newblock {\em Simulating {H}amiltonian dynamics}, volume~14.
	\newblock Cambridge university press, 2004.
	
	\bibitem{li2013new}
	H.~Li and J.~Sun.
	\newblock A new multi-symplectic {E}uler box scheme for the {B}{B}{M} equation.
	\newblock {\em Mathematical and Computer Modelling}, 58(7-8):1489--1501, 2013.
	
	\bibitem{LSXC2020}
	X.~Li, W.~Sun, Y.~Xing, and C.-S. Chou.
	\newblock Energy conserving local discontinuous {G}alerkin methods for the
	improved boussinesq equation.
	\newblock {\em Journal of Computational Physics}, 401:109002, 2020.
	
	\bibitem{li2019optimal}
	X.~Li, Y.~Xing, and C.-S. Chou.
	\newblock Optimal energy conserving and energy dissipative local discontinuous
	{G}alerkin methods for the {B}enjamin-{B}ona-{M}ahony equation.
	\newblock {\em submitted}, 2019.
	
	\bibitem{LKX2015}
	X.~Liang, A.~Q.~M. Khaliq, and Y.~Xing.
	\newblock Fourth order exponential time differencing method with local
	discontinuous {G}alerkin approximation for coupled nonlinear {S}chrodinger
	equations.
	\newblock {\em Communications in Computational Physics}, 17:510--541, 2015.
	
	\bibitem{liu2016invariant}
	H.~Liu and Y.~Xing.
	\newblock An invariant preserving discontinuous {G}alerkin method for the
	{C}amassa--{H}olm equation.
	\newblock {\em SIAM Journal on Scientific Computing}, 38(4):A1919--A1934, 2016.
	
	\bibitem{liu2016hamiltonian}
	H.~Liu and N.~Yi.
	\newblock A {H}amiltonian preserving discontinuous {G}alerkin method for the
	generalized {K}orteweg--de {V}ries equation.
	\newblock {\em Journal of Computational Physics}, 321:776--796, 2016.
	
	\bibitem{mclachlan1999geometric}
	R.~I. McLachlan, G.~Quispel, and N.~Robidoux.
	\newblock Geometric integration using discrete gradients.
	\newblock {\em Philosophical Transactions of the Royal Society of London.
		Series A: Mathematical, Physical and Engineering Sciences},
	357(1754):1021--1045, 1999.
	
	\bibitem{mclachlan2015multisymplectic}
	R.~I. McLachlan and M.~Wilkins.
	\newblock The multisymplectic diamond scheme.
	\newblock {\em SIAM Journal on Scientific Computing}, 37(1):A369--A390, 2015.
	
	\bibitem{moore2003backward}
	B.~Moore and S.~Reich.
	\newblock Backward error analysis for multi-symplectic integration methods.
	\newblock {\em Numerische Mathematik}, 95(4):625--652, 2003.
	
	\bibitem{quispel2008new}
	G.~Quispel and D.~I. McLaren.
	\newblock A new class of energy-preserving numerical integration methods.
	\newblock {\em Journal of Physics A: Mathematical and Theoretical},
	41(4):045206, 2008.
	
	\bibitem{reed1973triangular}
	W.~H. Reed and T.~Hill.
	\newblock Triangular mesh methods for the neutron transport equation.
	\newblock Technical report, Los Alamos Scientific Lab., N. Mex.(USA), 1973.
	
	\bibitem{reich2000multi}
	S.~Reich.
	\newblock Multi-symplectic {R}unge--{K}utta collocation methods for
	{H}amiltonian wave equations.
	\newblock {\em Journal of Computational Physics}, 157(2):473--499, 2000.
	
	\bibitem{ryland2008multisymplecticity}
	B.~N. Ryland and R.~I. McLachlan.
	\newblock On multisymplecticity of partitioned {R}unge--{K}utta methods.
	\newblock {\em SIAM Journal on Scientific Computing}, 30(3):1318--1340, 2008.
	
	\bibitem{sun2004multi}
	Y.-J. Sun and M.-Z. Qin.
	\newblock A multi-symplectic scheme for {R}{L}{W} equation.
	\newblock {\em Journal of Computational Mathematics}, 22(4):611--621, 2004.
	
	\bibitem{sun2018discontinuous}
	Z.~Sun, J.~A. Carrillo, and C.-W. Shu.
	\newblock A discontinuous {G}alerkin method for nonlinear parabolic equations
	and gradient flow problems with interaction potentials.
	\newblock {\em Journal of Computational Physics}, 352:76--104, 2018.
	
	\bibitem{sun2019enforcing}
	Z.~Sun and C.-W. Shu.
	\newblock Enforcing strong stability of explicit {R}unge--{K}utta methods with
	superviscosity.
	\newblock {\em submitted}, 2019. arXiv:
	\texttt{\href{https://arxiv.org/abs/1912.11596}{1912.11596 [math.NA]}}.
	
	\bibitem{tang2017discontinuous}
	W.~Tang, Y.~Sun, and W.~Cai.
	\newblock Discontinuous {G}alerkin methods for hamiltonian {O}{D}{E}s and
	{P}{D}{E}s.
	\newblock {\em Journal of Computational Physics}, 330:340--364, 2017.
	
	\bibitem{xing2013energy}
	Y.~Xing, C.-S. Chou, and C.-W. Shu.
	\newblock Energy conserving local discontinuous {G}alerkin methods for wave
	propagation problems.
	\newblock {\em Inverse Problems \& Imaging}, 7(3), 2013.
	
	\bibitem{xing2006high}
	Y.~Xing and C.-W. Shu.
	\newblock High order well-balanced finite volume {W}{E}{N}{O} schemes and
	discontinuous {G}alerkin methods for a class of hyperbolic systems with
	source terms.
	\newblock {\em Journal of Computational Physics}, 214(2):567--598, 2006.
	
	\bibitem{XZS2010}
	Y.~Xing, X.~Zhang, and C.-W. Shu.
	\newblock Positivity-preserving high order well-balanced discontinuous
	{G}alerkin methods for the shallow water equations.
	\newblock {\em Advances in Water Resources}, 33:1476--1493, 2010.
	
	\bibitem{xushuNLS}
	Y.~Xu and C.-W. Shu.
	\newblock Local discontinuous {G}alerkin methods for nonlinear
	{S}chr\"{o}dinger equations.
	\newblock {\em Journal of Computational Physics}, 205(1):72--97, 2005.
	
	\bibitem{zhang2010maximum}
	X.~Zhang and C.-W. Shu.
	\newblock On maximum-principle-satisfying high order schemes for scalar
	conservation laws.
	\newblock {\em Journal of Computational Physics}, 229(9):3091--3120, 2010.
	
	\bibitem{zhao2000multisymplectic}
	P.-F. Zhao and M.-Z. Qin.
	\newblock Multisymplectic geometry and multisymplectic {P}reissmann scheme for
	the {K}d{V} equation.
	\newblock {\em Journal of Physics A: Mathematical and General}, 33(18):3613,
	2000.
	
\end{thebibliography}
\end{document}